\documentclass[11pt,reqno]{amsart}

\usepackage{amsfonts}
\usepackage{amsmath}
\usepackage{amssymb, esint}
\usepackage{amsthm,color}
\usepackage{enumerate}
\usepackage{mathtools}
\usepackage{enumitem}
\usepackage{url}
\usepackage[hidelinks, bookmarksdepth=3]{hyperref}
\usepackage{xcolor}
\usepackage{cancel}
\usepackage{ulem}
\usepackage{caption,subcaption}
\usepackage{tikz}
\usepackage{tikz-3dplot}

\usepackage{longtable}

\usepackage{soul}

\tdplotsetmaincoords{70}{110}

\usepackage{geometry}
\def\haus{\,d{\mathcal H}^{n-1}}

\newcommand{\R}{{\mathbb R}}
\newcommand{\C}{{\mathbb C}}
\newcommand{\Z}{{\mathbb Z}}
\newcommand{\N}{{\mathbb N}}

\newcommand{\D}{{\mathbb D}}
\newcommand{\cE}{{\mathcal E}}
\newcommand{\cL}{{\mathcal L}}

\newcommand{\cN}{{\mathcal N}}

\newcommand{\cA}{{\mathcal A}}
\newcommand{\cB}{{\mathcal B}}
\newcommand{\cC}{{\mathcal C}}
\newcommand{\cD}{{\mathcal D}}

\newcommand{\cH}{{\mathcal H}}

\newcommand{\cF}{{\mathcal F}}
\newcommand{\Ss}{{\mathbb S}}

\def\0{{\mathbf 0}}
\def\loc{{\textup{loc}}}
\def\id{{\rm id}}
\def\Id{{\rm Id}}
\def \mb{\mathbb}

\newcommand{\e}{\varepsilon}
\newcommand{\vp}{\varphi}
\newcommand{\osc}{\operatornamewithlimits{osc}}

\newcommand{\supp}{\operatorname{supp}}
\newcommand{\ddiv}{\operatorname{div}}
\newcommand{\diam}{\operatorname{diam}}
\newcommand{\dist}{\operatorname{dist}}

\newcommand{\Hess}{\operatorname{Hess}}
\newcommand{\re}{\operatorname{Re}}

\newcommand{\esssup}{\operatornamewithlimits{ess\,sup}}

\newcommand{\sgn}{\operatorname{sgn}}

\theoremstyle{plain}
\newtheorem{thm}{Theorem}[section]

\newtheorem{cor}[thm]{Corollary}
\newtheorem{lem}[thm]{Lemma}
\newtheorem{prop}[thm]{Proposition}

\newcommand{\thistheoremnames}{}
\newtheorem*{genericthms}{\thistheoremnames}
\newenvironment{para*}[1]
  {\renewcommand{\thistheoremnames}{#1}%
   \begin{genericthms}}
  {\end{genericthms}}

\newtheorem{defn}[thm]{Definition}

\theoremstyle{remark}

\newtheorem{ex}[thm]{Example}

\newtheorem{case}{Case}
\newtheorem{step}{Step}

\newtheorem*{claim*}{Claim}
\newtheorem{notation}[thm]{Notation}
\newtheorem{rem}[thm]{Remark}

\numberwithin{equation}{section}

\title[Constraint maps: singularities vs free boundaries]{
Constraint maps: {singularities vs free boundaries}
}

\author[A.\ Figalli]{Alessio Figalli}
\email{alessio.figalli@math.ethz.ch}
\address{Department of Mathematics, ETH Z\"urich,  Raemistrasse 101, 8092 Z\"urich, Switzerland }

\author[A.\ Guerra]{Andr\'e Guerra}
\email{andre.guerra@eth-its.ethz.ch}
\address{Institute for Theoretical Studies, ETH Z\"urich,  CLV, Clausiusstrasse 47, 8006 Z\"urich, Switzerland}

\author[S.\ Kim]{Sunghan Kim}
\email{sunghan.kim@math.uu.se}
\address{Department of Mathematics, Uppsala University, S-751 06 Uppsala, Sweden}

\author[H.\ Shahgholian]{Henrik Shahgholian}
\email{henriksh@kth.se}
\address{Department of Mathematics, KTH Royal Institute of Technology, 100 44 Stockholm, Sweden}

\begin{document}

\begin{abstract}
Energy-minimizing constraint maps are a natural extension of the obstacle problem within a vectorial framework. Due to inherent topological constraints, these maps manifest a diverse structure that includes singularities similar to harmonic maps, branch points reminiscent of minimal surfaces, and the intricate free-boundary behavior of the obstacle problem. 
The complexity of these maps poses significant challenges to their analysis.

In this paper, we first focus on constraint maps with uniformly convex obstacles and establish continuity (and therefore higher-order regularity) within a uniform neighborhood of the free boundary. More precisely, thanks to a new quantitative unique continuation principle near singularities (which is new even in the setting of classical harmonic maps), we prove that, in the uniformly convex setting, topological singularities can only lie in the interior of the contact set. We then establish the optimality of this result.

Second, while exploring the structure of the free boundary, we investigate the presence of branch points and show how they lead to completely new types of singularities not present in the scalar case.

\end{abstract}

\maketitle
\setcounter{tocdepth}{1}

\tableofcontents


\section{Introduction}\label{sec:intro}

Let $M$ be a smooth domain with compact, connected complement in $\R^m$, $m \geq 2$, and let $\Omega$ be a bounded domain in $\R^n$, with $n\geq 1$. In this paper we consider maps $u\in W^{1,2}(\Omega;\overline M)$ minimizing the Dirichlet energy
\begin{equation}\label{eq:main}
\cE[v] := \int_\Omega |Dv|^2\,dx,
\end{equation}
among all admissible maps $v\in W^{1,2}(\Omega;\overline M)$ with $v - u \in W_0^{1,2}(\Omega;\R^m)$. To emphasize the target constraint $u \in \overline M$ as well as its energy minimality, we call such maps {\it minimizing constraint maps}. As shown in \cite{D}, minimizing constraint maps are weak solutions to the Euler--Lagrange system
\begin{equation}\label{eq:main-sys}
\Delta u = A_u(Du,Du)\chi_{u^{-1}(\partial M)}\quad\text{in }\Omega, 
\end{equation}
where $A$ is the second fundamental form of the boundary $\partial M$ relative to $M$.  

We shall call $M^c := \R^m\setminus M$ the {\it obstacle}, as the image of the constraint maps is not allowed to go inside it. The presence of the obstacle naturally induces the {\it free boundary}\footnote{The notion of free boundary here differs significantly from the one in the theories of harmonic maps and minimal surfaces, as our notion resides in the domain, while the latter appears in the target.}  $\Omega\cap \partial u^{-1}(M)$, which is the interface between the {\it non-coincidence set} $u^{-1}(M)$ and the {\it coincidence set} $u^{-1}(\partial M)$. System \eqref{eq:main-sys} shows that, when restricted to ${\rm int}(u^{-1}(\partial M))\subset \Omega$, constraint maps are harmonic maps into $\partial M$. At the same time, each component of a constraint map is harmonic in ${\rm int}(u^{-1}(M))$. This variational problem hence can be seen as a canonical extension of the scalar obstacle problem into the vectorial setting. 

Being vectorial maps, constraint maps develop in general discontinuities\footnote{Although the term ``singularities'' is more common in the Harmonic Maps literature, here we have chosen to use the word ``discontinuities'' to ensure that the Free Boundary community does not confuse them with singularities of the free boundaries.} 
due to topological differences between the domain and their image in the target.  Their \textit{singular set} is denoted   throughout the paper by 
\begin{equation}\label{eq:sing}
\begin{aligned}
\Sigma(u) := \{ x\in\Omega: u\text{ is not continuous at }x\}.
\end{aligned}
\end{equation}
Besides discontinuities, constraint maps can also exhibit branching, which is concerned with the regularity of the image of the map and is encapsulated by the \textit{branch set}
\begin{equation}
    \label{eq:branching}
    \cB(u) := \bigcup_{k=0}^{\min\{n,m\}-1} \cB_k(u),\qquad \cB_k(u) := \{x\in \Omega\setminus \Sigma(u): \mathrm{rank}(Du(x))=k \}.
\end{equation}
Here, the word branching is motivated by the connection with minimal surfaces and geometric function theory,  e.g.\  the map  $u\colon z\mapsto (z^2, \re z^{\ell/2})$, defined for $z\in \C=\R^2$, is such that $\cB(u)= \cB_0(u)= \{0\}$. Points in $\cB(u)$ can either be \textit{true branch points},
such as when $\ell$ is odd, or \textit{false branch points}, 
such as when $\ell$ is even.
More precisely, false branch points correspond to a degeneracy of the parametrization and not of the actual image that is being parametrized, while true branch points correspond to singularities of the image; see Figure \ref{fig:branches} for an illustration. We will be especially concerned with the case $k=0$ since the right-hand side in \eqref{eq:main-sys} vanishes at points in $\cB_0(u)$, making the obstacle problem degenerate.

\medskip

Interestingly, constraint maps reveal an extremely rare interplay between $\partial u^{-1}(M)$, $\Sigma(u)$, and $\cB(u)$ altogether. Despite the long history of constraint maps since their initial appearance in the literature, dating back to the 1970s by the work of Hilderbrandt \cite{Hildebrandt1977, Hildebrandt1977a} and Tomi \cite{T1, T2}, the interaction between these apparently correlated objects was addressed only recently by the authors \cite{FGKS} in a different but related setting. 
The goal of this paper is to address the following fundamental question: 
\begin{center}
    \textit{
    Does the free boundary meet the mapping singularities manifested by $\Sigma(u)$ or $\cB(u)$?
    }
\end{center}
This is an inherently difficult question,  
as answering it requires an understanding of the interplay between local analytic aspects and global geometric and topological  constraints.

One of the major findings of this paper is that there are no discontinuities of the maps near their free boundaries, provided that the obstacle is {\it uniformly convex}\footnote{Uniform convexity of $M^c$ means that there exists $R>0$ such that the following holds: for any point $x \in \partial M$ there exists a ball $B_R(y)\supset M^c$ such that $x \in \partial B_R(y)$.}.
This condition is essentially optimal: we find a large class of general convex obstacles which produce discontinuities on the free boundaries. 
To our knowledge, and despite the significance of this problem, the first suggestion of a possible resolution of discontinuities by free boundaries was raised in the authors’ recent work \cite{FGKS} in a different but related setting. 

The second major finding of this paper is that, even when the map is regular and the obstacle is a ball, the map may exhibit branching
on the free boundary in such a way that the free boundary itself becomes singular. Indeed, it is  important to observe that the right-hand side in \eqref{eq:main-sys} becomes zero at all points within $\cB_0(u)$. Consequently, at these points, the behavior of the free boundary may exhibit significant irregularities, which do not occur in the classical scalar obstacle problem.

As we shall explain later, once discontinuities and branch points are understood,  the free boundary regularity can be reduced to that for the scalar obstacle problem (cf.\ \cite{FKS, F1}).

\subsection{A brief overview of the literature}

The study of minimizing constraint maps as solutions to vectorial obstacle problems dates back to the 1970s. Pioneering work by Hildebrandt \cite{H1, H2}, and Tomi \cite{T1, T2} independently addressed the regularity problems for constraint maps of two independent variables, employing different methodologies.
Early advances in higher dimensions can be attributed to the work of Hildebrandt-Widman \cite{HW}, who considered convex target constraints. The first comprehensive result in higher dimensions was established by Duzaar \cite{D}, who proved an $\varepsilon$-regularity theorem. As a consequence, a partial regularity theorem was derived, affirming the $W^{2,p}$-regularity of mappings outside a closed singular set with dimension less than $n-2$. Subsequently, Duzaar--Fuchs \cite{DF} proved that the singular set has dimension at most $n-3$. Although this result is optimal in general, it was later observed in \cite{F1,FW} that, under certain topological and geometric conditions on the target constraint, discontinuities are not present.
We note that the development of these works closely followed the higher-dimensional theory of minimizing harmonic maps, pioneered by Schoen--Uhlenbeck~\cite{ScU1,ScU2}: as for harmonic maps, 
the theory of constraint maps has been generalized to the setting of the $p$-Dirichlet energy in \cite{F4,F3}, almost minimizing maps in \cite{L},  and weakly constraint maps as well as the associated heat flow problem in \cite{CM}. The regularity results in two space dimensions were recently extended to weakly constraint maps in \cite{KS}. We also refer the reader to \cite{DMM} for results in the case where the obstacle is \textit{thin}, in the sense that it has higher codimension.

It is worth noting that, in the simple case where the domain is one-dimensional, minimizing constraint maps are simply energy-minimizing geodesics on the manifold-with-boundary $\overline M$. The regularity and (non-)uniqueness of geodesics on manifolds-with-boundary were studied in depth in \cite{AB, ABB2,ABB1} and also more recently in \cite{LL}. 

Quite surprisingly, and despite all of the above results, the theory did not progress much since the 1990s, until recently three of the authors revived the problem in \cite{FKS} from the perspective of free boundary problems, opening up new and delicate issues on the exact regularity of the maps away from their singular set and around their free boundaries. Although the study of minimizing constraint maps had its origin as a vectorial obstacle problem, historically only basic properties of the free boundaries were shown in \cite{F2}. The result closest to this paper is, in fact, the very recent work \cite{FGKS} of ours concerning constraint maps minimizing a vectorial Alt--Caffarelli energy,  which corresponds to a vectorial version of the Bernoulli problem.

In recent years there has been a surge of interest from the free boundary community towards extending scalar problems into vectorial settings, see for instance \cite{CSY,DPESV,KL1, KL2,MTV1, MTV2} for Bernoulli-type problems from diverse angles,  \cite{ASUW,SY1, SY2, SY3} for obstacle-type problems,  and \cite{And,DST} for thin-obstacle-type problems. All the aforementioned works successfully extend aspects of the scalar theory and explore inherently interesting issues related to the branching of free boundaries from several components. 

However, our minimizing constraint maps exhibit very different characteristics compared to the solutions to the aforementioned vectorial problems, due to the possible development of \textit{discontinuities of the map itself} around the free boundaries. Hence a complete understanding of constraint maps requires bridging the theory of harmonic maps with that of free boundary problems.

\subsection{Main results}

Our first result asserts that, whenever the obstacle is uniformly convex, minimizing constraint maps are uniformly $C^{1,1}$-regular in a universal neighborhood of the non-coincidence set (hence, in a set which strictly includes the free boundary). 
It should be noted that $C^{1,1}$-regularity is {optimal}, due to the classical theory of the obstacle problem. 
In addition, the result cannot be improved: solutions cannot be continuous everywhere neither in the interior of the coincidence set,  due to the possible presence of singularities coming from minimizing harmonic maps with values in $\partial M$.

\begin{thm}\label{thm:uni-reg}
Let $M\subset \R^m$ be a smooth domain with uniformly convex complement. There exist constants $\delta \in (0,\frac{1}{2})$ and $c > 1$, depending only on $n$, $m$, $\partial M$, and $\Lambda$, such that for any minimizing constraint map $u\in W^{1,2}(B_2;\overline M)$ satisfying $\| u \|_{W^{1,2}(B_2)} \leq \Lambda$, we have $u \in C^{1,1}(B_1\cap B_\delta(u^{-1}(M)))$ with the estimate 
$$ 
\| u \|_{C^{1,1}(B_1\cap B_\delta(u^{-1}(M)))} \leq c.
$$
\end{thm}

To describe the strategy and challenges in the proof of Theorem \ref{thm:uni-reg}, it is fruitful to introduce the \textit{distance function} to $\partial M$, which we denote by $\rho$.
Note that the smoothness and convexity of $M^c$ guarantee that 
$\rho$ is smooth in $\overline M$. Then \eqref{eq:main-sys} 
yields an Euler--Lagrange equation for $\rho\circ u$, namely
\begin{align}
    \label{eq:ELrho}
\Delta(\rho\circ u ) & = \Hess \rho_u(Du,Du) \chi_{u^{-1}(M)},
\end{align}
where $\Hess \rho_u$ denotes the Hessian of $\rho$ evaluated at $u$. Noting that $u^{-1}(M)$ is precisely the positivity set of $\rho\circ u$, we observe that \eqref{eq:ELrho} is a \textit{scalar} obstacle-type equation, but with an important caveat. Indeed, the nonlinearity $\Hess \rho_u(Du,Du)$, which acts as a force in \eqref{eq:ELrho}, 
can become unbounded (when $u$ is discontinuous) and can also vanish (when either $\Hess\rho_u$, $Du$, or their bilinear combination vanishes). Hence $\rho \circ u$ can behave quite differently from a solution to the scalar obstacle problem, where one deals with forcing terms that are regular and bounded away from zero. In particular, $\rho \circ u$ has no inherent non-degeneracy: as we show in Theorem \ref{thm:axial} below, it may in fact vanish to arbitrarily high order on the free boundary.

In our previous work \cite{FGKS} (in a slightly different context) we observed that the distance map $\rho\circ u$ is continuous everywhere, even across $\Sigma(u)$, and that the non-coincidence set $u^{-1}(M)$ has \textit{null density} at any point of $\Sigma(u)$: at any (small) scale $r$, we have 
\begin{equation}\label{eq:nulldensity}
|B_r(x_0)\cap u^{-1}(M)| = o(r^n),\quad\forall x_0\in\Sigma(u). 
\end{equation}
Interestingly, the assertion \eqref{eq:nulldensity} is proved via the unique continuation principle (UCP) for constraint/harmonic maps. However, to show that $\Sigma(u)$ is completely disjoint from the free boundary, we need the much stronger assertion that $\rho\circ u$  vanishes identically in a full neighborhood of $\Sigma(u)$, i.e.\ there is some small $r_0>0$ such that 
\begin{equation}
    \label{eq:inc}
    B_{r_0}(x_0)\subset u^{-1}(\partial M),\quad \forall x_0\in\Sigma(u).
\end{equation} 
Our idea is then as follows. We look at the distance map $\rho\circ u$, which is supported exactly on $u^{-1}(M)$, and is subharmonic for convex obstacles $M^c$. A key property of $\rho\circ u$ is that it vanishes continuously at any point of $\Sigma(u)$ whence, upon translating and rescaling, we may assume that $0\leq \rho\circ u \leq \e$ in $B_2$ and $0\in\Sigma(u)$, for a given $\e > 0$. Then, by a lemma reminiscent of the De Giorgi argument (Lemma \ref{lem:degiorgi}), we observe that $\rho\circ u = 0$ in $B_1$, provided that $\e$ is sufficiently small and  
\begin{equation}
    \label{eq:holderdecay}
    |B_r\cap u^{-1}(M)|\leq \frac{c}{(R-r)^\beta} \|\rho\circ u\|_{L^\infty(B_R)}^\alpha,\quad \forall\, 1 \leq r< R \leq 2,
\end{equation}
for some $\alpha$, $\beta$, and $c$ that do not depend on the scales $r$ and $R$. Condition \eqref{eq:holderdecay} can be considered as the reverse of the local maximum principle for subharmonic functions. Since we are free to choose the smallness parameter $\e$ for the sup-norm of $\rho\circ u$ over $B_2$, our task boils down to verifying \eqref{eq:holderdecay}. 

What appears to be intriguing is the involvement in \eqref{eq:holderdecay} of the vectorial character of $u$ implied by the presence of $0\in\Sigma(u)$ (note that $u^{-1}(M)$ is simply the set where $\rho\circ u$ is positive). To exploit the hidden vectorial character in \eqref{eq:holderdecay}, we resort to the Euler-Lagrange equation \eqref{eq:ELrho} for $\rho\circ u$. Integrating both sides of \eqref{eq:ELrho}, and invoking the uniform convexity of $M^c$, we obtain around $0\in\Sigma(u)$,\footnote{That $0\in\Sigma(u)$ is needed to get $\Hess\rho_u(Du,Du)\approx |Du|^2$ in the $L^1$-sense. In general, we always have $\Hess\rho_u(Du,Du)\approx |D(\Pi\circ u)|^2$ a.e., with $\Pi$ being the nearest point projection to $\partial M$ (by \eqref{eq:xitop} and \eqref{eq:hessrho}). However, $\rho\circ u$ becomes negligible around any discontinuity point of $u$, so we obtain the last inequality in \eqref{eq:reverse1}.} along a chain of elementary tools (see \eqref{eq:key11} -- \eqref{eq:key1}), that 
\begin{equation}\label{eq:reverse1}
\frac{1}{(R-r)^2} \int_{B_R} (\rho\circ u)\,dx \gtrsim \int_{B_r\cap u^{-1}(M)} \Hess\rho_u(Du,Du) \,dx \gtrsim  \int_{B_r\cap u^{-1}(M)} |Du|^2 \,dx
\end{equation}
for every $1\leq r<R\leq 2$. 

Our analysis would have been softer if $|Du| \gtrsim 1$ uniformly in $B_r$, as the rightmost side of \eqref{eq:reverse1} would then readily yield \eqref{eq:holderdecay} with $w = \rho\circ u$. Since $|Du|$ blows up at $0\in\Sigma(u)$ possibly with a rate $1/r$ in $B_r$, $|Du|\gtrsim 1$ is expected to be true at least for a large portion of the neighborhood of $0\in\Sigma(u)$.
However, $|Du|$ may not behave so uniformly near its singularity. In principle, $|Du|$ may become very small in a set of tiny measure (and even vanish somewhere inside), and such a tiny set may just happen to be the non-coincidence set $u^{-1}(M)$ as we are close to $\Sigma(u)$ (recall \eqref{eq:nulldensity}): e.g., $u^{-1}(M)$ is a cusp with its tip on $\Sigma(u)$, such that one of the tangent maps $\phi$ of $u$ is of rank zero along the direction of the cusp. The latter would certainly be the worst case scenario for us, as the rightmost side of \eqref{eq:reverse1} sees $u^{-1}(M)$ only. 

Thus, our analysis calls for a more robust approach that works in the absence of pointwise control of the vanishing behavior of $|Du|$. For this purpose, our new key result is an $L^{-\gamma}$-estimate of $|Du|$ near $\Sigma(u)$. Note that the next statement does not require convexity of the obstacle.

 \begin{thm} 
 \label{thm:Le}
     Let $M\subset \R^m$ be a smooth domain with compact complement and let $u\in W^{1,2}(B_4;\overline M)$ be a minimizing constraint map satisfying $\| u \|_{W^{1,2}(B_4)}\leq \Lambda$. There exist $\e_0, \gamma > 0$, both depending only on $n$, $m$ and $\partial M$, and a constant $c>1$, depending further on $\Lambda$, such that if $\int_{B_1(y)} |Du|^2\,dx \geq \e_0^2$ for some $y\in B_2$ then 
\begin{equation}
\label{eq:Le1}
     \int_{B_1} |Du|^{-\gamma}\,dx \leq c.
\end{equation}
 \end{thm}

This estimate can be understood as a quantitative version of the UCP for the constraint maps, as it shows that $|Du|$ cannot vanish too fast around $\Sigma(u)$. Once such an estimate is available, one can use the generalized H\"older inequality to verify the reverse link \eqref{eq:holderdecay} (see \eqref{eq:key3}), which in turn yields \eqref{eq:inc} as desired. Therefore, Theorem \ref{thm:Le} lies at the heart of our analysis for Theorem \ref{thm:uni-reg}.

\begin{rem} \label{thm:Le-harm}
As already observed above, Theorem \ref{thm:Le} works for general smooth compact obstacles $M^c$ without any condition on their principal curvatures. Moreover, this theorem also holds for minimizing harmonic maps into smooth compact manifolds without boundary, and it is even new in the latter context. 
\end{rem}

    The above results of UCP-type are somewhat reminiscent of the classical quantization results for minimizing harmonic maps \cite{SaU,ScU1}, which assert  that there is a universal lower bound on the energy of smooth harmonic maps which are not homotopic to a constant.  The theorem  is also evocative of a result from \cite{HL}, where the authors prove an $\cH^{n-2}$-measure estimate on the set $\cB_0(u)$ for a minimizing harmonic map $u$; however,  this estimate is qualitatively different from \eqref{eq:Le1}, and neither of them implies the other.

    We emphasize that the structure of blow-ups of harmonic maps near their discontinuity points is only well-understood in some very specific instances, such as when $n=3$ and $N = \mathbb S^2$ \cite{BCL},  or $n=4$ and $N = \mathbb S ^3$ \cite{N}. In these cases, the structure of the blow-ups can be exploited to prove \eqref{eq:Le1}, as we will see in Section \ref{sec:axial-reg}. However, since blow-ups are poorly understood in general, we are forced to rely on a fundamentally different type of argument to prove Theorem \ref{thm:Le}.

    The key observation for proving Theorem \ref{thm:Le} is the existence of a \textit{critical scale} at every regular point such that,
    at this scale, we can control both the regularity and the frequency of $|Du|$ uniformly and in our favor. The critical scale plays a  significant role in two ways. Above this scale, the effect of the discontinuous singularity dominates, so that $|Du|\gtrsim 1$ in a large portion of the neighborhoods. Below the critical scale we have $|Du| \ll 1$, but we can still control the doubling property of $|Du|$ in a uniform way. This is because the regularity of $u$ is uniformly controlled on such a scale. In other words, $1/|Du|$ does not grow so fast in a large portion of the neighborhoods on these scales. Finally, we can patch these two different scales altogether via a Calder\'on--Zygmund decomposition argument, which leads us to the desired $L^{-\gamma}$-estimate for $|Du|$. 
    
    Let us address that the notion of the critical scale here is closely related to that of the \textit{regularity scale} introduced in \cite{CN,NV1} in the study of the rectifiability properties of $\Sigma(u)$.  However, the analysis in these works is different from ours,  as the authors there are concerned with obtaining (optimal) $L^p$-estimates for $|Du|$ with $p > 1$.

\medskip

Returning to Theorem \ref{thm:uni-reg}, let us now address the optimality of the uniform convexity condition. Our next result asserts that this condition is essentially optimal,  as the conclusion of Theorem \ref{thm:uni-reg} no longer holds as soon as $\partial M$ contains a portion of a hyperplane, \textit{even when $M^c$ is convex}.   We remark that by \cite[Theorem 6.1]{FGKS}, the convexity of $M^c$ ensures a well-defined free boundary even in the presence of discontinuities, as it implies that $\dist(u,\partial M)$ is continuous.  In fact,  we can show that for obstacles with flat sides, minimizing constraint maps whose boundary values have non-trivial topological degree (see Appendix \ref{ap:degree} for the definition) \textit{always} develop discontinuities on their free boundaries:

\begin{thm}\label{thm:flat}
Let $\Omega\subset \R^n$ be a  smooth uniformly convex  domain, and let $M\subset \R^n$ be a smooth domain with convex complement such that  $M^c\Subset \Omega$ and $\partial M$ contains a relatively open subset $T$ of a hyperplane. Let $n\geq 3$ and  consider boundary data $g\in C^\infty(\partial\Omega;\partial\Omega)$ with $\deg(g) \neq 0$. If $u \in W_g^{1,2}(\Omega;\overline M) = \{ u \in W^{1,2}(\Omega;\overline M); \ u = g \ \hbox{on } \partial \Omega \} $ 
is a minimizing constraint map, then 
$$\partial u^{-1}(M)\cap \Sigma(u) \neq \emptyset.$$
In fact,   for every $y\in T$ there is $\{x_k\}\subset \Omega\cap u^{-1}(M)\setminus \Sigma(u)$ such that $x_k\to x_0 \in \partial u^{-1}(M)\cap \Sigma(u)$ and $u(x_k)\to y$, and so we write
$$T\subseteq u( \partial u^{-1}(M)\cap \Sigma(u)).$$  
\end{thm}

\begin{rem}
The assumption $n\geq 3$ in Theorem \ref{thm:flat} is used only to guarantee the existence of a minimizing constraint map. See Remark \ref{rem:dirclass} below for further discussions.
\end{rem}

The situation for non-convex obstacles is even more exotic than the one described above: as a byproduct of Theorem \ref{thm:flat},  we show in Corollary \ref{cor:flat} that, for such obstacles,  the non-coincidence set $u^{-1}(M)$ may \textit{not} be open.\footnote{It is not immediately clear how to define $u^{-1}(M)$ in general, since $u$ may develop discontinuities. A robust definition, which is independent of the choice of representative for $u\in W^{1,2}(\Omega;\overline M)$, is to set 
$$u^{-1}(M):=\{x\in \overline \Omega: \liminf_{r\to 0}\|\rho \circ u\|_{L^\infty(B_r(x))}>0\}.$$
This set contains those points in $\Omega\setminus \Sigma(u)$ which are mapped by (the smooth representative of) $u|_{\Omega\setminus \Sigma(u)}$ to $M$.}
This yields in particular the sharpness of \cite[Theorem 6.1]{FGKS}.

Producing examples as those in Theorem \ref{thm:flat} is challenging, since we do not know how to construct \textit{explicit} minimizing maps  in non-perturbative regimes. Indeed,  due to the vectorial nature of the problem, minimizers generally do not inherit the symmetry of their boundary data; see e.g.\ \cite{AL} for a dramatic symmetry-breaking phenomenon in the classical case of harmonic maps from the unit ball $\mathbb B^3$ to $\mathbb S^2$.  In particular, our approach to Theorem \ref{thm:flat} is based on a topological degree argument and is necessarily global, in contrast to the typical local analysis performed in the context of free boundary problems.

\medskip

The above results give a rather complete picture of the behavior of minimizing constraint maps near their free boundaries but do not address the structure of the free boundary itself.  The main difficulty in the analysis is due to the possible occurrence of branch points on the free boundary,  and in particular of points in $\cB_0(u)$: at such points, no matter what the obstacle is, the right-hand side in \eqref{eq:main-sys} vanishes and so we are faced with a \textit{degenerate} obstacle problem. On the other hand, if the obstacle is uniformly convex, then by Theorem \ref{thm:uni-reg} and the classical theory for the obstacle problem, the structure of the free boundary is well-understood in $\partial u^{-1}(M)\setminus \cB_0(u)$, see \cite{FKS}. 
We thus arrive at the following questions: \textit{For uniformly convex obstacles, is the set $\partial u^{-1}(M)\cap \cB_0(u)$ non-empty? If so, how does the free boundary look like around such points?}

The structure of the free boundary for degenerate free boundary problems is a difficult topic. In this direction, very little is known even in the scalar setting, except e.g.\ in the particular examples studied in \cite{McCurdy2024,Yeressian2015,Yeressian2016}. However, the degeneracy in these references stems from an absorbing term acting as an additional force in the problem; see also the discussions in \cite{Caff98-revisited} for a forcing term given directly by the obstacle. Nevertheless, the potential appearance of branch points for vectorial maps is part of the intrinsic nature of the problem and is not dictated by external forces. 

Here we construct and analyze a large class of examples of mininimizing constraint maps $u\in W^{1,2}(\Omega;\R^3\backslash\mb B^3)$ with non-trivial free boundaries such that
$$
\partial u^{-1}(\R^3\backslash \mb B^3)\cap \cB_0(u)\neq \emptyset,
$$
thus showing that branch points can occur on free boundaries, even when the obstacle is a ball. 
In order to construct such examples, we shall consider maps that minimize the Dirichlet energy in the more restricted class of $k$-axially symmetric maps;
see Definition \ref{def:axial}. Briefly, a map is called $k$-axially symmetric if it maps each cylinder centered on the $z$-axis to another such cylinder while rotating it $k$-times. Our analysis shows that, for $k\geq 2$, each free boundary point on the $z$-axis is simultaneously a branch point and a \textit{singular free boundary point}: thus, branching can occur and can be responsible for free boundary singularities.

\begin{thm}\label{thm:axial}
    Fix  $k\geq 2$, $\lambda > 1$, and a $k$-axially symmetric map $g\in C^2(\Ss^2;\Ss^2)$ with $\lvert\deg(g)\rvert = k$. Let $u\in W_{\lambda g}^{1,2}(\mb B^3;\R^3\backslash \mb B^3)$ be an energy minimizer among $k$-axially symmetric maps. For any $x_0\in\partial\{|u|>1\}$, there is a ball $B\subset \mb B^3$ centered at $x_0$ such that $u\in C^{1,1}(B)$ is the unique minimizing constraint map in $W^{1,2}_u(B;\R^3\backslash \mb B^3)$. Moreover, the following hold for the set $\cB_0(u)$ of rank zero points and the free boundary $\partial\{|u|>1\}$:
    \begin{enumerate}
        \item\label{it:nontrivFB} $\cB_0(u)\cap\partial\{|u|>1\}\cap\partial({\rm int}\{|u|=1\})\neq\emptyset$;
        \item\label{it:vanishorderDu} $|Du|$ vanishes exactly with order $k-1$ at every point of $\cB_0(u)$;
        \item\label{it:singFB} $\partial\{|u|>1\}$ cannot be represented as a $C^1$-graph in any neighborhood of $\cB_0(u)$ on the $z$-axis.
    \end{enumerate}
    In particular, (\ref{it:nontrivFB}) and (\ref{it:singFB}) show that the free boundary is non-trivial and singular.
\end{thm}

We refer the reader to Section \ref{sec:axial-bp}, and in particular to Theorem \ref{thm:bp-precise} therein, for a much more precise description of the structure of the free boundary near points in $\cB_0(u)$. As we shall see, in the setting above we can give a complete classification of the possible blow-ups appearing in the free boundary analysis.


\subsection{Organization of the paper}

The paper is organized as follows. In Section~\ref{sec:prelim}, we describe the setting and notation of this paper and collect preliminary results on minimizing constraint maps. In Section~\ref{sec:dist}, in preparation for our first main result, we prove the uniform continuity of the distance of a minimizing constraint map to the convex hull of the obstacle. Section~\ref{sec:Le} is devoted to the proof of Theorem~\ref{thm:Le}. In Section~\ref{sec:reg} we prove Theorem~\ref{thm:uni-reg}. In Section~\ref{sec:flat}, we study minimizing constraint maps for non-uniformly convex obstacles, and in particular we prove Theorem~\ref{thm:flat}. 
Section \ref{sec:branch} collects some remarks and examples concerning branch points and their interaction with free boundaries. In Section \ref{sec:axial-reg} we develop the regularity theory of axially symmetric constraint maps, and in particular we prove analogues of Theorems \ref{thm:uni-reg} and \ref{thm:Le}  in this setting. Finally, in Section \ref{sec:axial-bp} we study branch points on free boundaries for axially symmetric maps and prove Theorem \ref{thm:axial}. 

The paper also contains four appendices consisting of auxiliary results that are either known or follow with some relatively mild modifications from standard arguments. In Appendix \ref{ap:degree} we recall the definition and some basic properties of the topological degree, and its local variant. In Appendix \ref{app:freq} we recall standard facts about the Almgren frequency of a solution to a regular linear elliptic system. In Appendix \ref{app:strong-ax} we prove a version of Luckhaus' lemma in the setting of axially symmetric constraint maps. Finally, in Appendix \ref{app:unique} we prove the uniqueness and regularity of weakly constraint maps into small geodesic balls.

\subsection{Notation}\label{sec:notation1}

For the reader's convenience, we gather here a list of notation to be used throughout this paper.

\begin{longtable}[l]{l l }
\,\,\,\,$\Omega$ &\qquad  a bounded domain in $\R^n$ ($n\geq 1$)  \medskip \\ 
\,\,\,\,$\Omega_\eta $ &\qquad  $\{ x \in\Omega:\dist(x,\partial \Omega)>\eta\}$  ($\eta \geq 0$)    \medskip \\ 
\,\,\,\,$M$  &\qquad smooth domain in $\R^m$ ($m\geq 2$)  \medskip \\ 
\,\,\,\,$u$  &\qquad $(u^1, \cdots , u^m)$    \medskip \\
\,\,\,\,$W^{1,p}(\Omega;\overline M)$&\qquad $\{u\in W^{1,p}(\Omega;\R^m): u\in \overline M\text{ a.e.\ in }\Omega\}$      \medskip \\
\,\,\,\,$W^{1,p}_g(\Omega;\overline M)$&\qquad  the Dirichlet class $W^{1,p}(\Omega;\overline M)\cap (g+W^{1,p}_0(\Omega;\R^m))$, for $g\in W^{1,p}(\Omega;\overline M)$   \medskip \\
\,\,\,\,$\mb B^m$  &\qquad the open unit ball in $\R^m$  \medskip \\  
\,\,\,\,$\mb S^{m-1}$  &\qquad the unit sphere in $\R^m$, $\mb S^m=\partial \mb B^m$ \medskip \\ 
\,\,\,\,$A_p (\cdot , \cdot)$  &\qquad     the second fundamental form of $\partial M$ at $p$, see \eqref{eq:hessrho}   \medskip \\   
\,\,\,\,$\cN(\partial M)$  &\qquad    tubular neighborhood  of $\partial M$ 
    \medskip\\   
\,\,\,\,$\rho$  &\qquad   signed distance function to $\partial M$, smooth on $\cN(\partial M)$    \medskip \\   
\,\,\,\,$\nu=\nabla \rho$  &\qquad inward unit normal to $M$, extended to $\cN(\partial M)$ \medskip \\   
\,\,\,\,$\Pi$   &\qquad  nearest point projection: $ \cN(\partial M)\to \partial M$
       \medskip\\ 
\,\,\,\,$\nu\otimes \nu$   &\qquad        the tensor product matrix with entries $(\nu^i\nu^j)_{ij}$ \medskip\\   
\,\,\,\,$\xi^\top$  &\qquad   $
  (I - \nu\otimes \nu)\xi
$,   \ orthogonal projection of $\xi$ into $T(\partial M)$, see \eqref{eq:xitop}      \medskip \\ 
\,\,\,\,$\cE[u]$  &\qquad $\int_\Omega |Du|^2 \, dx$   \medskip \\ 
\,\,\,\,$E(u,x_0,r)$  &\qquad   $ r^{2-n} \int_{B_r(x_0)} |Du|^2\,dx$, \ normalized energy     \medskip\\
\,\,\,\,$N(f,x_0,r)$  &\qquad   $r\int_{B_r(x_0)} |Df|^2\,dx/\int_{\partial B_r(x_0)} |f - (f)_{x_0,r}|^2 \haus$, the frequency function     \medskip\\
\,\,\,\,$\Sigma(u)$   &\qquad  $ \{ x\in\Omega: u\text{ is not continuous at }x\} $    \medskip \\   
\,\,\,\,$ \cB_k(u)$   &\qquad  $\{x\in \Omega\setminus \Sigma(u):Du\text{ has rank $k$ at }x\}$  \medskip \\ 
\,\,\,\,$\cH^k$  &\qquad  $k$-dimensional Hausdorff measure    \medskip \\   
\,\,\,\,$\hbox{dim}_\cH$  &\qquad  Hausdorff dimension     \medskip\\    
\,\,\,\,$F:G$  &\qquad  $ f_\alpha^i g_\alpha^i$, \ where  $F = (f_\alpha^i)$ and  $G= (g_\alpha^i)$ are matrices in $\R^{m\times n}$    \medskip \\ 
\,\,\,\,$D^k$ &\qquad the $k$-th order differential operator in the ambient space, $D_\alpha := \frac{\partial}{\partial x_\alpha}$ \medskip \\
\,\,\,\,$\nabla^k$ &\qquad the $k$-th order differential operator in the target space, $\partial_i := \frac{\partial}{\partial y^i}$ \medskip \\
\,\,\,\,$\Hess f_y(\xi,\zeta)$  &\qquad $\partial_{ij} f(y)\xi_\alpha^i\zeta_\alpha^j$, for a function $f$ from the target space, and $\xi,\zeta \in \R^{m\times n}$ \medskip\\ 
\end{longtable}


\section{Preliminaries}\label{sec:prelim} 


\subsection{Problem setting}\label{sec:notation2}

\begin{figure}[ht]
    \centering
    \includegraphics[scale=0.25]{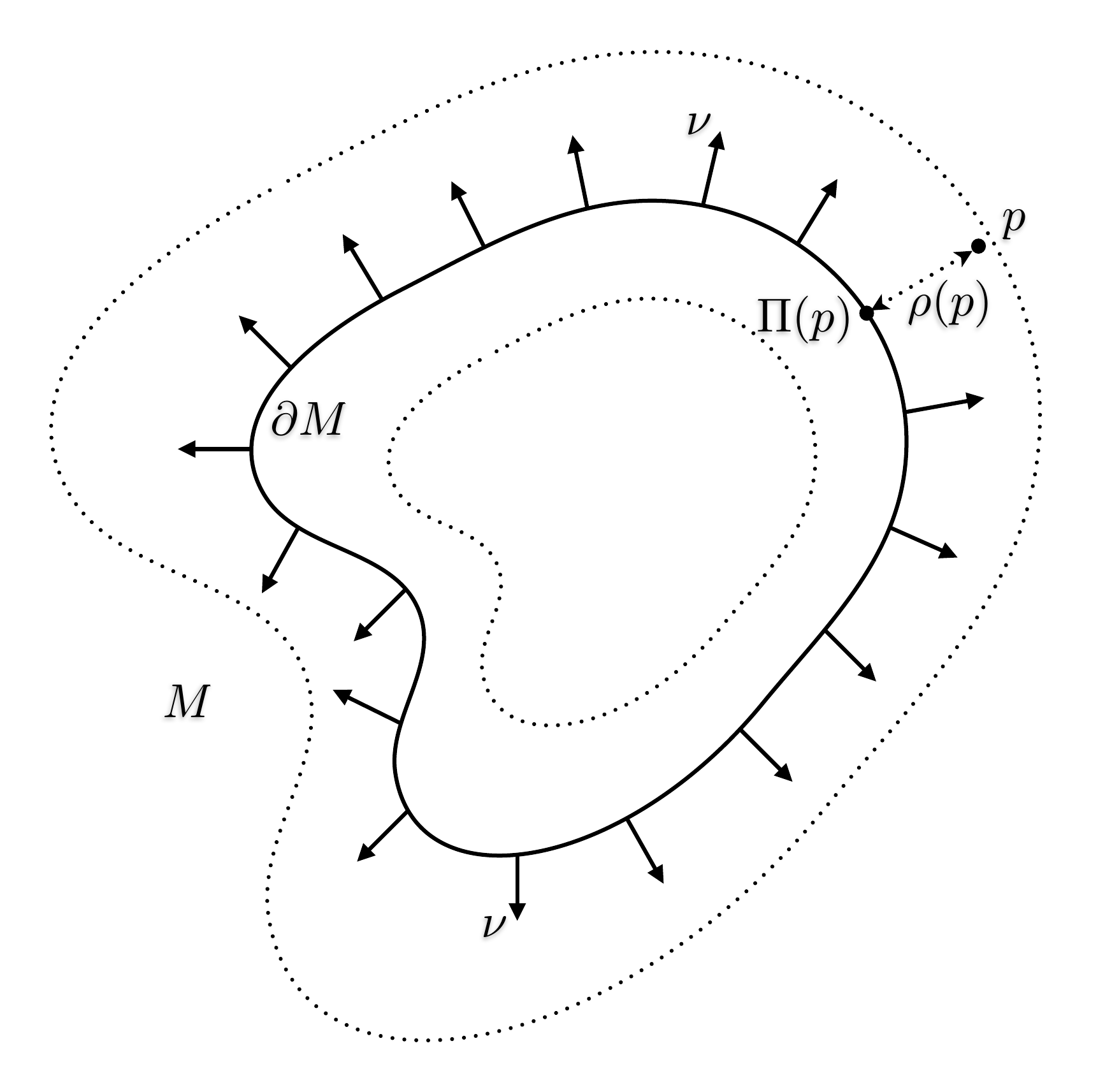}
    \caption{A representation of the tubular neighborhood $\cN(\partial M)$ of the boundary of a smooth domain $M\subset \R^m$ with compact complement.}
    \label{fig:manifold}
\end{figure}

Throughout this paper, $\Omega\subset \R^n$ is a bounded domain and $M\subset \R^m$ is a smooth domain; we always assume that $n,m\geq 2$. The  complement of $M$, $M^c$, is called the \textit{obstacle}, and we always assume that it is compact. We denote by $\rho$ the (signed) distance function to $\partial M$ that takes positive values in $M$, and by $\Pi:\R^m \to \partial M$ the nearest point projection onto $\partial M$. We also write $\cN(\partial M)$ for the largest open neighborhood of $\partial M$ where $\rho$ and $\Pi$ are smooth; in particular, in the case where $M^c$ is convex, the set $\cN(\partial M)$ strictly contains $\overline M$.  We set $\nu:=\nabla \rho$ on $\cN(\partial M)$ and note that, on $\partial M$, $\nu$ is the inward-pointing unit normal to $\partial M$. See Figure \ref{fig:manifold} for an illustration.

We denote by $B_r(x_0)$ the ball of radius $r$ centered at the point $x_0$, and more generally $B_r(F)$ denotes the neighborhood of width $r$ around the set $F$.  If $B=B_{r}(x_0)$ is a ball, we write $\lambda B := B_{\lambda r}(x_0)$
for the dilation of $B$ by a factor $\lambda>0$ around its center; we use a similar notation for cubes. We also denote by $\dim_\cH(S)$ the Hausdorff dimension of a set $S$, by $\cH^k(S)$ the $k$-dimensional Hausdorff measure of $S$, and by $|S|$ the Lebesgue measure of $S$. We write $D$ and $\nabla$ for the usual differential operators in the ambient space and target spaces, respectively.  
Also, for a function $f,$ we write $\Hess f_y(\xi,\zeta):= \partial_{ij} f(y)\xi_\alpha^i\zeta_\alpha^j.$


\subsection{Known results for minimizing constraint maps}\label{sec:known}

To keep our exposition as self-contained as possible,  we gather here some basic results about minimizing constraint maps. These can be seen as generalizations of classical theorems concerning minimizing harmonic maps into closed manifolds \cite{LW,Simon}. 

A basic quantity in the study of minimizing constraint maps is the \textit{normalized energy}
\begin{equation}\label{eq:Eq}
E(u,x_0,r) := r^{2-n} \int_{B_r(x_0)} |Du|^2\,dx.
\end{equation}
Minimizing constraint maps are stationary upon inner variations, and this yields a monotonicity formula for the normalized energy. A more general result is proved in \cite{FGKS}; here, we only present the statement adapted to the current setting. 

\begin{lem}[Monotonicity formula]\label{lem:monot}
Let $u\in W^{1,2}(\Omega;\overline M)$ be a minimizing constraint map. For all $x_0\in\Omega$ we have
$$
\begin{aligned}
E (u, x_0,{s}) - E (u, x_0,r) \geq 2 \int_{B_{s}(x_0)\setminus B_r(x_0)} |x-x_0|^{2-n} \left|\frac{\partial u}{\partial N}\right|^2 dx 
\qquad \forall \ \, 0<r<s<\dist(x_0,\partial\Omega),
\end{aligned}
$$
where $\frac{\partial}{\partial N}$ denotes the directional derivative in the radial direction. In particular, if $E(u,x_0,{s}) = E(u,x_0,r)$ then $u$ is $0$-homogeneous around $x_0$ inside $B_{s}(x_0)\setminus B_r(x_0)$, i.e.,\ 
$$(x-x_0) \cdot Du = 0\quad \text{ a.e.\ in } B_{s}(x_0)\setminus B_r(x_0).$$
\end{lem}

A crucial tool in our analysis is the compactness of minimizing constraint maps, originally established in \cite{L} within a broader context. Since we will frequently reference it here, we provide the adapted statement for the reader's convenience.

\begin{lem}[Compactness]\label{lem:strong}
Let $\{u_k\}_{k=1}^\infty\subset W^{1,2}(\Omega;\overline M)$ be a bounded sequence of minimizing constraint maps. Then there exists a minimizing constraint map $u \in W^{1,2}(\Omega;\overline M)$ such that, up to a subsequence, $u_k \to u$ strongly in $W_\loc^{1,2}(\Omega;\R^m)$. 
\end{lem}

The monotonicity of the normalized energy yields a well-defined limit $E(u,x_0,0^+)$ for every $x_0\in\Omega$. 
As an immediate consequence of Lemmas~\ref{lem:monot} and~\ref{lem:strong}, we obtain the upper semicontinuity of the energy density. 

\begin{cor}[Upper semicontinuity]\label{cor:semi}
Let $\{u_k\}_{k=1}^\infty\subset W^{1,2}(\Omega;\overline M)$ be a sequence of minimizing constraint maps such that $u_k\to u$ strongly in $W^{1,2}(\Omega;\R^m)$. Then for any sequence $\{x_k\}_{k=1}^\infty\subset \Omega$ with $x_k \to x\in\Omega$ we have 
$$
E(u,x,0^+) \geq \limsup_{k\to\infty} E(u_k,x_k,0^+). 
$$
\end{cor}

We will also make use of the $\e$-regularity theorem, which states that minimizing constraint maps are regular in the vicinity of points with small energy. The optimal version of this result is due to the recent work \cite{FKS} (see also \cite{D} for a previous version of this statement).

 \begin{thm}[$\e$-regularity]\label{thm:e-reg}
There exist a small constant $\e_0 > 0$ 
depending only on $n$, $m$, and $\partial M$, such that the following holds: if $u\in W^{1,2}(\Omega;\overline M)$ is a minimizing constraint map with
 $$E(u,x_0,2{r}) \leq \e_0^2, \quad \text{ where } B_{2{r}}(x_0)\Subset\Omega, $$ 
 then $u\in C^{1,1}(B_r(x_0))$, with the estimate  
 $$
 \| D^j u\|_{L^\infty(B_{r}(x_0))} \leq \frac{1}{{r}^j},\quad j=1,2.
 $$
 \end{thm}

Recall from \eqref{eq:sing} that $\Sigma(u)$ consists of all discontinuity points of the map $u$ in the prescribed domain. A useful corollary to Theorem \ref{thm:e-reg} is the energy quantization at points of $\Sigma(u)$:

\begin{cor}[Energy quantization]\label{cor:e-reg}
There exists a constant $\e_0 > 0$, depending only on $n$, $m$, and $\partial M$, such that for any minimizing constraint map $u\in W^{1,2}(\Omega;\overline M)$ we have that $x_0\in \Sigma(u)$ if and only if $E(u,x_0,0^+) \geq \e_0^2$. 
\end{cor}

Through the usual dimension-reduction principle,  the above results yield an optimal partial regularity theorem, see \cite{DF,FKS}.

\begin{thm}[Partial regularity]\label{thm:dim-sing}
Let $u\in W^{1,2}(\Omega;\overline M)$ be a minimizing constraint map. Then $u\in C_\loc^{1,1}(\Omega\setminus\Sigma(u))$ and $\dim_\cH(\Sigma(u))\leq n-3$.
\end{thm}

The following lemma can be found in \cite[Lemma 2.2]{D}, and will be used frequently in our subsequent analysis.

\begin{lem}\label{lem:dist-pde}
Let $u\in W^{1,2}(\Omega;\overline M)$ be a minimizing constraint map such that $u(B)\subset \cN(\partial M)$ a.e.\ in a ball $B\subset\Omega$. Then 
    \begin{equation}\label{eq:dist-pde}
        \Delta (\rho\circ u) =  \Hess \rho_u(Du,Du)\chi_{u^{-1}(M)}\quad\text{in }B,
    \end{equation}
    in the weak sense. 
\end{lem}

\subsection{Vector calculus identities and geometric bounds}

In this subsection we prove some elementary identities and estimates related to the geometry of $\partial M$. 
Since $\nu = \nu \circ \Pi$ in $\cN(M)$, we have the decomposition 
\begin{equation}
\label{eq:decomp}
\id = \Pi+ \rho \nu = \Pi+\rho(\nu \circ \Pi).
\end{equation}
In particular, differentiating both sides of this identity, we have $\Id - \nu\otimes \nu = (\Id + \rho\nabla \nu)\nabla\Pi$; note that $\Id -\nu \otimes\nu$ is the orthogonal projection onto the tangent bundle $T(\partial M)$.
For a vector $\xi\in \R^m$, we thus write
    \begin{equation}
\label{eq:xitop}
\xi^\top := \xi - (\nu\cdot\xi)\nu= [(\Id + \rho \nabla\nu)\nabla\Pi]\xi
\end{equation}
 for the orthogonal projections of $\xi$ into $T(\partial M)$.
 
The second fundamental form of $\partial M$, appearing in \eqref{eq:main-sys}, can be expressed at points in $\partial M$ as

\begin{equation}
\label{eq:hessrho}
\Hess \rho(\xi,\xi) = \Hess \rho(\xi^\top,\xi^\top) = -\nu \cdot A(\xi^\top,\xi^\top).
\end{equation}
When the obstacle $M^c$ is uniformly convex, we obtain a simple relation between $\Hess\rho$ and the differential $|\nabla\Pi|$:

\begin{lem}\label{lem:Hessrho}
Let $M$ be a smooth domain with uniformly convex complement. Then 
\begin{equation}\label{eq:Hessrho}
\Hess\rho(\xi,\xi) \geq \frac{\kappa}{1+\kappa\rho} |(\nabla\Pi)\xi|^2\quad \forall\,\xi\in\R^m, 
\end{equation}
where $\kappa$ is the curvature lower bound for $\partial M$. 
\end{lem}

\begin{proof}
As $M^c$ is convex, $\rho$, $\nu$, and $\Pi$ are well-defined and smooth in a neighborhood of $\overline M$. Let $\{\tau^1,\cdots,\tau^{m-1},\nu\}$ be the principal coordinate system of $\partial M$, and $\kappa_i$ be the $i$-th principal curvature. As $\nu$ points outward to $M^c$, it follows from \cite[Lemma 14.17]{GT} that 
\begin{equation}\label{eq:hessrho-1}
\Hess\rho(\tau^i,\tau^i) = \frac{\kappa_i\circ\Pi}{1 + (\kappa_i\circ \Pi)\rho}\quad\text{in }\overline M. 
\end{equation}
Also, as $M^c$ is uniformly convex, $\kappa_i\circ\Pi\geq \kappa$ in $\overline M$, for some $\kappa>0$ which depends only. Since the map $t\mapsto \frac{t}{1+ at}$ is strictly increasing in $t>0$, for each $a \geq 0$, we deduce from \eqref{eq:hessrho-1} that
\begin{equation}\label{eq:hessrho-2}
\Hess\rho(\xi,\xi)\geq \frac{\kappa}{1+\kappa\rho} |\xi^\top|^2\quad\text{in }\overline M\quad\forall\, \xi\in\R^m. 
\end{equation} 
Since $\nabla \nu = \Hess\rho$ (as a consequence of the relation $\nu = \nabla\rho$), using again \eqref{eq:hessrho-1} we see that
\begin{equation}\label{eq:hessrho-3}
(\Id + \rho\nabla\nu)(\tau^i,\tau^i) = \frac{1 + 2(\kappa_i\circ\Pi)\rho}{1+ (\kappa_i\circ\Pi)\rho} \in [1,2]\quad\text{in }\overline M.
\end{equation}
In particular, combining \eqref{eq:hessrho-3} with \eqref{eq:xitop} yields 
\begin{equation}\label{eq:hessrho-4}
|\xi^\top|^2\geq  |(\nabla \Pi)\xi|^2\quad\text{in }\overline{M}\quad\forall\,\xi\in\R^m.
\end{equation}
Thus, \eqref{eq:Hessrho}  follows from \eqref{eq:hessrho-2} and \eqref{eq:hessrho-4}. 
\end{proof}


\subsection{Existence and examples of constraint maps}

Consider boundary data $g\in C^\infty(\partial\Omega; \overline M)$. In this setting, the direct method yields existence of minimizing constraint maps in $W^{1,2}_g(\Omega;\overline M)$, \textit{provided this space is non-empty}. 
We discuss whether this holds in the next remark and for simplicity we take $\Omega=\mb B^n$ to be the open unit ball.

\begin{rem}\label{rem:dirclass}
In studying $W^{1,2}_g(\mb B^n;\overline M)$, there are two different cases to consider.
    \begin{enumerate}
    \item If $ n\geq 3$, then $W_g^{1,2}(\mb B^n;\overline M)\neq \emptyset$. Indeed, it suffices to consider the 0-homogeneous extension $x\mapsto g(x/|x|)$ to obtain a map in that class. The crucial calculation is that the map $x\mapsto x/|x|$ is in $W^{1,p}$ near the origin if and only if $n>p$. Thus, in order to be able to take $p=2$, we require that $n\geq 3$.
    \item If $n=2$, the situation is more subtle and it can be concisely expressed in terms of the homotopy extension property: $W^{1,2}_g(\mb B^2;\overline M)\neq \emptyset $ if and only if
    \begin{equation}
        \label{eq:homotopyextension}
      g \text{ is homotopic to } G|_{\mb S^{1}}\text{ for some continuous map } G\colon \overline{\mb B^2}\to \overline M.
    \end{equation}
    Note that this characterization holds for boundary data $g$ that are not necessarily smooth, but are just in the natural trace space $H^{\frac 1 2}(\mb S^{1};\overline M)$ \cite[Theorem 2]{Bethuel1995}.
    
    In the simplest possible case where 
    $g(\mb S^{1})$ is diffeomorphic to $\mb S^1$, \eqref{eq:homotopyextension} can be characterized elegantly in terms of the topological degree of $g$ (we refer the reader to Appendix \ref{ap:degree} for a brief introduction to this notion). Indeed, according to Hopf's Theorem\footnote{The Hopf Theorem is a result within differential topology asserting that the topological degree stands as the sole homotopy invariant among continuous mappings taking values in a sphere.}, in this case
    $$\eqref{eq:homotopyextension} \quad \iff\quad \deg (g) = 0.$$
    For data $g\in H^{\frac 1 2}(\mb S^1;\overline M)$ that are not continuous, such equivalence still holds, but one needs to interpret the degree in the sense of the VMO-degree of Brezis and Nirenberg \cite{Brezis1995}.
\end{enumerate}
\end{rem}

We now construct an explicit example of a radially symmetric constraint map with a discontinuity at the origin.

\begin{ex}\label{ex:radial}
Given $n=m\geq 3$ and a radius $a \in (0,1)$, we aim to find a map $u_a\in C^{1,1}(\overline{B_1}\setminus\{0\}; B_a^c)\cap W^{1,2}(B_1;B_a^c)$ such that $u_a = \id$ on $\partial B_1$
and which solves
        $$
    \Delta u_a = - |Du_a|^2 u_a\chi_{\{|u_a|=a\}}\quad\text{in }B_1,
    $$    
i.e., the Euler--Lagrange system for constraint maps with  obstacle $\overline{ B_a}$, see \eqref{eq:main}.
   
   We look for a map $u_a\colon \overline{B_1}\setminus \{0\}\to B_a^c$ of the form 
    $$
    u_a (x) := w_a(|x|) \frac{x}{|x|}.
    $$
The above system yields an overdetermined ODE for $w_a$, namely
     $$
        \begin{cases}   \ddot w_a + \frac{n-1}{r} \dot w_a = \frac{n-1}{r^2} w_a \quad  \text{ for $r \in (r_a,1)$},\\
        w_a(r_a) = a, \quad \dot w_a(r_a) = 0, \quad w_a(1) = 1, \\
        w_a(r)=a\quad \text{ for $r \in [0,r_a]$,}\qquad 
    \end{cases}
    $$
    where $r_a \in (0,1)$ is a free parameter corresponding to the free boundary constraint. 
Ignoring the  boundary conditions, the above ODE  has two  particular solutions, namely  $r$ and $r^{1-n}$.
We thus  look for a solution of the form
  $$
    w_a (r) = 
    \left\{\begin{array}{ll}
a & \text{ for $r \in [0,r_a]$,}\\
t_a r + (1-t_a)r^{1-n} & \text{ for } r \in (r_a,1),
    \end{array}\right.
    $$
    where the parameters $t_a,r_a \in (0,1)$ satisfy
$$
t_ar_a+(1-t_a)r_a^{1-n}=a,\qquad t_a+(1-n)(1-t_a)r_a^{-n}=0.
$$
One can easily check that the system above has a unique solution. In addition, one can observe that $r_a\to 0^+$  as $a\to 0^+$, while $r_a\to 1^-$ as $a\to 1^-$.  
Returning to the original map $u_a$,  the coincidence set is given by $u_a^{-1}(\partial B_a) = \overline{B_{r_a}}$, with free boundary $\partial u_a^{-1}((\overline B_a)^c) = \partial B_{r_a}$. 

\end{ex}

\begin{rem}
An interesting problem is to understand whether the map $u_a$ is a \textit{globally minimizing} constraint map. We recall that the radial solution\footnote{Note that the function $u_1(x)= x/|x|$ corresponds to the limit of the functions $u_a$ as $a\to 1^-.$} $u_1(x)= x/|x|$ of the Dirichlet problem for harmonic maps into a sphere with identity boundary values is a minimizer, i.e.
$$\cE[u_1]\leq \cE[v] \qquad \text{for all } v\in W^{1,2}_\id(B^n;\mathbb S^{n-1}),$$
see e.g.\ \cite[Proposition 2.2.1]{LW} and the references therein; here, as before, we assume that $n\geq 3$.   However, despite the intuitive nature of this result, it is quite non-trivial to show the minimality of $u_1$.
  By analogy with the case of harmonic maps, it seems plausible  that the maps $u_a$ above should be minimizers. 
  
\end{rem}

We now give a few more different examples of minimizing constraint maps.

\begin{ex}
    As mentioned in the introduction, (energy-minimizing) geodesic curves in Riemannian manifolds-with-boundary are (minimizing) constraint maps, see e.g.\ \cite{AB, ABB2,ABB1} for a study of their properties.
\end{ex}

\begin{ex}\label{ex:hull}
    Consider a smooth boundary condition $g\colon \partial \Omega\to \partial M$,  suppose that the obstacle $M^c$ is convex, and let $u\in W^{1,2}_g(\Omega;\overline M)$ be a minimizing constraint map. Then, according to \cite[Theorem 6.1]{FGKS}, see also Corollary \ref{cor:dist} below, the non-negative function $\rho \circ u\colon \overline \Omega\to [0,+\infty)$ is subharmonic and vanishes on $\partial \Omega$. Therefore, by the maximum principle, $\rho\circ u=0$ in $\Omega$ and thus $u^{-1}(\partial M)=\overline \Omega$. It follows that $u$ is a \textit{minimizing harmonic map} into $\partial M.$

An alternative  geometric  proof (which avoids subharmonicity of $\rho \circ u$)  is to take any supporting plane $H(y) = (y-y_0) \cdot \nu_{y_0} $, with $y_0 \in \partial M$, with 
$M \subset \{ H < 0 \}$,  
and consider the restriction of the map to 
$\{ H \geq 0\}$. Applying the minimum principle to 
the harmonic function $H(u (x))$ in the set 
$\{ H > 0 \}$ allows us to conclude that $u^{-1}(\partial M) = \overline \Omega.$
\end{ex}

\begin{ex} In analogy with parametric minimal surfaces, one can show that 
    if the convex hull of the image of the boundary satisfies $\textup{conv}(u(\partial \Omega )) \subset M$ then $u(\Omega )\subset M$.   
Indeed,  as in Example \ref{ex:hull}, the distance function to the convex hull is subharmonic (as $u$ is harmonic in each component), so the maximum principle shows that the distance function vanishes everywhere, proving that $u(\Omega)\subset \textup{conv}(u(\partial\Omega))\subset M$.
\end{ex}

We conclude this section with a comparison between minimizing constraint maps and solutions to the scalar obstacle problem.

\subsection{Comparison with the scalar obstacle problem}
Given a scalar function $\vp\in C^\infty(\Omega)$, the classical obstacle problem consists in minimizing the Dirichlet energy $\cE$ among all functions $v \in W^{1,2}(\Omega)$ satisfying the constraint $v\geq \varphi$ in $\Omega$. Minimizers of this problem are precisely the solutions of the equation
\begin{equation}
    \label{eq:ELobstacle}
    \Delta v = (\Delta \varphi)\chi_{\{v=\vp\}}\quad\text{in }\Omega.
\end{equation}
This problem and many variations of it have been the subject of intensive research over several decades and from various perspectives. Several influential works that have shaped the recent landscape of the scalar theory are, for instance,
\cite{MR4054360,  MR3448773, MR4689552,  MR4695505, MR3904453, MR3840912, MR2511747, MR3388870, MR2601491,  MR3855748,  MR4655359, MR4595308, MR3980852}.

Here, we consider graphical obstacles and compare the solutions of the classical scalar obstacle problem to minimizing constraint maps. As we shall prove below, these two objects are always different unless the solution does not ``actively'' touch the obstacle.  Informally speaking, the reason for this difference comes from the fact that the graph of a solution to the obstacle problem is pushed \textit{upwards} by the obstacle, while the hypersurface parametrized by a minimizing constraint map is pushed \textit{away from} the obstacle.

To begin our analysis, 
we denote by $M_\vp$ the subgraph of $\vp$, i.e., $M_\vp = \{(\bar y,y_{n+1})\in\R^{n+1}: \vp(\bar y) < y_{n+1}\}$, fix a solution $v\in W^{1,2}(\Omega)$ to \eqref{eq:ELobstacle}, and consider the vectorial map $u\in W^{1,2}(\Omega;\overline{M_\vp})$ defined by 
\begin{equation}\label{eq:u-graph}
u(x) := (x,v(x)).
\end{equation}

Assume by contradiction that $u$ is a minimizing constraint map. Since $\vp\in C^2(\Omega)$, we have $v\in C_\loc^{1,1}(\Omega)$, whence $u\in C_\loc^{1,1}(\Omega;\overline{M_\vp})$. Since $Dv = D\vp$ on $\Omega\cap \{v = \vp\} = \Omega\cap u^{-1}(\partial M)$, we also have 
\begin{equation}\label{eq:Du}
D_\alpha u = e_\alpha + (D_\alpha \vp) e_{n+1}\quad\text{on }\Omega\cap u^{-1}(\partial M). 
\end{equation}
for every $\alpha \in\{1,2,\cdots,n\}$, where $e_\alpha$ is the standard $\alpha$-th basis vector in $\R^{n+1}$. Thus, denoting by $A_y(\cdot,\cdot)$ the second fundamental form of $\partial M_\vp$ (which is precisely the graph of $\vp$) at $y$, 
$$
A_u(Du,Du) = \frac{(-D\vp,1)}{1 + |D\vp|^2}\Delta \vp\quad\text{on }\Omega\cap u^{-1}(\partial M),
$$
with the specific choice $u$ as in \eqref{eq:u-graph}\footnote{In fact, setting $\psi:\Omega\to \R^{n+1}$ as the canonical parametrization of the graph, $\partial M_\vp$, of $\vp$, (i.e., $\psi(x) := (x,\vp(x))$) we obtain, for any $1\leq \alpha,\beta\leq n$,
$$
A_\psi(D_\alpha \psi, D_\beta \psi) = \frac{(-D\vp, 1 )}{1 + |D\vp|^2}D_{\alpha\beta}\vp\quad\text{in }\Omega.
$$} and using \eqref{eq:Du}.
Thus, if $u$ is a minimizing constraint map, it follows from the Euler-Lagrange system \eqref{eq:main-sys} as well as the definition \eqref{eq:u-graph} of $u$ that 
\begin{equation}\label{eq:ELugraph}
(\Delta v)e_{n+1} = \Delta u = -\frac{\Delta\vp}{1 + |D\vp|^2}\big((D_\alpha \vp)e_\alpha - e_{n+1}\big)\chi_{\{v = \vp\}}\quad\text{in }\Omega;    
\end{equation}
here, as usual, we used the summation convention for the repeated index $\alpha$ running through $\{1,2,\cdots,n\}$. Comparing \eqref{eq:ELobstacle} with \eqref{eq:ELugraph}, it follows that 
\begin{equation}\label{eq:graph-cond}
|D\vp|^2\Delta\vp = 0 \quad\text{a.e.\ on }\Omega\cap\{v=\vp\}.    
\end{equation}
Note that \eqref{eq:graph-cond} is a very strong restriction.
Indeed, since $D^2\varphi=0$ (and so $\Delta\varphi=0$) a.e. on the set $\{D\varphi=0\}$, \eqref{eq:graph-cond} is equivalent to implies
\begin{equation}\label{eq:graph-cond2}
\Delta\vp = 0 \quad\text{a.e.\ on }\Omega\cap\{v=\vp\}.    
\end{equation}
Combining this information with \eqref{eq:ELobstacle}, we deduce that $\Delta v=0$ a.e., thus $v$ is harmonic.

In other words, except in the trivial case when  $v$ is harmonic (which means that the obstacle is not acting on the function $v$), the graph of the solution to the obstacle problem does not coincide with the corresponding minimizing constraint map.


\section{Regularity of the distance to the convex hull of the obstacle}\label{sec:dist}

In this section, we establish the uniform regularity of the distance of a minimizing constraint map to the convex hull of the given obstacle. This result significantly improves \cite[Theorem 6.1]{FGKS}.
For convenience, we state the result below with $Y$ being an arbitrary closed convex set containing  $M^c$, but the result is particularly relevant when $Y$ is the convex hull of $M^c$.

\begin{thm}\label{thm:dist-reg}
Let $M$ be a smooth domain in $\R^m$ with compact complement, ${Y}$  a closed convex set containing  $M^c$, and $u\in W^{1,2}(\Omega;\overline M)$ a minimizing constraint map with $\| u \|_{W^{1,2}(\Omega)}\leq\Lambda$. Then for every $\delta > 0$ there exists a modulus of continuity $\omega$, depending only on $n$, $m$, $\partial M$, $\Lambda$, and $\delta$, such that 
\begin{equation}
    \label{eq:oscdist}
    \osc_{B_r(x_0)} \dist(u,Y) \leq \omega(r)\qquad\forall\,r \in(0,\delta),
\end{equation}
whenever $B_{2\delta}(x_0)\subset\Omega$. 
\end{thm}

We begin by showing that the distance function is subharmonic.

\begin{lem}\label{lem:dist-subharm}
 Let $Y$ be a closed convex set containing  $M^c$.   Then the function $\dist(u,Y)$ is weakly subharmonic in $\Omega$.
\end{lem}

\begin{proof}
  Define the function $w:\Omega\to \R$,  $w (x)  :=  \dist(u(x)  ,Y)$. Since $\{w  > 0\}\setminus \Sigma(u)\subset u^{-1}(M)\setminus \Sigma(u)$, it follows by \eqref{eq:main} that  $\Delta u = 0$ inside ${\rm int}(\{ w > 0\})\setminus\Sigma(u)$. Hence, since $y\mapsto \dist(y,Y)$ is convex (by the convexity of $Y$), we conclude that 
   $ \Delta w \geq 0$   inside  ${\rm int}(\{ w > 0\})\setminus\Sigma(u).$
 As $\dim_\cH \Sigma(u) \leq n-3$ by Theorem \ref{thm:dim-sing} and $w\in W^{1,2}(\Omega)$, by a capacity argument it follows that $\Delta w \geq 0$ weakly in $\Omega$. 
\end{proof}

As a corollary, we obtain an interior $L^\infty$-estimate for $u$:

\begin{cor}\label{cor:sup}
Under the hypothesis of Theorem \ref{thm:dist-reg}
there is  $a\in M^c$ such that $|u - a| \in L_\loc^\infty(\Omega)$ and 
$$
\| u -  a \|_{L^\infty(\Omega')}\leq c\diam M^c + \frac{ c\| u - a\|_{L^2(\Omega)}}{\dist(\Omega',\partial\Omega)^{n/2}} 
$$
for every subdomain $\Omega'\Subset\Omega$, where $c$ depends only on $n$.  
\end{cor}

\begin{proof}
    Since $M^c$ is compact, there is $a\in M^c$ such that $M^c\subset \overline B_R(a)$, with $R := \diam M^c$. Set ${Y} := \overline B_R(a)$, so that $\dist(u,Y) = (|u-a| - R)_+$. By Lemma \ref{lem:dist-subharm}, $\dist(u,Y)$ is weakly subharmonic in $\Omega$. The conclusion then follows from the local $L^\infty$-estimate for weakly subharmonic functions. 
\end{proof}

We also immediately obtain the following conclusion on the range of constraint maps:

\begin{cor}\label{cor:image}
Let $u \in W^{1,2}(\Omega;\overline M)$ be a minimizing constraint map. If $Y$ is the closed convex hull of $u(\partial\Omega)\cup M^c$, then $u(x)\in Y$ for a.e.\ $x\in \Omega$.
\end{cor}

\begin{proof}
This follows from the weak maximum principle by applying Lemma \ref{lem:dist-subharm} to $\dist(u,Y)$.
\end{proof}

In what follows, we shall fix a closed convex set ${Y}$ containing $M^c$. We shall denote by $\e_0= \e_0(n,m,\partial M)$ the small energy threshold from Corollary \ref{cor:e-reg}. In addition, $c$ will denote a generic positive constant depending at most on $n$, $m$, $\partial M$, and $\Lambda$, for a given $\Lambda > 1$. 

The proof of Theorem \ref{thm:dist-reg} will be based on the two subsequent lemmas. The first lemma asserts that if a map mostly takes values outside an $\eta$-neighborhood of the convex hull of the obstacle, then the mapping is uniformly regular in the interior of the domain.

\begin{lem}\label{lem:dist-reg}
    Let $u\in W^{1,2}(B_4;\overline M)$ be a minimizing constraint map such that $\| u \|_{W^{1,2}(B_4)}\leq \Lambda$. For every $\eta > 0$, there exists a constant $\sigma \equiv \sigma(n,m,\partial M,\Lambda,\eta)\in(0,1)$, such that if $|\{\dist(u,Y) \leq \eta\}\cap B_2|\leq \sigma$, then $u \in C^{1,1}(B_1)$ with the estimate 
    \begin{equation}\label{eq:dist-reg}
    \| u \|_{C^{1,1}(B_1)} \leq c,
    \end{equation}
    where $c > 1$ depends only on $n$, $m$, $\partial M$, ${Y}$, and $\Lambda$. 
\end{lem}

\begin{proof}
In view of the $\e$-regularity theorem, Theorem \ref{thm:e-reg}, it suffices to find, under the assumption of the statement, a radius $r_0\equiv r_0(n,m,\partial M,\Lambda) > 0$, such that 
$$
E(u,x,2r_0) < \e_0^2\qquad\forall\, x_0\in B_1. 
$$
It suffices to check this claim when $Y$ is the convex hull of $M^c$, as the convex hull is the smallest (closed) convex set containing $M^c$. 

Suppose that the claim is not true for the convex hull $Y$, for some $\eta > 0$. Then, for each $k=1,2,\cdots$, there exists a minimizing constraint map $u_k \in W^{1,2}(B_4;\overline M)$ such that $\| u_k \|_{W^{1,2}(B_4)}\leq \Lambda$ and 
\begin{equation}\label{eq:distukY}
|\{\dist(u_k,Y) \leq \eta\}\cap B_2|\leq \sigma_k \to 0,
\end{equation}
but for some $x_k\in B_1$ and $r_k \to 0^+$ it holds
\begin{equation}\label{eq:Eukxk2rk}
E(u_k,x_k,2r_k) \geq \e_0^2.
\end{equation}
By Lemma \ref{lem:strong},  up to a subsequence $u_k\to u$ strongly in $W^{1,2}(B_3;\R^m)$ as $k\to\infty$, for some minimizing constraint map $u\in W^{1,2}(B_3;\overline M)$. In particular, extracting a further subsequence if necessary, we have $u_k\to u$ a.e.\ in $B_3$, so 
$$
\dist(u_k,Y)\to \dist(u,Y)\quad\text{a.e.\ in }B_3. 
$$
Passing to the limit in \eqref{eq:distukY}, we observe that
$$
\dist(u,Y) \geq \eta\quad\text{a.e.\ in }B_2. 
$$
In particular, $B_2\cap u^{-1}(\partial M) = \emptyset$, which combined with  \eqref{eq:main-sys} yields $\Delta u = 0$ in $B_2$, and therefore
\begin{equation}\label{eq:singuB2}
\Sigma(u)\cap B_2 = \emptyset.
\end{equation}
However, by \eqref{eq:Eukxk2rk}, it follows from Corollary \ref{cor:semi} that $$E(u,x_0,0^+)\geq\e_0^2 > 0,$$ for some $x_0\in\overline B_1$; by Corollary \ref{cor:e-reg} this implies $x_0\in\Sigma(u)$, a contradiction to \eqref{eq:singuB2}. This finishes the proof. 
\end{proof}

\begin{rem}
{In the setting of Lemma \ref{lem:dist-reg}, if $\sigma$ is chosen sufficiently small (for instance, if $\sigma \ll \eta^{n/2}$), then in fact $\Delta u =0$ in $B_1$. Indeed, for such $\sigma$ we must have $\dist(u,Y)>0$ in $B_1$, since otherwise by \eqref{eq:dist-reg} we would have $\dist(u,Y)\leq \eta$ in a ball of radius $c\sqrt{\eta}$, contradicting our assumption.}
\end{rem}

Our second lemma shows the uniform decay of the distance map at a point of high energy.

\begin{lem}\label{lem:dist-reg2}
    Let $u\in W^{1,2}(B_4;\overline M)$ be a minimizing constraint map, with $\| u \|_{W^{1,2}(B_4)} \leq \Lambda$. Then for every $\eta >0$, there exist radii $0 < s < r \ll 1$, both depending only on $n$, $m$, $\partial M$, $\Lambda$, and $\eta$, such that   
    $$
    E(u,0,s) \geq \e_0^2 \quad\implies\quad \dist(u,Y) \leq \eta\text{ in }B_r. 
    $$
\end{lem}

\begin{proof}
    Since $\| u \|_{L^2(B_4)} \leq \Lambda$, we have $\|\dist(u,Y)\|_{L^2(B_4)} \leq c$. By Lemma \ref{lem:dist-subharm}, $\dist(u,Y)$ is weakly subharmonic in $B_4$, therefore $\dist(u,Y) \leq c_0c$ in $B_3$ for some $c_0 = c_0(n)$, cf.\ also Corollary \ref{cor:sup}.

    Let $\eta \in (0,\frac{1}{2}c_0c)$ be given. Let $\ell$ be a positive integer, and suppose that 
    \begin{equation}\label{eq:ru}
    \lambda_\ell := \| \dist(u_\ell,Y)\|_{L^\infty(B_4)} > \eta,
    \end{equation}
    where $u_\ell(x):= u(2^{-\ell}x)$ is still  a minimizing constraint map in $W^{1,2}(B_4;\overline M)$. By the monotonicity of the normalized energy, Lemma \ref{lem:monot}, we have 
    $$E(u_\ell,0,4) = E(u,0, 2^{2-\ell}) \leq E(u,0,4) \leq 4^{2-n}\Lambda^2.$$ Furthermore, since $\dist(u,Y)\leq c_0c$ in $B_3$, $\dist(u_\ell,Y) \leq c_0c$ in $B_4$ as $2^{2-\ell} \leq 3$, for any $\ell \geq 1$.   
    
    Let $\bar c$ be the $C^{1,1}$-bound in \eqref{eq:dist-reg}, so $\bar c$ depends only on  the parameters $n$, $m$, $\partial M$ and $\max\{\Lambda^2,c_0c\}$, and choose $\sigma$ to be the measure threshold in the same lemma, determined further by $\frac{1}{2}\eta$. Since we assume $E(u,0,s) \geq \e_0^2$, we have $E(u_\ell,0,2^\ell s) \geq \e_0^2$. Hence, choosing $s$ sufficiently small so that 
    \begin{equation}\label{eq:r1}
    2^{2\ell+1} s^2|B_1|\bar c < \e_0^2, 
    \end{equation}
    there must exist a set $F\subset B_1$ of positive measure such that $|Du_\ell| > 2\bar c$ on $F$. In particular, the $C^{1,1}$-estimate \eqref{eq:dist-reg} cannot hold for $u_\ell$. Then it follows as the contraposition of Lemma \ref{lem:dist-reg} that $|\{\dist(u_\ell,Y) \leq \frac{1}{2}\eta\}\cap B_2| > \sigma$. Combining this with \eqref{eq:ru} gives us $|\{ \dist(u_\ell,Y) \leq\frac{1}{2}\lambda_\ell\}\cap B_2| > \sigma.$
    As $\dist(u_\ell,Y)$ is weakly subharmonic in $B_4$ and $\dist(u_\ell,Y) \leq \lambda_\ell$ in $B_4$, the De Giorgi oscillation lemma (see e.g.\ \cite[Proposition 2]{CV}) now implies that 
    \begin{equation}\label{eq:ru-re}
    \lambda_{\ell+2} = \| \dist(u_\ell,Y) \|_{L^\infty(B_1)} \leq (1-\mu)^2\lambda_\ell, 
    \end{equation}
    for some $\mu \in (0,\frac{1}{2})$, depending only on $n$ and $\sigma_\eta$, hence on $n$, $m$, $\partial M$, ${Y}$, $\Lambda$ and $\eta$. 

    Iterating the implication \eqref{eq:ru}--\eqref{eq:ru-re}, and recalling that $\lambda_\ell \leq \lambda_1 \leq c_0c$, we obtain 
    $$
    \lambda_\ell \leq \max\{ 4(1-\mu)^\ell c_0c, \eta\},
    $$
    for every $\ell\in\N$ for which \eqref{eq:r1} holds. Taking $\bar\ell := \min\{\ell\in\N: 4(1-\mu)^\ell c_0c \leq \eta\}$, and selecting $s$ to be sufficiently small so that \eqref{eq:r1} holds for all $\ell\leq \bar\ell$, we obtain $\lambda_{\bar\ell} \leq \eta$, that is, 
    $$
    \dist(u,Y)\leq \eta\quad\text{in }B_{4\cdot 2^{-\bar\ell}}.
    $$
    Finally, we choose $r := 4\cdot 2^{-\bar\ell}$; it is not hard to see from \eqref{eq:r1} that $r > s$. Tracking down the dependence of these radii, we also verify that both $r$ and $s$ depend only on $n$, $m$, $\partial M$, ${Y}$, $\Lambda$ and $\eta$ as desired. 
\end{proof}

We are now ready to prove the universal modulus of continuity of $\dist(u,Y)$.

\begin{proof}[Proof of Theorem \ref{thm:dist-reg}]
It suffices to prove the theorem when $\Omega = B_4$, $x_0 = 0$, $\delta = 1$. The general case follows by a simple scaling argument. 

    Let $\eta > 0$ be given and let $\lambda = \frac{1}{8}$. We fix a constant $c(n)$  to be fixed later, and we choose radii $0<s<r\ll 1$ as in Lemma \ref{lem:dist-reg2}, corresponding to the parameters $n$, $m$, $\partial M$, $c(n)\Lambda$, and $\eta$.     
    
    If $E(u,0,2\lambda s) \geq \e_0^2$, then Lemma \ref{lem:dist-reg2}, applied to $u_{2\lambda}(x) := u(2\lambda x)$, shows that $0\leq\dist(u,Y) \leq \eta$ in $B_{2\lambda r}$, whence 
    \begin{equation}\label{eq:dist-reg-case1} \osc_{B_{2\lambda r}} \dist(u,Y)\leq \eta.
    \end{equation}
    This application of Lemma \ref{lem:dist-reg2} determines the choice of $c(n)$.

    Alternatively, if $E(u,0,2\lambda s) < \e_0^2$, then by the $\e$-regularity theorem, Theorem \ref{thm:e-reg}, we have $u\in C^{1,1}(B_{\lambda s})$ and $|Du|\leq \frac{c_0}{\lambda s}$ in $B_{\lambda s}$, for some constant $c_0\equiv c_0(n,m,\partial M)$. Combined with the fact that $\dist(\cdot,Y)$ is 1-Lipschitz, we find that
    \begin{equation}\label{eq:dist-reg-case2}
    \osc_{B_{\eta \lambda s/(2c_0)}} \dist(u,Y)\leq \eta.
    \end{equation}
    By virtue of \eqref{eq:dist-reg-case1} and \eqref{eq:dist-reg-case2}, we derive in either case that the oscillation of $\dist(u,Y)$ is bounded by $\eta$ in $B_t(x_0)$, for some radius $t \equiv t(n,m,\partial M,\Lambda,\eta)$. This finishes the proof. 
\end{proof}

A direct consequence of Theorem \ref{thm:dist-reg} is that, when the obstacle $M^c$ is convex, there exists a universal continuity estimate on $\rho\circ u$, thus improving \cite[Theorem 6.1]{FGKS}.

\begin{cor}\label{cor:dist}
    Let $M\subset \R^m$ be a smooth domain with compact convex complement, and let $u\in W^{1,2}(\Omega;\overline M)$ be a minimizing constraint map with $\| u \|_{W^{1,2}(\Omega)}\leq\Lambda$. Then $\rho\circ u$ is subharmonic and continuous in $\Omega$. In addition, its modulus of continuity depends on $n$, $m$, $\partial M$, $\Lambda$, and $\delta$ only.
\end{cor}


\section{Universal \texorpdfstring{$L^{-\gamma}$}{L(-epsilon)}-estimate near discontinuities}\label{sec:Le}

In this section we establish a universal $L^{-\gamma}$-estimate of the energy density of minimizing constraint maps near {\it effective} discontinuity points,\footnote{By {\it effective} discontinuity points, we mean the points which have normalized energy above the $\e$-regularity threshold. The set of effective discontinuity points necessarily includes a neighborhood of the set $\Sigma(u)$ of the {\it actual} discontinuity points, due to  the monotonicity of the energy. For further discussion, we refer the reader to \cite{NV1}.} see Theorem \ref{thm:Le}.
As mentioned before, this result also applies to classical harmonic maps, and we believe it to be of its own interest.

Throughout this section, we shall fix $\e_0 \equiv \e_0 (n,m,\partial M)$ as the energy threshold provided by Corollary \ref{cor:e-reg}. Let us begin our analysis by finding the critical scale, labeled  $r_u (x)$ below.

\begin{lem}[Critical scale]
\label{lem:crit}
Let  $u\in W^{1,2}(B_4;\overline M)$ be a minimizing constraint map satisfying $\| u \|_{W^{1,2}(B_4)}\leq \Lambda$. 
There exists a large integer $\ell_0 > 1$, depending only on $n$, $m$, $\partial M$, and $\Lambda$, such that the following holds:  setting 
\begin{equation}\label{eq:crit}
r_u(x) := \begin{cases}
\sup\left\{ r\in (0,2]: E(u,x,r)\leq \ell_0^{-2}\e_0^2\right\} & \text{if }x\in B_2\setminus \Sigma(u),\\
0 &  \text{if }x\in B_2\cap \Sigma(u),
\end{cases}
\end{equation}
we have that $\Sigma(u) = r_u^{-1}(0)$ and that if $E(u,y,1)\geq\e_0^2$ for some $y\in B_2$, then $r_u(x) \leq 1$ for all $x\in B_2$.\end{lem}

\begin{proof}
The assertion $\Sigma(u) = r_u^{-1}(0)$ follows immediately from Corollary~\ref{cor:e-reg}, as $\e_0$ is chosen as the energy threshold there.

Assume now that $E(u,y,1)\geq\e_0^2$ for some $y\in B_2$. To prove that $r_u(x)\leq 1$ for all $x\in B_2$, we claim that 
$$
E(u,x,1) >\frac{1}{\ell_0^2}\e_0^2,\quad\forall \,x\in B_2,
$$
for some positive integer $\ell_0 \equiv \ell_0(n,m,\partial M,\Lambda)$. Once the claim holds, the conclusion follows from the monotonicity of the normalized energy, Lemma \ref{lem:monot}. 

Suppose that the claim does not hold. Then there is, for every $k\in\N$, a minimizing constraint map $u_k \in W^{1,2}(B_4;\overline M)$ with $\| u_k \|_{W^{1,2}(B_4)} \leq \Lambda$ such that $E(u_k,y_k,1)\geq \e_0^2$ for some $y_k \in B_2$, but 
$$
E(u_k,x_k,1) \leq \frac{1}{k^2} \e_0^2,
$$
for some $x_k \in B_2$. We may assume without loss of generality that $x_k\to x_\infty$ and $y_k\to y_\infty$ as $k\to\infty$, for some $x_0,y_0\in\overline B_2$.
Up to a subsequence, the compactness of minimizing constraint maps (Lemma \ref{lem:strong}) implies that $u_k$ converges strongly to $u$ in $W^{1,2}(B_3;\overline{M})$, where $u$ is a minimizing constraint map. Given that $B_1(y_k) \cup B_1(y_\infty) \subset B_3$, this strong convergence implies that 
\begin{equation}\label{eq:bound}
E(u,y_\infty,1) \geq \e_0^2 > 0. 
\end{equation}
However, since $B_1(x_k)\cup B_1(x_\infty)\subset B_3$, we also have
$$
E(u,x_0,1) = 0, 
$$
thus $|Du| = 0$ a.e.\ in $B_1(x_0)$. 

Note now that, in view of \eqref{eq:main-sys}, we have $|\Delta u|\leq c|Du|^2$. Also, $u$ is locally Lipschitz in $B_3\setminus \Sigma(u)$. Furthermore, Theorem \ref{thm:dim-sing} implies that $\dim_\cH\Sigma(u)\leq n-3 < n-1$, and so by \cite[Theorem IV.4]{HuWa} the open set $B_3\setminus\Sigma(u)$ is connected. Thus, the classical unique continuation principle (see \cite[Remark 3]{A}) yields $|Du| = 0$ everywhere outside $\Sigma(u)$, whence also $|Du|=0$ a.e.\  in $B_3$; in particular $E(u,y_0,1) = 0$,  contradicting \eqref{eq:bound}. This finishes the proof.
\end{proof}

In what follows, we shall also fix the large constant $\ell_0 \equiv \ell_0(n,m,\partial M,\Lambda)$ chosen in Lemma \ref{lem:crit}. Given a smooth non-constant map $f\colon \Omega\to \R^m$ and a ball $B_r(x_0)\subset\Omega$, let us write $N(f,x_0,r)$ for the \textit{frequency} of $f$ in the ball $B_r(x_0)$, i.e.,
 \begin{equation}\label{freq}
      N(f,x_0,r) := \frac{r\int_{B_r(x_0)} |Df|^2\,dx}{\int_{\partial B_r(x_0)} |f - (f)_{x_0,r}|^2 \haus},
  \end{equation}
 where $(f)_{x_0,r}$ denotes the integral average of $f$ over the $(n-1)$-dimensional sphere $\partial B_r(x_0)$. When $x_0 = 0$, we shall write $(f)_r \equiv (f)_{0,r}$ for simplicity. 

We are ready to make the key observation of this section: the frequency is uniformly bounded below the critical scale, provided that the mapping has discontinuities nearby (or, more generally, points of sufficiently large energy).
 
 \begin{lem}[Frequency bound at the critical scale]
\label{lem:freq}
Let $u\in W^{1,2}(B_4;\overline M)$ be a minimizing constraint map satisfying $\| u \|_{W^{1,2}(B_4)} \leq \Lambda$. There exist a frequency threshold $\lambda_0 > 1$, depending only on $n$, $m$, $\partial M$, and $\Lambda$, such that if $E(u,y,1)\geq \e_0^2$ for some $y\in B_2$, then 
     $$
     N(u,x,r_u(x)) \leq \lambda_0,
     $$
     for any $x\in B_1\setminus\Sigma(u)$, where $r_u(x)$ is as in \eqref{eq:crit}.
 \end{lem}

 \begin{proof} 
Throughout the proof, $c$ will denote a generic constant depending at most on $n$, $m$ and $\partial M$. Assume that the conclusion of the lemma is false. Then for every $k \in \N$, there must exist a minimizing constraint map $u_k \in  W^{1,2}(B_4;\overline M)$ such that $\| u_k \|_{W^{1,2}(B_4)} \leq \Lambda$, $E(u_k,y_k,1)\geq\e_0^2$ for some $y_k\in B_2$, but at some $x_k \in B_1\setminus\Sigma(u_k)$ it satisfies, with $r_k := r_{u_k}(x_k)$,
     \begin{equation}\label{eq:freq-false}
     N(u_k,x_k,r_k) \geq k. 
     \end{equation}    
By Lemma \ref{lem:crit}, $r_k \leq 1$ for all $k$. This, along with $x_k\in B_2$, yields 
     \begin{equation}\label{eq:rk-re}
     B_{\frac{3}{2}r_k}(x_k)\subset B_{\frac{3}{2}}(x_k)\subset B_4. 
     \end{equation}
In what follows, let us write  
     \begin{equation}\label{eq:vk}
     v_k (y) := u_k\bigg(\frac{1}{2}r_ky + x_k\bigg),
     \end{equation}
 which is a minimizing constraint map in $W^{1,2}(B_3;\overline M)$. By Lemma \ref{lem:monot}, \eqref{eq:rk-re} and the a priori bound $\int_{B_1}|Du_k|^2\,dx \leq \Lambda^2$, we deduce that
     \begin{equation}\label{eq:Evk}
     \begin{aligned}
     E(v_k,0,3) = E\bigg(u_k,x_k,\frac{3}{2}r_k\bigg)\leq E\bigg(u_k,x_k,\frac{3}{2}\bigg) \leq c\int_{B_4} |Du_k|^2\,dx \leq c\Lambda^2.
     \end{aligned}
     \end{equation}
In addition, since $M^c$ is compact and $\| u_k \|_{L^{2}(B_4)} \leq \Lambda$, Corollary \ref{cor:sup} implies that $|u_k|\leq c (1+\Lambda)$ a.e.\ in $B_{7/2}$, thus
 $$
 |v_k|\leq c(1 + \Lambda)\quad\text{a.e.\ in }B_{5/2}.
 $$
In particular $\{v_k\}_{k=1}^\infty$ is bounded in $W^{1,2}(B_{5/2};\overline M)$, so Lemma \ref{lem:strong} yields a minimizing constraint map $v \in W^{1,2}(B_2;\overline M)$ such that 
     \begin{equation}\label{eq:vk-v0-W12} 
     v_k\to v\quad\text{strongly in }W^{1,2}(B_2;\R^m),
     \end{equation}
after potentially extracting a subsequence. 

However, as \eqref{eq:freq-false} yields $N(v_k,0,2) = N(u_k,x_k,r_k) \geq k$, it follows from \eqref{eq:Evk} that
     \begin{equation}\label{eq:Pivk-L2}
     \int_{\partial B_2} |v_k - (v_k)_1|^2\,d\cH^{n-1} \leq \frac{c}{k}\Lambda^2,
     \end{equation}
for every $k\in\N$ .
By \eqref{eq:vk-v0-W12}, the boundedness of the trace operator \cite[\S3 Theorem 4.6]{EG}   implies 
     $$
     v_k \to v \quad\text{strongly in }L^2(\partial B_2;\R^m),
     $$
and so we find, via \eqref{eq:Pivk-L2}, that  
     \begin{equation}\label{eq:V0}
     v  = a\quad\text{a.e.\ on }\partial B_2,
     \end{equation}
for some point $a\in \partial M$. 

We now observe that \eqref{eq:V0} implies $v = a$ a.e.\ in $B_3$: taking $\vp : B_3\to \R^m$ defined by $\vp = a$ in $B_2$ and $\vp = v$ in $B_3\setminus B_2$, $\vp$ is an admissible map and so the minimality of $v$ implies 
$$
\int_{B_2}|Dv|^2\,dx \leq \int_{B_2} |D\vp|^2\,dx = 0, 
$$
i.e.\ $|Dv| = 0$ a.e.\ in $B_2$. As in the proof of Lemma \ref{lem:crit}, the  unique continuation principle \cite[Remark 3]{A} allows us to deduce that $v = a$ in 
$B_3\setminus\Sigma(v)$, which then readily implies that $v = a$ in $B_3$ and $\Sigma(v) = \emptyset$. This is however a contradiction, as the choice of $r_k$ and \eqref{eq:vk-v0-W12} yield
$$
E(v,0,2) = \lim_{k\to\infty} E(v_k,0,2) = \lim_{k\to\infty} E(u_k,x_k,r_k) \geq \frac{1}{\ell_0^2}\e_0^2. 
$$
This finishes the proof. 
 \end{proof}

Next we show that, above the critical scale, the differential of the map has a uniform lower bound on a universal fraction of each cube.

 \begin{lem}[Above the critical scale]
 \label{lem:above}
For any $\sigma\in (0,1)$, there are a small constant $\delta > 0$, depending only on $n$, $m$, $\partial M$, $\Lambda$, and $\sigma$, and a large constant $c_0 > 1$, depending solely on $n$, $m$, and $\partial M$, such that for any minimizing constraint map $u\in W^{1,2}(B_8;\overline M)$ satisfying $\|u \|_{W^{1,2}(B_8)}\leq\Lambda$, if $E(u,0,4) \geq \frac{1}{\ell_0^2}\e_0^2$, then
     $$
     |\{|Du| \geq \delta |Q|^{-\frac{1}{n}}\}\cap Q| \geq (1-\sigma)|Q|,
     $$
where $Q=\left[-1/\sqrt{n},1/\sqrt{n}\right]^n$ is the largest cube contained in $B_1$.
\end{lem} 

 \begin{proof}

Fix $\sigma\in(0,1)$. Assume, for the sake of contradiction, that for each $k\in\N$, there exists a minimizing constraint map $u_k \in W^{1,2}(B_8;\overline M)$ satisfying $\|u_k\|_{W^{1,2}(B_8)} \leq \Lambda$ and $E(u_k,0,4) \geq\frac{1}{\ell_0^2}{\e_0^2}$, but 
 \begin{equation}\label{eq:alt-false}
 |\{ |Du_k| >\delta_k |Q|^{-\frac{1}{n}} \} \cap Q| < (1-\sigma)|Q|,
 \end{equation}
 along a sequence $\delta_k\to 0$.

 As in the proof of Lemma \ref{lem:freq}, one can find, via Lemma \ref{lem:strong}, a minimizing constraint map $u \in W^{1,2}(B_4;\overline M)$ such that 
\begin{equation}\label{eq:uk-u-W12}
u_k \to u \quad\text{strongly in }W^{1,2}(B_4;\R^m),
\end{equation}
along a subsequence. To simplify our notation, let us assume without loss of generality that the convergence holds along the full sequence. 
Therefore, by Egorov's theorem, there exists a closed set $F\subset Q$, with 
\begin{equation}\label{eq:F}
|F| < \frac{1}{2}\sigma|Q|,
\end{equation}
such that 
\begin{equation}\label{eq:DPivk-DPv}
Du_k \to Du \quad\text{uniformly on }Q\setminus F.
\end{equation}
Thus, by \eqref{eq:alt-false} and \eqref{eq:F},
\begin{equation}\label{eq:DPivk-msr}
|\{ |Du_k| \leq \delta_k |Q|^{-\frac{1}{n}}\}\cap (Q\setminus F)| \geq \sigma |Q| - |F| \geq \frac{1}{2} \sigma |Q|. 
\end{equation}
Recalling that $\delta_k\to 0$, we obtain from \eqref{eq:DPivk-DPv} and \eqref{eq:DPivk-msr} that 
\begin{equation*}\label{eq:DPiv-msr}
 |\{ |Du| = 0 \}\cap (Q\setminus F)\}| \geq \frac{1}{2}\sigma |Q|. 
\end{equation*}
Arguing as in the proof of Lemma \ref{lem:freq}, we again deduce from the minimality of $u$ and the unique continuation principle that $|Du| = 0$ in $B_4$. However, by \eqref{eq:uk-u-W12},
$$
    E(u,0,4) = \lim_{k\to\infty} E(u_k,0,4) \geq \frac{1}{\ell_0^2}\e_0^2,
$$
which is a contradiction. This completes the proof.

 \end{proof}

Next, we look below the critical scale $r_u(x) $. There, our system \eqref{eq:main-sys} can be considered as a linear system with  bounded coefficients for the first-order term, due to the $\e$-regularity theorem. Therefore, the energy density $Du$ admits a uniform doubling property depending on the frequency of the mapping $u$.  For the reader's convenience, we recall these standard arguments in Lemma \ref{lem:freq-monot} in the appendix.

 \begin{lem}[Below the critical scale]
 \label{lem:below}
Let $\lambda > 1$. Then there are constants $\delta = \delta(n,m,\partial M,\lambda)$ and $\sigma = \sigma(n,m,\partial M,\lambda)$, $0 < \delta,\sigma < 1$, such that if $u\in W^{1,2}(B_4;\overline M)$ is a minimizing constraint map satisfying $E(u,0,4) \leq \e_0^2$ and $N(u,0,4) \leq \lambda$, then 
 $$
 |\{ |Du| \geq \delta \| Du\|_{L^\infty(3Q)} \}\cap Q|\geq (1-\sigma)|Q|,
 $$
 for every cube $Q$ such that $3Q\subset B_2$. 
 \end{lem}

 \begin{proof}
By our assumptions and Theorem \ref{thm:e-reg}, we deduce from \eqref{eq:main-sys} that 
 \begin{equation}\label{eq:V-pde-re}
 \Delta u^i = a_{ij}^\alpha D_\alpha u^j\quad\text{in }B_2,
 \end{equation}
with $\| a_{ij}^\alpha \|_{L^\infty(B_2)} \leq c$, where $c$ depends only on $n$, $m$ and $\partial M$. Hence, the assertion follows immediately by applying Lemma \ref{lem:freq-monot} with  $f = u$. 
 \end{proof}

 Combining the last three lemmas, we arrive at the following conclusion.

 \begin{lem}
 \label{lem:gen}
There are constants $\delta_0,\sigma_0 \in (0,1)$, both depending only on $n$, $m$, $\partial M$, and $\Lambda$, such that for any minimizing constraint map $u\in W^{1,2}(B_4;\overline M)$ satisfying $\| u \|_{W^{1,2}(B_4)} \leq \Lambda$, if $E(u,y,1)\geq \e_0^2$ for some $y\in B_2$, then 
 $$
 |\{ |Du| \geq \delta_0 \min\{|Q|^{-\frac{1}{n}}, \| Du\|_{L^\infty(3Q)}\} \}\cap Q| \geq (1-\sigma_0)|Q|, $$
 for every cube $Q$ centered in $\overline B_2$ such that $\diam Q\leq\frac{1}{8}$. 
 \end{lem}

 \begin{proof}
We may assume without loss of generality that $\Lambda > 1$. Let $\lambda_0$ be as in Lemma \ref{lem:freq}, and choose $\delta_1$ and $\sigma_0$ as in Lemma \ref{lem:below} corresponding to $\lambda =\lambda_0$. We then select $\delta_2$ as in Lemma \ref{lem:above} corresponding to the chosen $\sigma_0$ (and with $c\Lambda$ in place of $\Lambda$ for some $c\equiv c(n,m,\partial M)$). Finally, we set $\delta_0 := \min\{\delta_1,\delta_2\} > 0$. 

Let $Q$ be a cube centered in $\overline B_2$ with $\diam Q\leq\frac{1}{8}$. Denote by $s$ its sidelength and by $x_0$ its center, so that $Q = Q_s(x_0)$, $\sqrt{n}s \leq \frac{1}{8}$, and $|x_0|\leq 2$. Let $B_r(x_0)$ be the smallest ball containing $Q$, i.e.\ $r=\frac{1}{2}\sqrt{n}s \leq \frac{1}{16}$. Then $B_{8r}(x_0)\subset B_1(x_0)\subset B_3$. Hence, the monotonicity of the normalized energy, Lemma \ref{lem:monot}, along with our assumption $\| u \|_{W^{1,2}(B_4)}\leq \Lambda$, yields  
 \begin{equation}\label{eq:Eux0r}
E(u,x_0,8r)\leq E(u,x_0,2) \leq 2^{2-n}\|Du\|_{L^2(B_4)}^2 \leq 2^{2-n}\Lambda^2.  
 \end{equation}
Moreover, by the local $L^\infty$-estimate, Corollary \ref{cor:sup}, we have 
\begin{equation}\label{eq:supux0r}
\|u \|_{L^\infty(B_{8r}(x_0))} \leq \| u\|_{L^\infty(B_3)} \leq c + c \|u\|_{L^2(B_4)} \leq c\Lambda;
\end{equation}
in the last inequality, we also used $\Lambda > 1$. 
Recall the critical scale $r_u(x_0)$ defined in \eqref{eq:crit}. By Theorem \ref{thm:e-reg} and Lemma \ref{lem:crit}, it holds $0 \leq r_u(x_0) \leq 1$.  
 Let us divide the proof into two cases.
 \setcounter{case}{0}

 \begin{case}
 $r \geq \frac{1}{4}r_u(x_0)$.
 \end{case}
In view of \eqref{eq:crit}, we have $E(u,x_0,4r) \geq \frac{1}{\ell_0^2}{\e_0^2}$. Also, by virtue of \eqref{eq:Eux0r} and \eqref{eq:supux0r}, we have 
$$
\| u_{x_0,r} \|_{W^{1,2}(B_8)} \leq c\Lambda,
$$
where $u_{x_0,r}(y) := u(ry + x_0)$ with $c\equiv c(n,m,\partial M)$. Hence, Lemma \ref{lem:above} (with $\sigma = \sigma_0$ and $\delta = \delta_2$ there) applies to $u_{x_0,r}$. Rescaling back and recalling that $Q$ is the largest cube contained in $B_r(x_0)$, since $\delta_0 \leq \delta_1$ we obtain  
 $$
 |\{ |Du| \geq \delta_0 |Q|^{-\frac{1}{n}}\}\cap Q| \geq (1-\sigma_0)|Q|.
 $$
 Thus, the conclusion of the lemma follows.

 \begin{case}
 $r < \frac{1}{4}r_u(x_0)$.
 \end{case} 
In this case, we have $x_0\not \in \Sigma(u)$ and $E(u,x_0,4r) < \frac{1}{\ell_0^2}\e_0^2$. By Lemma \ref{lem:freq} we have $N(u,x_0,r_u(x_0)) \leq \lambda_0$, and by definition $E(u,x_0,r_u(x_0)) \leq \frac{1}{\ell_0^2}\e_0^2 < \e_0^2$. Thus Lemma \ref{lem:below} applies to $u_{x_0,r_u(x_0)} (y) := u(r_u(x_0)y + x_0)$. Since $4r \leq r_u(x_0)$, and $Q\subset B_r(x_0)$, we have $3Q\subset B_{3r}(x_0)\subset B_{r_u(x_0)}(x_0)$. Thus, we deduce (after rescaling back) that 
 $$
 |\{ |Du| \geq \delta_0\|Du\|_{L^\infty(3Q)}\}\cap Q| \geq (1-\sigma_0)|Q|. 
 $$
 Therefore, the conclusion of this lemma is verified for the alternative case as well. This finishes the proof. 
 \end{proof}

 We are now ready to prove Theorem \ref{thm:Le}, which is a direct consequence of the last lemma and the Calder\'on--Zygmund cube decomposition, whose statement we recall here for the reader's convenience (see e.g.\ \cite{CC} for a proof).

 \begin{lem}[Calder\'on-Zygmund]\label{lem:CZ-re}
 Let $A,B\subset Q_1$ be measurable sets, and let $\sigma > 0$ be given. Assume the following:  
 \begin{enumerate}
 \item $|A| \leq \sigma$.
 \item If $|A\cap Q| > \sigma |Q|$ for some dyadic cube $Q$, then $\tilde Q\subset B$, where $\tilde Q$ is the predecessor of $Q$. 
 \end{enumerate}
 Then $|A| \leq \sigma |B|$. 
 \end{lem}

\begin{proof}[Proof of Theorem \ref{thm:Le}]
Consider a collection of disjoint open cubes $Q^i$, centered in $\overline B_2$, with the same side-length $s$, i.e.\  $n s^2 = (\diam Q^i)^2$, such that $B_1\subset \bigcup_{i=1}^{i_0} \overline Q^i$ and $\diam Q^i\leq\frac{1}{8}$ for each $i$. We can take the collection in such a way that its cardinality $i_0$ is bounded by a constant depending only on $n$. By the assumption $\diam Q^i\leq\frac{1}{8}$, we must have $64ns^2 \leq 1$. 

Take $\e_0$, $\ell_0$, $\delta_0$, and $\sigma_0$ as in Corollary \ref{cor:e-reg}, Lemma \ref{lem:crit} and Lemma \ref{lem:gen}. We claim that 
 \begin{equation}\label{eq:CZ-claim}
 |\{|Du| \leq t\}\cap Q^i| \leq \min\{ct^{2\gamma},1\}s^n,
 \end{equation}
 for every $t > 0$, for some $c > 1$ and $\gamma > 0$, both depending only on $n$, $m$, $\partial M$, and $\Lambda$, but not on $i$. Then by \eqref{eq:CZ-claim} and elementary real analysis (integrating over all volume of super-level sets), we deduce that 
 $$
 \begin{aligned}
\int_{B_1} |Du|^{-\gamma}\,dx &\leq \sum_{i=1}^{i_0} \int_{Q^i} |Du|^{-\gamma}\,dx = \gamma \sum_{i=1}^{i_0} \int_0^\infty \tau^{-1 + \gamma} |\{ |Du| \leq \tau^{-1}\}\cap Q^i|\,d\tau \\
& \leq \gamma i_0 s^n \bigg[\int_0^1 \tau^{-1+\gamma}\,d\tau + \int_1^\infty c\tau^{-1-\gamma}\,d\tau \bigg] \leq c,
 \end{aligned}
 $$
where the last inequality follows from our choice (specified above) of both $i_0$ and $s$, in such a way that they are bounded above by a constant depending only on $n$.

To prove \eqref{eq:CZ-claim},  let us define, for each $k \in \mathbb N$, the set
 $$
 A_k^i := \{ |Du| \leq \delta_0^k\}\cap Q^i. 
 $$
 It is clear that $A_{k+1}^i\subset A_k^i$. Moreover, Lemma \ref{lem:gen}
 implies that $|A_1^i|\leq \sigma s^n$, therefore 
 \begin{equation}
\label{eq:CZ1}
 |A_{k}^i|\leq \sigma s^n, \quad \forall\,k \in \mathbb N.  
\end{equation}
 Next we claim that if $Q$ is a member of the dyadic cube decomposition of $Q^i$ then 
\begin{equation}
\label{eq:CZ2}
|A_{k+1}^i\cap Q| > \sigma_0 |Q| \quad \implies\quad \tilde Q\subset A_k^i,
\end{equation} where $\tilde Q$ is the predecessor of $Q$.
 
To prove this, assume by contradiction that $\tilde Q\setminus A_k^i\neq\emptyset$. By the definition of $A_k^i$, we have 
 $$
 \|D u\|_{L^\infty(\tilde Q)} > \delta_0^k. 
 $$
 Since $\tilde Q\subset 3Q \subset 3Q^i\subset B_2$ and $|Q|^{-1/n} \geq |Q^i|^{-1/n} = s^{-1} > \sqrt n > 1 > \delta_0^{k+1}$, Lemma \ref{lem:gen} yields that 
 $$
 |\{ |Du| \geq \delta_0^{k+1}\}\cap Q| \geq (1-\sigma_0)|Q|.
 $$
This implies that $|A_{k+1}^i\cap Q| \leq \sigma_0 |Q|$, a contradiction that proves the result.

Thanks to \eqref{eq:CZ1} and \eqref{eq:CZ2}, it follows from Lemma~\ref{lem:CZ-re} that $|A_{k+1}^i| \leq \sigma_0 |A_k^i|$. Since $k$ was an arbitrary positive integer and $|A_1^i|\leq \sigma_0s^n$, we  conclude that 
 $$
 |A_k^i| \leq \sigma_0^k s^n=s^n (\delta_0^k)^{2\gamma},\quad\forall\, k\in\N,
 $$
 with $\gamma:=\frac{\log\sigma_0}{2\log\delta_0}$.
Due to this bound, \eqref{eq:CZ-claim} follows. Since both $\delta_0$ and $\sigma_0$ depend only on $n$, $m$, $\Lambda$, and $\partial M$, so do the constants $\gamma$ and $c$. This completes the proof. 
 \end{proof}

 
\section{Regularity near free boundaries}\label{sec:reg}

The purpose of this section is to prove Theorem \ref{thm:uni-reg}, which asserts that,  whenever the obstacle is uniformly convex, minimizing constraint maps are $C^{1,1}$ in a universal neighborhood of the non-coincidence set. Therefore, throughout this section, we assume that $M$ is a smooth domain with uniformly convex complement. Moreover, we fix $\e_0$ as the energy threshold chosen as in Lemma~\ref{lem:crit}, which depends only on $n$, $m$, and $\partial M$.

As we shall see, the result follows from Theorem \ref{thm:Le}, proved in the previous section,  with the following elementary   De Giorgi-type lemma, which we will apply to the subharmonic function $\rho \circ u$.

 \begin{lem}\label{lem:degiorgi}
 Let $w\in W^{1,2}(B_2)$ be a nonnegative subharmonic function, and suppose that there are positive constants $a$, $p$, $q$ such that 
 \begin{equation}\label{eq:decay}
 |\{ w > 0\}\cap B_r|\leq \frac{a}{(R-r)^q}  \| w \|_{L^\infty(B_R)}^p,
 \end{equation}
 whenever $1\leq r <R \leq 2$. Then there is a constant $\eta\in(0,1)$, depending only on $n$, $a$, $p$, and $q$, such that if $w \leq \eta$ a.e.\ in $B_2$, then $w = 0$ a.e.\ in $B_1$. 
 \end{lem} 

 \begin{proof}
 Let $r_k := 1 + 4^{-k}$ and $s_k := 1 + 2\cdot 4^{-k-1}$, for each $k\in \N$. Then 
 $r_k - s_k = 2\cdot 4^{-k-1}$ and $s_k - r_{k+1} = 4^{-k-1}$. Denote $\delta_k := |\{ w > 0\}\cap B_{s_k}|$ and $\lambda_k := \|w\|_{L^\infty(B_{r_k})}$. Then our hypothesis \eqref{eq:decay} yields
 $$\delta_k \leq a (2\cdot 4^{k+1})^q \lambda_k^p.$$
This, combined with the local $L^\infty$-estimate for subharmonic functions, implies that 
 $$
 \lambda_{k+1} \leq \frac{c}{(s_k-r_{k+1})^n}\|w\|_{L^2(B_{s_k})} \leq c\cdot 4^{n(k+1)} \lambda_k \delta_k \leq  C a \cdot 4^{ k(n+q)} \lambda_k^{1+p},
 $$
 where $c= c(n)$ and $C=2^{3 \beta + 2n} c$. Thanks to this recursive relation, it follows that
 $$
 \lambda_0 \leq \eta \quad\Longrightarrow \quad \lim_{k\to\infty} \lambda_k = 0,
 $$
 if $\eta$ is sufficiently small; the smallness condition can be determined {\it a priori} by $n$, $a$ and $\gamma$ only. Since $\| w \|_{L^\infty(B_1)} \leq \lambda_k$ for all $k$, we conclude that $w = 0$ a.e.\ in $B_1$.
 \end{proof}

We can now prove the following key step towards Theorem \ref{thm:uni-reg}. 

\begin{lem}\label{lem:reg}
There exists a constant $\eta>0$, depending only on $n$, $m$, $\partial M$, and $\Lambda$, such that for any minimizing constraint map $u \in W^{1,2}(B_4;\overline M)$ satisfying $\int_{B_4} |Du|^2\,dx \leq \Lambda^2$ and $\rho\circ u \leq \eta$ in $B_4$, 
$$
\text{$E(u,y,1)\geq\e_0^2$ for some $y\in B_2$}\quad\implies\quad B_1\subset u^{-1}(\partial M).
$$
\end{lem} 

 \begin{proof}
Fix a pair of radii with $1 \leq r<R \leq  2$, and let $\psi\in C_0^\infty(B_R)$ be a smooth cutoff function such that $\psi = 1$ in $B_r$, $|D^2\psi|\leq c/(R-r)^2$, and $0\leq \psi\leq 1$ in $B_R$, with $c = c(n)$. In view of the equation \eqref{eq:dist-pde} and by the geometric bounds \eqref{eq:Hessrho} with $\xi = D_\alpha u(x)$ (note that, by chain rule, $\nabla\Pi(u(x)) \cdot D_\alpha u(x) = D_\alpha (\Pi\circ u)(x)$), we find
 \begin{equation}\label{eq:key11}
 \begin{aligned}
 \frac{c}{(R-r)^2}\int_{B_R} (\rho\circ u)\,dx &\geq \int_{B_R} (\rho\circ u)\Delta\psi\,dx \\
 &= \int_{B_R\cap u^{-1}(M)} \psi \Hess\rho_u(Du,Du)\,dx \\
 &\geq c_0\int_{B_r\cap u^{-1}(M)} |D(\Pi\circ u)|^2\, dx.
 \end{aligned}
 \end{equation}
On the other hand, as $\rho\circ u$ is weakly subharmonic in $B_4$ due to the convexity of $M^c$, the Caccioppoli inequality along with $\rho\circ u \leq \eta < 1$ in $B_4$ and $|D(\rho\circ u)| = 0$ a.e.\ in $B_r\setminus u^{-1}(M)$ shows that 
\begin{equation}\label{eq:key12}
\frac{c}{(R-r)^2}\int_{B_R} (\rho\circ u)\,dx \geq \int_{B_r\cap u^{-1}(M)} |D(\rho\circ u)|^2\,dx.
\end{equation}
Adding \eqref{eq:key11} and \eqref{eq:key12}, and utilizing $|Du|^2\leq c|D(\Pi\circ u)|^2 + c|D(\rho\circ u)|^2$ a.e., we get
\begin{equation}\label{eq:key1}
\frac{c}{(R-r)^2}\int_{B_R} (\rho\circ u)\,dx \geq c_0 \int_{B_r\cap u^{-1}(M)} |Du|^2\,dx,
\end{equation}
with a possibly larger constant $c>1$. 

Let $\e ,\eta >0$ and $c_1 > 1$ be the constants determined {\it a priori} by $n$, $m$, $\partial M$, and $\Lambda$, as in Theorem \ref{thm:Le}. Since by assumption $E(u,y,1)\geq\e_0^2$ for some $y\in B_2$ and $\rho\circ u\leq \eta$ in $B_4$, we invoke Theorem \ref{thm:Le} to obtain 
 \begin{equation}\label{eq:key2}
 \int_{B_2} |Du|^{-\gamma}\,dx \leq c.
 \end{equation}
 Choosing $p > 1$ such that $\frac{1}{p} = \frac{1}{2} + \frac 1 \gamma$, H\"older inequality and \eqref{eq:key2} yield
 \begin{equation}\label{eq:key3}
 |B_r\cap u^{-1}(M)|^{\frac{2}{p}}\leq \biggl(\int_{B_r\cap u^{-1}(M)} |Du|^{-\gamma}\,dx\biggr)^{\frac{2}{\gamma}}\int_{B_r\cap u^{-1}(M)} |Du|^2\,dx \leq c^{\frac{2}{\gamma}}\int_{B_r\cap u^{-1}(M)} |Du|^2\,dx.   
 \end{equation}
 By \eqref{eq:key1} and \eqref{eq:key3}, this implies
 \begin{equation}\label{eq:key4}
|B_r\cap u^{-1}(M)| \leq \frac{c}{(R-r)^{p}} \| \rho\circ u\|_{L^\infty(B_R)}^{\frac{p}{2}};
 \end{equation}
 note that $p = \frac{2\gamma}{2+\gamma}$ depends only on $n$, $m$, $\partial M$, and $\Lambda$. Since $u^{-1}(M) = \{\rho\circ u > 0\}$, condition  \eqref{eq:decay} is verified. By choosing $\eta$ even smaller if necessary,  Lemma \ref{lem:degiorgi} applies, so from the assumption $\rho\circ u \leq \eta$ in $B_4$ we conclude that $\rho\circ u = 0$ in $B_1$. In other words,  $B_1\subset u^{-1}(\partial M)$, as desired. 
 \end{proof}

 With the help of Lemma \ref{lem:reg}  we establish the following rigidity theorem, which may be of independent interest.

 \begin{thm}\label{thm:rigid}
 Let $u\in W_\loc^{1,2}\cap L^\infty(\R^n;\overline M)$ be a locally minimizing constraint map. Suppose that \begin{equation}\label{eq:rigid}
 \limsup_{k\to\infty}  \left[ k^{2-n} \int_{B_k}  |Du|^2 \,dx \right] < + \infty\quad\text{and}\quad 0\in\Sigma(u) .
 \end{equation}
 Then $u \in W_\loc^{1,2}(\R^n;\partial M)$ and in particular $u$ is a locally minimizing harmonic map into $\partial M$. 
 \end{thm}

 \begin{proof} 
By Lemma \ref{lem:monot}, Corollary \ref{cor:e-reg}, $0\in\Sigma(u)$, and \eqref{eq:rigid}, there is a constant $\theta_0 > \e_0^2$ such that 
 \begin{equation}\label{eq:rigid-re}
 \lim_{k\to\infty} \left[ (kR)^{2-n}\int_{B_{kR}} |Du|^2\,dy \right] = \theta_0,
 \end{equation}
 for every $R> 0$.  Let us write $u_k (y) := u(ky)$, which is a minimizing constraint map in $W^{1,2}(B_8;\overline M)$. By \eqref{eq:rigid-re} and the bound $|u|\in L^\infty(\R^n)$, the sequence $\{u_k\}_{k=1}^\infty$ is bounded in $W^{1,2}\cap L^\infty(B_8;\overline M)$. Thus, by Lemma \ref{lem:strong}, there exists a minimizing constraint map $u_\infty \in W^{1,2}(B_4;\overline M)$ and a subsequence $k_i\to\infty$ such that $u_{k_i}\to u_\infty$ strongly in $W^{1,2}(B_4;\R^m)$.  Moreover, by Lemma \ref{lem:monot}, $u_\infty$ is 0-homogeneous in $B_4$ since, by \eqref{eq:rigid-re}, we have 
 $$
R^{2-n}\int_{B_R} |Du_\infty|^2\,dy = \lim_{i\to\infty} \bigg[ R^{2-n}\int_{B_R} |Du_{k_i}|^2\,dy \bigg] = \lim_{i\to\infty} \bigg[ (k_i R)^{2-n}\int_{B_{k_i R}} |Du|^2\,dx\bigg] = \theta_0,
 $$
 for a.e.\ $R\in(0,4)$. Since $0\in\Sigma(u_k)$, it follows from Corollaries \ref{cor:semi} and \ref{cor:e-reg} that $0\in \Sigma(u_\infty)$. 

 As $u_\infty$ is a minimizing constraint map, it follows from Corollary \ref{cor:dist} that $\rho\circ u_\infty$ is weakly subharmonic and continuous in $B_4$. Then the $0$-homogeneity of $\rho\circ u_\infty$ implies that $\rho\circ u_\infty$ is constant in $B_4$, cf.\ \cite[Lemma 2.6]{FGKS}. Since $0\in\Sigma(u_\infty)$, in fact we must have
 \begin{equation}\label{eq:ruinf}
 \rho\circ u_\infty = 0\quad\text{in }B_4.
 \end{equation}
Indeed, if $\rho\circ u_\infty =a$ for some $a>0$,  then we would have $|u_\infty^{-1}(\partial M)| = 0$,  and so by \eqref{eq:main-sys}, $\Delta u_\infty = 0$ in $B_4$, a contradiction to $0\in\Sigma(u_\infty)$. We remark that this argument can also be found in the proof of \cite[Theorem 6.1]{FGKS}.

 Owing to \eqref{eq:ruinf} and the strong convergence of $u_{k_i}\to u_\infty$ in $L^2(B_4;\R^m)$, we have $\rho\circ u_{k_i} \to 0$ in $L^2(B_4)$. By the subharmonicity of $\rho\circ u_{k_i}$ from Corollary \ref{cor:dist}, we have 
 \begin{equation}\label{eq:rigid1} 
 \rho\circ u_{k_i}\to 0\quad\text{in }L^\infty(B_4).
 \end{equation}
 Now we choose $\eta$ to be a small constant, depending only on $n$, $m$, $\partial M$, and $\Lambda$, as in Lemma \ref{lem:reg}. Then by \eqref{eq:rigid1} we can find a large $i_0$ such that $\rho\circ u_{k_i} \leq \eta$ in $B_4$ for all $i\geq i_0$. Recall  that $0\in \Sigma(u_{k_i})$, which implies $E(u,0,1)\geq\e_0^2$ by Lemma~\ref{lem:monot} and Corollary~\ref{cor:e-reg}. Hence  Lemma~\ref{lem:reg} implies that $B_1 \subset u_{k_i}^{-1}(\partial M)$. Rescaling back, we have proved that 
 $$
 B_{k_i}\subset u^{-1}(\partial M),
 $$
 for every $i\geq i_0$. Since $k_i\to\infty$, we deduce that $u^{-1}(\partial M) = \R^n$, which proves $u\in W_\loc^{1,2}(\R^n;\partial M)$. Finally, $u$ is a minimizing harmonic map, since  $W^{1,2}(D;\partial M)\subset W^{1,2}(D;\overline M)$ for any bounded domain $D\subset\R^n$. 
 \end{proof}

We are now ready to prove the main result of this section.

\begin{proof}[Proof of Theorem \ref{thm:uni-reg}]
As $M^c$ is convex, we can take $Y = M^c$ in Theorem \ref{thm:dist-reg}, which along with $\|u \|_{W^{1,2}(B_2)}\leq\Lambda$ yields a modulus of continuity $\omega$, determined {\it a priori} by $n$, $m$, $\partial M$, and $\Lambda$, such that for every $r\in(0,\frac{1}{2})$ and $x_0\in B_1$, 
\begin{equation}\label{eq:oscBrx0}
\osc_{B_r(x_0)} (\rho\circ u)\leq \omega(r).
\end{equation}
By Corollary~\ref{cor:dist}, the free boundary $B_1\cap\partial u^{-1}(M)$ is well-defined, so let us take an arbitrary point $x_0$ there. Now let $\eta$ be as in Lemma \ref{lem:reg} (with $c\Lambda^2$ in place of $\Lambda$ there), and take $\delta= \delta(n,m,\partial M,\Lambda) \in (0,\frac{1}{8})$ small such that $\omega(8\delta) \leq \eta$. By \eqref{eq:oscBrx0} and the fact that $\rho\circ u(x_0) = 0$, we have 
\begin{equation}\label{eq:ruBdx0}
\rho\circ u \leq \eta \quad\text{in }B_{8\delta}(x_0). 
\end{equation}
On the other hand, since $B_{8\delta}(x_0)\subset B_1(x_0)\subset B_2$, Lemma \ref{lem:monot} yields
\begin{equation}\label{eq:Eux0d}
E(u,x_0,8\delta)\leq E(u,x_0,1) \leq \int_{B_2} |Du|^2\,dx \leq \Lambda^2. 
\end{equation}
Now, as $x_0\in B_1\cap\partial u^{-1}(M)$, $B_{2\delta}(x_0)\setminus u^{-1}(\partial M) \neq\emptyset$. By \eqref{eq:ruBdx0} and \eqref{eq:Eux0d}, the contraposition of Lemma \ref{lem:reg} (applied to $u_{x_0,2\delta}(y) := u(2\delta y + x_0)$) implies that 
\begin{equation}
    E(u,x_0,2\delta) < \e_0^2. 
\end{equation}
Then, by Theorem~\ref{thm:e-reg}, $u\in C^{1,1}(B_\delta(x_0))$, and since $\delta\equiv \delta(n,m,\Lambda,\partial M)$, we have the estimate
$$
\|D^j u\|_{L^\infty(B_{\delta}(x_0))} \leq c,\quad j\in\{1,2\},
$$
where $c \equiv c(n,m,\Lambda,\partial M)$. 
This proves the validity of a uniform $C^{1,1}$ bound in a neighborhood of the free boundary.
On the other hand, since $\Delta u = 0$ in $B_1\cap u^{-1}(M)$, for any $y_0\in B_1\cap u^{-1}(M)$ with $\dist(y_0,B_1\cap\partial u^{-1}(M)) > \delta$, we may apply the interior estimates for harmonic functions, and use the universality of $\delta$ as well as $\int_{B_2} |Du|^2\,dx \leq\Lambda^2$, to obtain that 
$$
|D^j u(y_0)| \leq c, \quad j\in\{1,2\}.
$$
Combining the above two inequalities, we arrive at the desired conclusion. 
\end{proof}

\section{Discontinuities on free boundaries produced by degenerate obstacles}\label{sec:flat}

This section stems from the search for the minimal geometric conditions on the obstacle under which there are no discontinuities of $u$ near the free boundary. In Theorem \ref{thm:uni-reg} we proved the absence of discontinuities under the assumption that the obstacle is uniformly convex.
The main result in this section is Theorem \ref{thm:flat}, which shows the sharpness of the previous result: for convex obstacles, there may be discontinuities on the free boundary if we allow the obstacle to have flat sides. More precisely, we shall hereby call a subset $T\subset \partial M$ a \textit{flat piece} if $T$ is relatively open, connected, and all the principal curvatures of $T$ vanish.  Equivalently, $T$ is a relatively open connected subset of a hyperplane.

Before beginning the main proofs of this section, we collect some useful remarks.

\begin{rem}\label{rem:bdryreg}
Consider the setting of Theorem \ref{thm:flat}. By the boundary regularity theory \cite{DF}, whenever both $\partial\Omega$ and $u|_{\partial\Omega}$ are of class $C^2$, we have that $u \in C^{1,\alpha}(\{x\in \overline \Omega: \dist(x,\partial \Omega)< \delta\})$ for  every $\alpha\in(0,1)$, where $\delta >0$ depends only on $n$, $\alpha$, the $C^2$-character of $\partial\Omega$, and $u|_{\partial\Omega}$; in particular, 
$$
\dist(\Sigma(u),\partial\Omega)\geq\delta.
$$
Now, if $M^c$ is convex then $\Omega\cap u^{-1}(M)$  is an open set according to Corollary \ref{cor:dist}, and hence $\Delta u =0$ in $\Omega\cap u^{-1}(M)$. In particular, we have $u\in C^\infty(\Omega\cap u^{-1}(M))$. In conclusion, we obtain 
$$
u\in C_\loc^{1,\alpha}(\overline\Omega\cap u^{-1}(M))\cap C^\infty(\Omega\cap u^{-1}(M)). 
$$
\end{rem}

The previous remark shows that there are no discontinuities in a neighborhood of the fixed boundary $\partial\Omega$. As we now observe,  discontinuities are unavoidable in the interior.

\begin{rem}
As mentioned in the introduction, constraint maps, just as harmonic maps, generally develop discontinuities for \textit{topological reasons}. For the convenience of the reader, we briefly recall why this is the case. Given $\Omega\subset \R^n$, we  have the following topological fact:
\begin{equation}
    \label{eq:inclusion}
    \begin{rcases}
    u\in C^0(\overline \Omega)\\
    u=\id\text{ on } \partial \Omega
\end{rcases}
\quad \implies \quad 
\overline \Omega \subseteq u(\overline \Omega).
\end{equation}
In particular, if $M^c\Subset \Omega$ and $u\in W^{1,2}(\Omega;\overline M)$ is a minimizing constraint map with identity boundary values, then the conclusion in \eqref{eq:inclusion} does not hold, and so we must have $\Sigma(u)\neq \emptyset.$ This is essentially the setting of Theorem \ref{thm:flat}.

Implication \eqref{eq:inclusion} is most easily proved using the notion of \textit{topological degree}, which is  a way of counting (with multiplicity) the number of solutions $x\in \Omega$ to an equation
\begin{equation}
    \label{eq:eq}
    u(x)=y,
\end{equation}
The topological degree shows that in the setting of Theorem \ref{thm:flat}, we must have $\Sigma(u)\neq \emptyset$, but it can also be used to yield information about the image of $u$, since  $y\in u(\Omega)$ whenever $\deg(u,\Omega,y)\neq 0$. 
\end{rem}

\begin{rem}\label{rem:deg g}
In the setting of Theorem \ref{thm:flat}, we have

\begin{equation}
\label{eq:assumptions}
 \qquad u|_{\Omega\cap u^{-1}(M)} \text{ is surjective onto } \Omega \cap M.
\end{equation}
The proof of this fact is somewhat involved and relies on degree theory. In fact, in Proposition \ref{prop:deg} below we will prove a much stronger degree theoretic result, which we believe to be independently significant. {We also note that, from well-known properties of the topological degree, see in particular Theorem \ref{thm:degprops}\eqref{it:im}, our assumption $\deg(g) \neq 0$ on the boundary data implies that $g$ is surjective, namely $g(\partial\Omega)=\partial\Omega$.}
\end{rem}

Assuming for now that \eqref{eq:assumptions} holds,  let us show how it implies Theorem \ref{thm:flat}.  First, we state a simple consequence of the maximum principle:

\begin{lem}\label{lem:conn}
Let $u\in W^{1,2}(\Omega;\overline M)$ be a minimizing constraint map, and let $M$ and $\Omega$ be as in Theorem \ref{thm:flat}.  If $u=g$ on $\partial \Omega$ and $g(\partial \Omega)\subseteq\partial \Omega$,  then the set $\Omega\cap u^{-1}(M)$ is open and connected.
\end{lem}

\begin{proof}
Note that our assumptions  ensure that $\dist(\partial\Omega,\partial M) > 0$.  Thus, since $u = g$ on $\partial\Omega$ and $g(\partial\Omega) \subseteq \partial\Omega$,  we see that the open set $u^{-1}(M)$ 
contains a neighborhood of $\partial \Omega$, according to the continuity provided by Remark \ref{rem:bdryreg}. 

If the set $u^{-1}(M)$ was not connected, then there would be a connected component $C\subset u^{-1}(M)$ that is disjoint from $\partial \Omega$, thus $\partial C\subset \partial u^{-1}(M)\cap \Omega$. We would then have $\dist(u,\partial M)=0$ on $\partial C$ and $\dist(u,\partial M)>0$ on $C$. However, since the function $\dist(u,\partial M)$ is subharmonic in $u^{-1}(M)$ (since $u$ is harmonic in $u^{-1}(M)$ and $\dist(\cdot,\partial M)$ is convex), this would contradict the weak maximum principle.

\end{proof}

\begin{proof}[Proof of Theorem \ref{thm:flat}, assuming \eqref{eq:assumptions}]
Let $T$ be a flat piece of $\partial M$, that is, 
$$T\subset \{x\in \R^n:e\cdot x=c\}$$  for some $e\in \R^n$ and $c\in \R$,
and assume without loss of generality that $\{e\cdot x> c\}\subset M$. Let $M'\subset M$ be the open half-cylinder with base $T$, namely
$$M':=\{x=t e+\bar x\in \R^n: \bar x \in T \text{ and } t> 0 \}.$$
Fix any point $y \in T$.  According to the surjectivity asserted in \eqref{eq:assumptions} (see Proposition \ref{prop:deg} below),  there is a sequence $\{x_k\}_{k=1}^\infty\subset u^{-1}(M)\cap \Omega$ such that $u(x_k)\to y$ as $k\to\infty$. We may also assume that $x_k\to x_\infty$ for some $x_\infty \in \overline{u^{-1}(M)}\cap \Omega$. We claim that
$$x_\infty\in \Sigma(u);$$
once this is shown, we must have $x_\infty \in  \partial u^{-1}(M)$, since $u$ is harmonic in the open set $u^{-1}(M)$ according to \eqref{eq:main-sys} and Theorem \ref{thm:dist-reg}. 

Suppose, for the sake of contradiction, that $x_\infty\not \in \Sigma(u)$,  so there is $r>0$ with $B_r(x_\infty)\subset \Omega\setminus\Sigma(u)$. Since $T$ is a relatively open set in $\partial M$ and $u\in C^0(B_r(x_\infty))$, we have $u(B_r(x_\infty))\subset M'\cup T$, by taking $r$ smaller if necessary. In particular, \eqref{eq:main-sys} shows that $u$ is harmonic in $B_r(x_\infty)$. Hence, 
$e\cdot u \geq c \text{ and } \Delta(e\cdot u) = 0 \text{ in } B_r(x_\infty).$
Also, we have
$$e\cdot  u(x_\infty) = \lim_{k\to \infty}  e\cdot u(x_k) = e\cdot u(y) = c,$$ so the strong minimum principle asserts that $e \cdot u = c$ in $B_r(x_\infty)$.  
Since $e\cdot u$ is a harmonic function in $\Omega\cap u^{-1}(M)$ and $\Omega\cap u^{-1}(M)$ is connected (see Lemma \ref{lem:conn}),  the unique continuation principle shows that $e\cdot u=c$ in $\Omega\cap u^{-1}(M)$.  On the other hand, $u$ is continuous near $\partial \Omega$ and $u(\partial \Omega)=g(\partial \Omega)=\partial \Omega$ (since $\deg(g)\neq 0$, see Remark~\ref{rem:deg g}), so a contradiction is reached.  
 \end{proof}

\begin{figure}
    \centering
\includegraphics[scale=0.3]{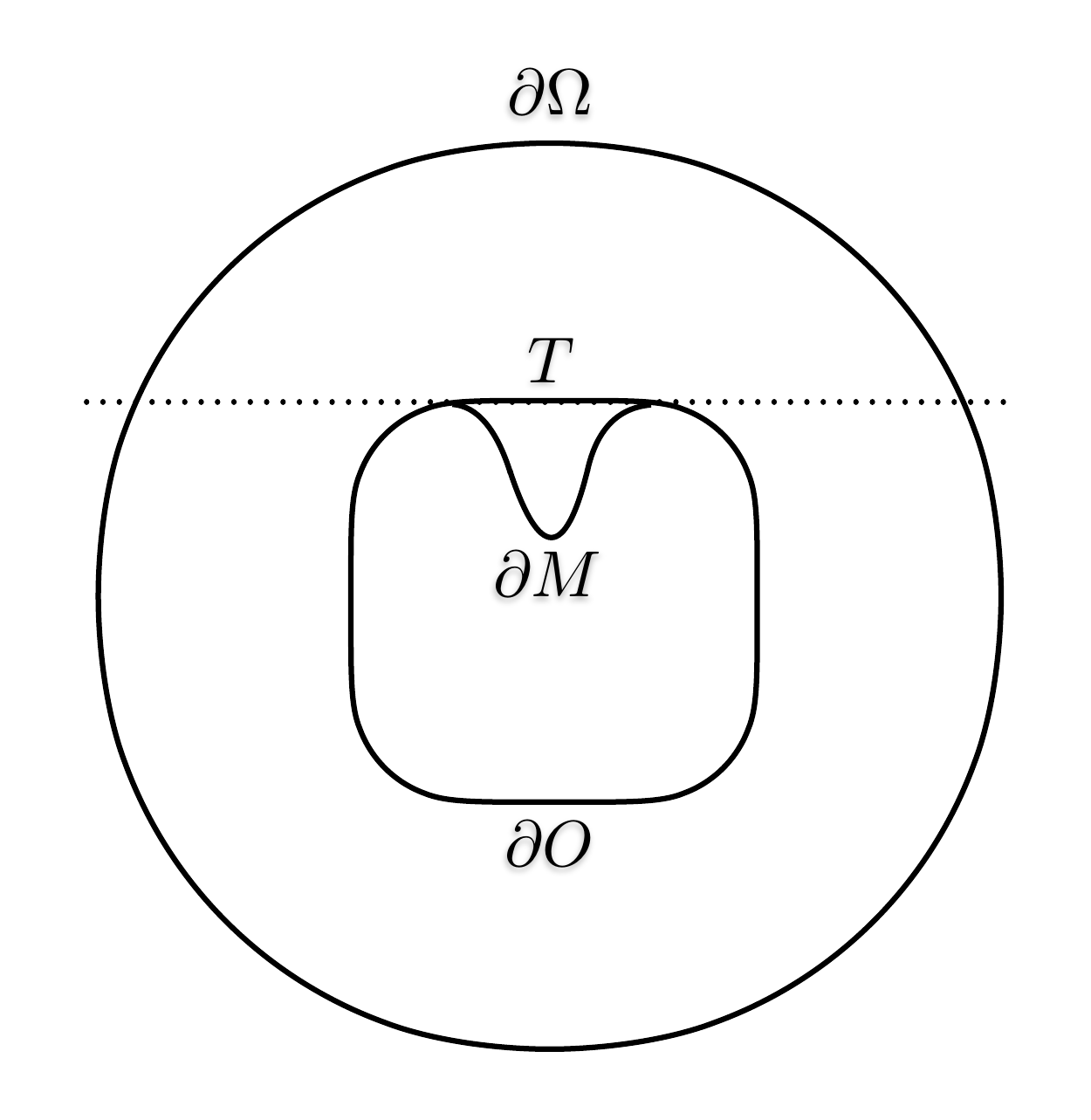}
    \caption{An example where Proposition \ref{prop:deg} and Corollary \ref{cor:flat} apply.}
    \label{fig:non-convex}
\end{figure}

The rest of this section is dedicated to proving  \eqref{eq:assumptions}. In fact, we will show a much more general result concerning the topological degree of minimizing constraint mappings, valid for non-convex obstacles, such as the one in Figure \ref{fig:non-convex}. We believe this more general result is independently interesting, and it is relevant in Corollary \ref{cor:flat} below.

\begin{prop}\label{prop:deg}
Let $\Omega\subset \R^n$ be a smooth uniformly convex domain, and  $M\subset \R^n$ be such that $M^c$ is connected and $M^c\Subset \Omega$. Let $O$ be the complement of the convex hull of $M^c$. Consider $g\in C^2(\partial\Omega;\partial\Omega)$ with $\deg(g) \neq 0$ and let $u \in W^{1,2}_g(\Omega;\overline M)$ be a minimizing constraint map. 
For any domain $O'\Subset O$,  we have
\begin{equation}
    \label{eq:Omega'}
    \deg(u,\Omega',\cdot) = \deg(g) \quad \text{on }  \Omega\cap O', \qquad \text{where } \Omega' := \Omega\cap u^{-1}( O').
\end{equation}
In particular,  $u|_{\Omega\cap u^{-1}(O)}$ is a surjective componentwise-harmonic map onto $\Omega\cap O$.
\end{prop}

Before beginning the proof of Proposition \ref{prop:deg}, we prove a simple topological lemma. This lemma does not rely on the specific geometric structure assumed in Proposition \ref{prop:deg}, and we use only that $O\subseteq M$.

\begin{lem}\label{lem:bdryvalue}
Let $u\colon \overline \Omega\to \overline M$ be a map such that $u|_{u^{-1}(O)}$  is continuous and $u(\partial \Omega)\subseteq \partial \Omega$. For every $ O'\Subset O$ we have $u(\partial\Omega')\subseteq \partial\Omega\cup \partial  O'$, where $\Omega'=\Omega\cap u^{-1}(O')$ is as in \eqref{eq:Omega'}.
\end{lem}

\begin{proof}
Note that $\partial\Omega'= \partial\Omega\cup [\Omega\cap \partial u^{-1}(O')]$. Since we assume $u(\partial\Omega)\subseteq \partial\Omega$, it suffices to prove that $$u(\Omega\cap \partial u^{-1}(O'))\subseteq \partial  O'.$$ For this purpose,  given $x_0 \in \Omega\cap \partial u^{-1}(O')$, there is a sequence $x_k\in \Omega\cap u^{-1}(O')$ such that $x_k\to x_0$, and since $u(x_k)\in O'$ we have $u(x_0)\in \overline{O'}$.

Assume by contradiction that $u(x_0) \in  O'$. Then $x_0 \in \Omega\cap u^{-1}(O')$, and since $\Omega\cap u^{-1}(O')$ is an open set, there is a ball $B_r(x_0)\subset \Omega\cap u^{-1}(O')$, impossible.

This proves that $u(x_0)\in \partial O'$, as desired. 
\end{proof}

\begin{proof}[Proof of Proposition \ref{prop:deg}]
Throughout the proof, we write $f$ for the signed distance function to $\partial \Omega$, with $f<0$ in $\Omega$.  Since $\Omega$ is smooth and uniformly convex,  there is a tubular neighborhood $\cN(\partial\Omega)$ of $\partial \Omega$ where $f$ is smooth and satisfies 
\begin{equation}
\label{eq:unicvx}
\Hess f(\xi,\xi)\geq c|\xi^\top|^2 \quad \text{on } \cN(\partial \Omega),
\end{equation} for some $c>0$ and all $\xi\in \R^n$. Here we write $\xi^\top:= \xi - (\nu\cdot \xi)\nu$ as in \eqref{eq:xitop}, where now $\nu$ (with a slight abuse of notation deviating from the introduction) is the outwards unit normal to $\partial \Omega$. For $\delta>0$ small enough, we write
$$O_\delta := \{-\delta < f\leq 0\}\subseteq \overline \Omega\cap \cN(\partial \Omega).$$
We split the proof into four steps.

\begin{step}
There is no ball $B\subset u^{-1}(O)$ such that $f\circ u=0$ in $B$.
\end{step}

As in the proof of Theorem \ref{thm:flat}, since we assume that $\deg(g)\neq 0$ we have that $g$ is a surjective map, i.e.\ $g(\partial \Omega)=\partial \Omega$.
Now assume towards a contradiction that $f\circ u=0$ in some ball $B\subset \Omega\cap u^{-1}(O)$. Since  $\Delta u=0$ in $\Omega\cap u^{-1}(O)$, we have
$$
0 = \Delta(f\circ u) =   \Hess f_u (Du,Du) \quad \text{ in } B.
$$
Since $u(B)\subseteq \partial \Omega$, $Du$ lies in the tangent bundle $T(\partial \Omega)$, i.e.\ $Du = (Du)^\top$. Thus, by \eqref{eq:unicvx}  we deduce that $Du=0$ in $B$. According to Lemma \ref{lem:conn}, the open set $\Omega\cap u^{-1}(O)$ is connected, and as $u$ is harmonic there, 
by unique continuation
we must have that $u$ is constant in $\Omega\cap u^{-1}(O)$. This contradicts the fact that $u(\partial \Omega)= g(\partial \Omega)= \partial \Omega$.

\begin{figure}
    \centering
\includegraphics[scale=0.3]{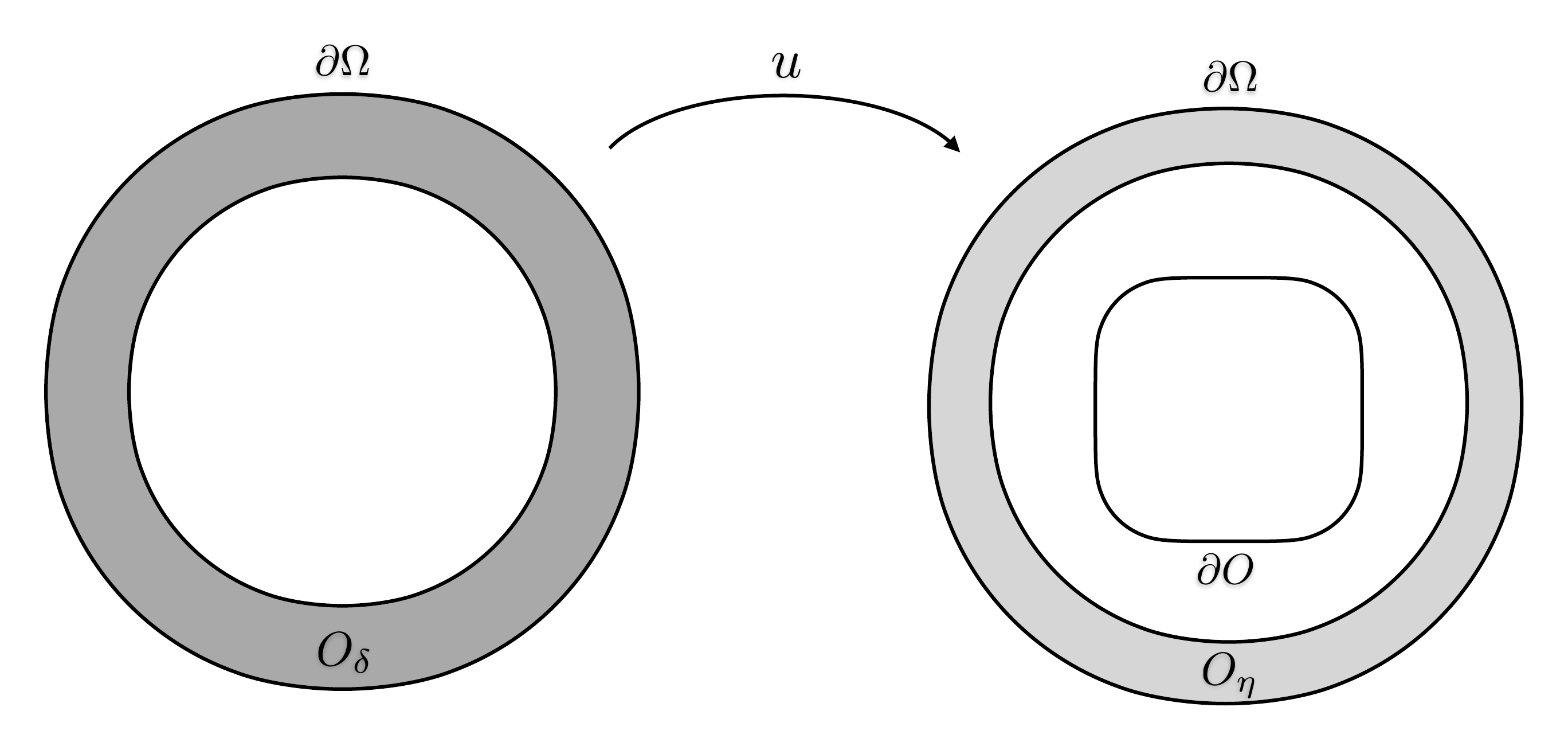}
    \caption{An illustration of Step 2 in the proof of Proposition \ref{prop:deg}. All pre-images of the thin boundary layer $O_\eta$ in the target are contained in another thin boundary layer $O_\delta$ in the domain.}
    \label{fig:non-convex2}
\end{figure}

\begin{step}
    For every sufficiently small $\delta>0$, there exists $\eta > 0$ such that $\Omega\cap u^{-1}(O_{\eta}) \subset \Omega\cap O_{\delta}$. 
\end{step}

Suppose, for the sake of contradiction, that there are some small $\delta>0$ and a sequence $\eta_k\to 0$ for which points $y_k\in \partial O_{\eta_k}$ and $x_k\in u^{-1}(y_k)\cap \Omega\setminus O_{\delta}$ can be found. Without loss of generality, we can assume that $x_k\to x_\infty$ and $u(x_k)=y_k\to y_\infty$.  Clearly, $x_\infty\in \Omega\setminus O_{\delta}$ and $y_\infty\in \partial\Omega = f^{-1}(0)$. Since Theorem \ref{thm:dist-reg} applies with $Y=O^c$, we have 
\begin{align*}
\dist(u(x_\infty),\partial O) &=\dist(u(x_\infty),O^c) =\lim_{k\to\infty}\dist(u(x_k),O^c)\\
&=\lim_{k\to\infty}\dist(u(x_k),\partial O) 
= \dist(y_\infty,\partial O)\geq \dist(\partial \Omega,\partial O)>0,
\end{align*}
and so $x_\infty\in \Omega\cap u^{-1}(O)$. In particular,  $u$ is continuous at $x_\infty$ and $u(x_\infty)=y_\infty$. Since $f$ is convex in $\cN(\partial \Omega)$,  by continuity and harmonicity of $u$ in $\Omega\cap u^{-1}(O)$ we conclude that the function $f\circ u$ is subharmonic in a ball $B$ centered at $x_\infty$.  However,  Corollary \ref{cor:image} shows that
\begin{equation*}
\label{eq:range-fu}
f\circ u \leq 0 \quad \text{in } \Omega
\end{equation*}
therefore $x_\infty$ is an interior maximum for $f\circ u$ in $B$, since $f\circ u(x_\infty)=f(y_\infty)=0$.  Therefore, by the strong maximum principle we deduce that $f\circ u = 0$ in $B$.  This contradicts Step 1,  completing the proof of Step 2.

\begin{step}
There is $d\in \Z$ such that    $\deg(u,\Omega',\cdot) = d$ on $\Omega\cap O'$ for every $O'\Subset O$.
\end{step}

Since $u\in C^{0}(\overline\Omega\cap u^{-1}(O))$ by Remark \ref{rem:bdryreg},  we have $u \in C^{0}(\overline{\Omega'})$ (recall $\Omega' = \Omega\cap u^{-1}(O')$).  By Lemma \ref{lem:bdryvalue} we have $\Omega\cap O' \subset \R^n\setminus u(\partial\Omega')$ and so,  by Lemma \ref{lem:conn}, $\Omega\cap O'$ is contained in a connected component of $\R^n\setminus u(\partial\Omega')$.  By Theorem \ref{thm:degprops}\eqref{it:const}, $\deg(u,\Omega',\cdot)$ is constant in each connected component of $\R^n\setminus u(\partial\Omega')$, and so there is $d\in \Z$ such that
$$
\deg(u,\Omega',\cdot) = d  \quad\text{in }\Omega\cap O',
$$
as desired. 

\begin{step}
We have $d=\deg(g)$. 
\end{step}

Note that $f\circ u$  is subharmonic in a neighborhood of $\partial \Omega$ inside $\Omega$, and it achieves a maximum on $\partial \Omega$.  Thus, the Hopf lemma combined with Step 1 shows that 
\begin{equation*}
    \label{eq:hopf}
D_\nu (f\circ u) \geq 2 c \quad \text{on }\partial \Omega,
\end{equation*} for some constant $c>0$, where we recall that $\nu$ is the outwards unit normal to $\partial \Omega$, extended to $\cN(\partial \Omega)$ via the formula $\nu=\nabla f$.  Since $f$ and $u$ are $C^1$ near $\partial \Omega$,  we deduce that 
\begin{equation}
\label{eq:lowerboundDnuu}
D_\nu(f\circ u) \geq c \quad \text{on } \Omega\cap O_{\delta}, 
\end{equation}
provided that $\delta>0$ is sufficiently small.  

Consider an oriented orthonormal basis $\mathcal B_x:=\{\tau^1_x,\dots,\tau^{n-1}_x,\nu_x\}$ at a point $x\in \partial \Omega$, and write $\mathcal{T}_x :=I-\nu_x\otimes \nu_x$ for the matrix with columns $\tau_x^1,\dots, \tau_x^{n-1}$. 
Recall that $\deg(g)$ is an integer which can be computed as follows: if $y\in \partial \Omega$ is a regular value\footnote{Recall that $y$ is a regular value of $g$ if the differential $dg_x$ is non-singular at all points $x\in g^{-1}(y)$. In particular, by the inverse function theorem, a regular value of $g$ has only a finite number of pre-images.} of $g$ and $g^{-1}(y)=\{x_1,\dots, x_{i_0}\}$ for some points $\{x_i\}\subset \partial \Omega$, then
$$
\deg(g) = \sum_{x_i\in g^{-1}(y)} \sgn(\det dg_{x_i}),
$$
where we see $dg_{x_i}\colon T_{x_i}\partial \Omega\to T_{y}\partial \Omega$ as an invertible linear map between $(n-1)$-dimensional vector spaces. More precisely, 
since $u=g$ on $\partial \Omega$, for each $x \in \partial\Omega$ we have
$$dg_{x} = \mathcal{T}_{u(x)}^t Du(x)\mathcal{T}_{x},$$
where the superscript ``$t$'' denotes the transpose.

Moreover, there are  relatively open subset $U_i, V\subset \overline \Omega$ such that $x_i\in U_i$ and $y\in V$, 
\begin{equation}
\label{eq:glocdiffeo}
g^{-1}(V\cap \partial \Omega)=\cup_{i=1}^{i_0} U_i\cap \partial \Omega , \quad\text{where } g\colon U_i \cap \partial \Omega \to V\cap \partial \Omega  \text{ is a diffeomorphism,}
\end{equation}
for all $i$. Note that, on $\partial \Omega,$ we also have
$$\nu_{u(x)}^t D u(x) \mathcal{T}_x = D_\tau(f\circ u)(x)=0, \qquad \nu_{u(x)}^t Du(x) \nu_x = D_\nu(f\circ u)(x).$$
Thus, for $x\in \partial \Omega$, we can write $Du(x)$ with respect to the bases $\mathcal B_x$ and $\mathcal B_{u(x)}$ as follows:
$$Du(x) = 
\begin{bmatrix}
    dg_x & * \\ 0 & D_{\nu}(f\circ u)(x)
\end{bmatrix},
$$
where $*$ denotes an unspecified entry of the matrix.
Since $u$ is $C^{1,\alpha}$ near $\partial \Omega$, by shrinking $\delta,$ $U_i$, and $V$ if needed,  we see from \eqref{eq:lowerboundDnuu} that
$$\sgn(\det Du(x))=\sgn(\det dg_{x_i}) \quad \text{for all } x\in U_i\cap \Omega.$$  

Let $\eta$ be determined in terms of $\delta$ through Step 2. By shrinking $U_i$ and $V$ further if necessary, we can also assume that $V\subset \overline \Omega\cap O_{\eta}$, so $u^{-1}(V)\subset \Omega\cap O_{\delta}$, and  by \eqref{eq:lowerboundDnuu} and \eqref{eq:glocdiffeo} we have
$$
u^{-1}(V) = \cup_{i=1}^{i_0} U_i, \qquad u\colon U_i \to V \text{ is a diffeomorphism.}
$$
Thus any point $y'\in V\cap \Omega$ is a regular value for $u$ which is not in $u(\partial \Omega)=\partial \Omega$ and so
$$\deg(u,\Omega',y')=
\sum_{\tilde x_i\in u^{-1}(y')} \sgn(\det Du(\tilde x_i))=
 \sum_{x_i\in g^{-1}(y)} \sgn(\det dg_{x_i})= \deg(g)$$
for every $\e>0$ sufficiently small, as wished.

\medskip

For the surjectivity claim, note that Theorem \ref{thm:degprops}\eqref{it:im} and the previous two steps show that $\Omega\cap O'\subseteq u(\Omega')$ for all $O'\Subset O$. The conclusion follows by taking the union over $O'$.
\end{proof}

We conclude this section by noting that, for general non-convex obstacles, the non-coincidence set $u^{-1}(M)$ is not necessarily open:

\begin{cor}\label{cor:flat}
Let $\Omega,M,O$ and $u$ be as in Proposition \ref{prop:deg}. Assume that $M^c$ is non-convex and that $\partial O$ contains a flat piece $T$. Suppose in addition that
\begin{equation}
    \label{eq:range} u(\Omega)\subseteq \overline O.
\end{equation}
Then $u^{-1}(M)$ is not an open set.
\end{cor}

\begin{proof}
Note that, by \eqref{eq:range}, $u\in W^{1,2}_g(\Omega;\overline O)$ is a minimizing constraint map. Since $O^c$ is convex, according to Theorem \ref{thm:flat} there is a point $\bar x\in \partial u^{-1}(O)\cap \Sigma(u)$ and a sequence of points $x_k\to \bar x$ such that $u(x_k)\in O\subset M$  and $u(x_k)\to \bar y$ for some $\bar y\in \partial O\setminus \partial M$. In particular, since $u$ is continuous over $u^{-1}(O)$, we have 
$$
\liminf_{r\to 0}\|\dist(u,\partial M)\|_{L^\infty(B_r(\bar x))} >0,
$$
thus $\bar x\in u^{-1}(M)$. On the other hand,  since $\bar x\in \Sigma(u)$, there is no ball around $\bar x$ that is contained in $u^{-1}(M)$, as otherwise $u$ would be harmonic near $\bar x$, according to \eqref{eq:main-sys}. Thus, $u^{-1}(M)$ is not open.
\end{proof}

\begin{rem}
We note that assumption \eqref{eq:range} is satisfied for many non-convex obstacles, such as the one in Figure \ref{fig:non-convex}. Indeed, for such an obstacle the nearest-point projection map $\pi\colon M\to O$ is 1-Lipschitz and so, by comparing a minimizing constraint map $u\in W^{1,2}_g(\Omega;\overline M)$ with a competitor $\pi \circ u$, we can deduce that in fact $u\in W^{1,2}_g(\Omega;\overline O).$
\end{rem}


\section{Branch points and free boundaries}\label{sec:branch}

\subsection{General remarks on branch points}\label{subsec:remks}
Our goal now is to understand the interaction between branch points and the free boundary. We begin by collecting some general remarks. Similarly to \eqref{eq:branching}, we state the following definition.

\begin{defn}\label{def:branching}
Let $u\colon \Omega \to \R^m$ be a minimizing constraint map.  We call $x_0\in \Omega\setminus \Sigma(u)$ a {\normalfont branch point} of $u$, and we write $x_0\in \cB(u)$, if $Du(x_0)$ does not have full rank. This leads to a decomposition
$$\cB(u) = \bigcup_{k=0}^{\min\{m,n\}-1} \cB_k(u),\qquad  \cB_k(u) := \{x\in \Omega\setminus \Sigma(u): \mathrm{rank}(Du(x))=k \}.$$
Given $x_0\in \cB(u)$, we say that:
\begin{enumerate}
\item $x_0$ is a {\normalfont false branch point} if $u(B_r(x_0))$ is a $C^1$ manifold, for all $r>0$ sufficiently small;
\item $x_0$ is a {\normalfont true branch point} if it is not a false branch point.
\end{enumerate}
\end{defn}

\begin{figure}
\centering
\begin{subfigure}{.5\textwidth}
  \centering
  \includegraphics[width=.6\linewidth]{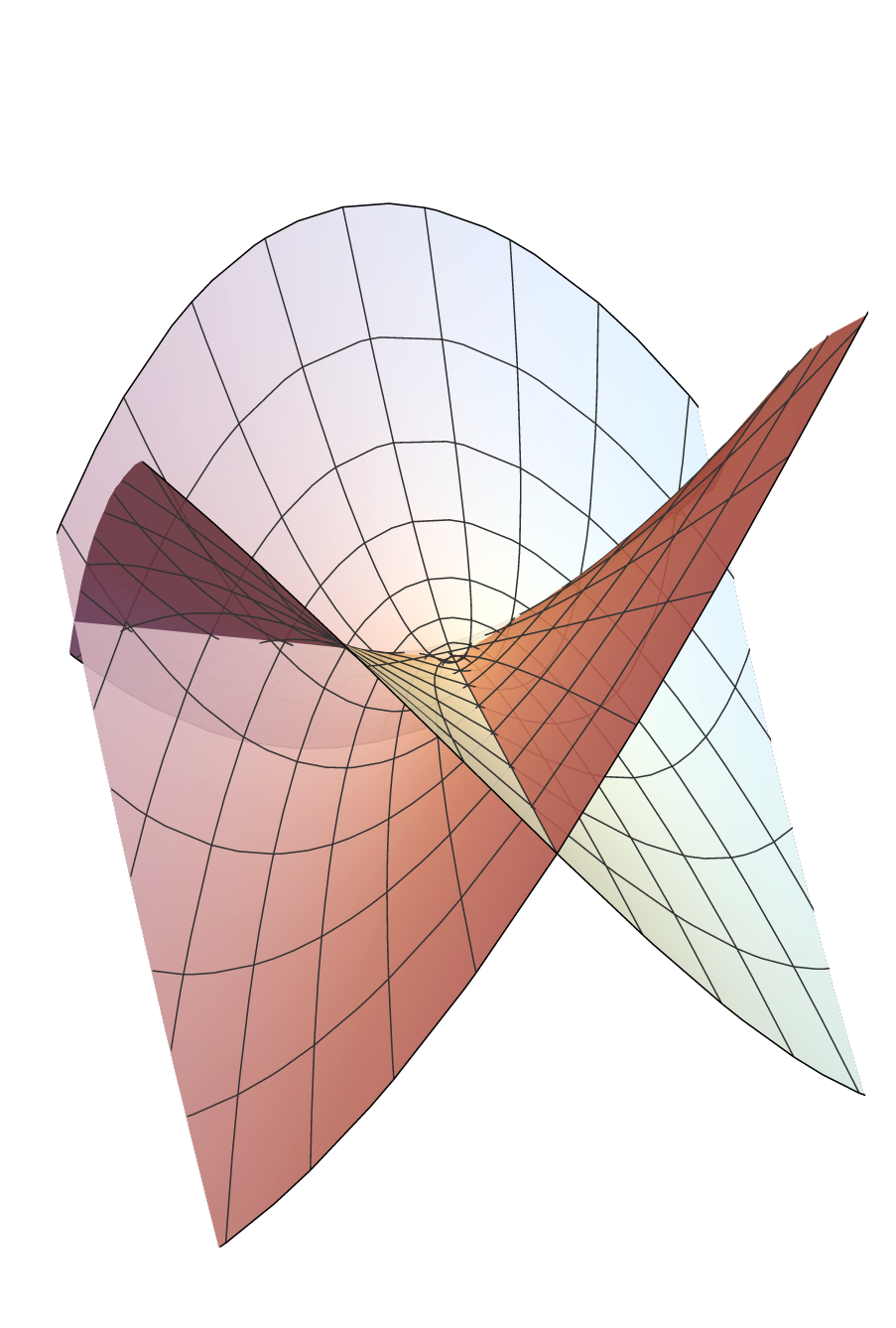}
  \caption{$k=3$: true branching}
  \label{fig:sub1}
\end{subfigure}%
\begin{subfigure}{.5\textwidth}
  \centering
  \includegraphics[width=.65\linewidth]{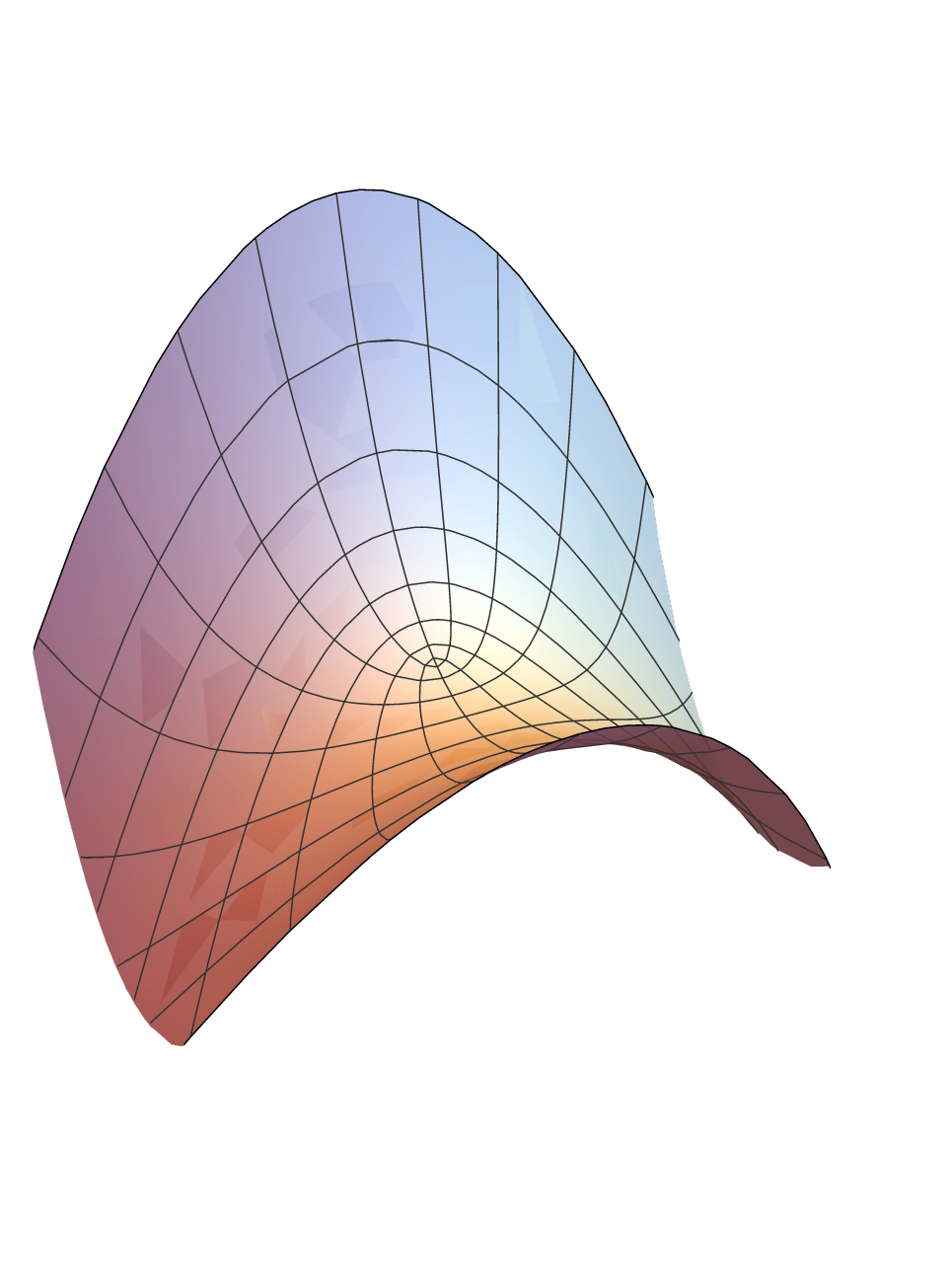}
  \caption{$k=4$: false branching}
  \label{fig:sub2}
\end{subfigure}
\caption{Plots of $u_k([-1,1]^2)$.}
\label{fig:branches}
\end{figure}

In order to illustrate the above concepts, consider the maps 
\begin{equation}
    \label{eq:defuk}
    u_k\colon \C\simeq \R^2\to \R^3, \quad z\mapsto  (z^2,\mathrm{Re}\,z^k),
\end{equation}
which are component-wise harmonic. In particular, if $M\subset \R^3$ is a smooth domain such that $u(\C)\subseteq \overline M$, then $u$ is a minimizing constraint map in $W^{1,2}(\Omega;\overline M)$ for any $\Omega\subset \C$. 
These maps also satisfy $\cB(u_k)=\cB_0(u_k)=\{0\}$,  and moreover $0$ is a true branch point if $k$ is odd, and it is a false branch point if $k$ is even, see Figure \ref{fig:branches} for an illustration.

Before proceeding further, we compare Definition \ref{def:branching} with other definitions of branch points occurring in the literature; we believe that this will also clarify the nature of our definition.

\begin{rem}[Branching in Geometric Function Theory] In the function theory literature \cite{Rickman1993}, one typically considers maps between domains of the same dimension, i.e., $m = n$. The branch set $\cB^{\textup{gft}}(u)$ is defined as the set of points near which a map is not a local homeomorphism.  Clearly, for a $C^1$ map $u$ we have
\begin{equation*}
    \cB^\textup{gft}(u)\subseteq \cB(u) ,
\end{equation*}
according to the Inverse Function Theorem. This inclusion is generally strict, and this highlights the fact that our definition of branching is \textit{analytical}, rather than topological. In fact, for open mappings of class $C^1$, we always have that
$$\cB^\textup{gft}(u) \subseteq \bigcup_{k=0}^{n-2} \cB_k(u);$$
in other words, $u$ is a local homeomorphism near any point in $\cB_{n-1}(u)$ \cite[Theorem 1.4]{Church1963}. 
For example, the homeomorphism $x\mapsto x^3$ of the real line has a rank zero point at the origin.
\end{rem}

\begin{rem}[Branching for parametric minimal surfaces]
    In the theory of parametric minimal surfaces, one considers weakly conformal parametrizations $u\colon \R^2\to \R^m$ of a (possibly branched) minimal surface $S=u(\R^2)$. Precisely, $u$ must satisfy:
    \begin{enumerate}
        \item \label{it:wconf} weak conformality: $|u_x|=|u_y|$ and $u_x\cdot u_y=0$;
        \item \label{it:min} minimality: $\Delta u = 0.$
    \end{enumerate}
    If $u$ satisfies \eqref{it:wconf} then $\cB(u)=\cB_0(u)$. Also, it is known that if, in addition, \eqref{it:min} holds and $u$ is minimizing, then $\cB(u) = \emptyset$, under appropriate conditions on the boundary data, see \cite{Gulliver1973} and the references therein. However, this result is not applicable in our setting, since minimizing constraint maps generally do not satisfy \eqref{it:wconf}, a simple example being given by \eqref{eq:defuk}.
\end{rem}

 We also note that, for harmonic maps between two-dimensional domains, a fine description of their behavior at the branch points was presented in \cite{Wood1977}.

We now want to obtain some information on the sets $\cB_k(u)$, in particular concerning their dimension and rectifiability properties. For reasons that will become apparent later, we will focus on the case $k=0.$ The next result follows from essentially standard considerations, cf.\ \cite{HL} for the case of harmonic maps. As we shall explain below, the following result is essentially contained in \cite{NV2}.

\begin{prop}\label{prop:dimB0}
Let $M\subset\R^m$ be a smooth domain and let $u\in W^{1,2}(\Omega;\overline M)$ be a non-constant minimizing constraint map.  For any ball $B$ such that $2B\Subset \Omega\setminus \Sigma(u)$,  we have the estimates
    \begin{align}
    \label{eq:measestS}
    \cH^{n-2}(B\cap \cB_0(u))& \leq c ,\\
    \label{eq:measestimage}
    \cH^{\frac{n-2}{2}} (u(B\cap \cB_0(u)))& \leq c,
\end{align} 
for some constant $c$ depending only on $n, \partial M$, $\|Du\|_{L^\infty(2B)}$, and $\int_{2B} |u-(u)_{2B}|^2 dx$.
\end{prop} 

\begin{proof}
Recall that, by definition, $\cB_0(u)$ is a closed set disjoint from $\Sigma(u)$.
By the Euler--Lagrange system \eqref{eq:main-sys}, any minimizing constraint map is a solution of an elliptic system of the type
\begin{equation}\label{eq:main-sys-re}
    \Delta u^i = a_{ij}^\alpha D_\alpha u^j\quad\text{in }\Omega, \qquad a_{ij}^\alpha \in L^2(\Omega)\cap L_\loc^\infty(\Omega\setminus\Sigma(u)),
\end{equation}
where the regularity of the coefficients follows from Theorem \ref{thm:e-reg}. In particular, as remarked in \cite{HL,NV2} (and as we formulate more precisely in Appendix \ref{app:freq}), due to the diagonal structure of this system, the standard scalar methods concerning the critical set of elliptic equations are applicable. Therefore the measure estimate follows as in the scalar case \cite{NV2}.  We remark here that the methods of \cite{HL} are not directly applicable, since minimizing constraint maps are not smooth near regular points; their methods require increased smoothness of the solution as the frequency increases. Nonetheless, the work \cite{NV2} deals precisely with elliptic equations with rough coefficients.

Concerning estimate \eqref{eq:measestimage}, we argue as follows. 
By \eqref{eq:measestS},  for any $\delta>0$ we can take points $x_i\in \cB_0(u)$ and radii $r_i \in (0,\delta)$, for $i\in I$, such that $$B\cap \cB_0(u)\subset\bigcup_{i\in I} B_{r_i}(x_i), \qquad 
\sum_{i\in I} r_i^{n-2} \leq c.
$$
Since $u\in C^{1,1}(B)$ and $|Du| = 0$ on $B\cap \cB_0(u)$, we can find a uniform $\lambda>0$ such that $u(B_{r_i}(x_i))\subset B_{\lambda r_i^2}(u(x_i))$, for every $i\in I$.  Since $\{B_{\lambda r_i^2}(u(x_i)) : i \in I\}$ is a covering of $u(B\cap \cB_0(u))$, letting $\delta\to 0$ we get 
$$
\cH^{\frac{n-2}{2}}(u(B\cap \cB_0(u)))\leq \lambda^{\frac{n-2}{2}}c.
$$
as desired. 
\end{proof}

\begin{rem}\label{rem:measbranch2}
Note that the estimate \eqref{eq:measestimage} on the image of $\cB_0(u)$ is stronger than what one could deduce from Sard's Theorem whenever $n \leq 3 \leq m$ or $n > m$, since minimizing constraint maps are only $C^{1,1}$.
\end{rem}

\begin{rem}[Branch points with higher rank]
For a smooth harmonic map $u$, Sampson \cite{S} proved the following unique continuation result: if $\cB_1(u)$ has non-empty interior, then $\cB_1(u)= \Omega$ and $u$ maps into a geodesic. Then, it was shown in \cite{HL} that if $u$ does not map into a geodesic, the set $\cB_1(u)$ has a locally finite $\cH^{n-1}$-measure.  
Perhaps surprisingly, the same unique continuation type properties do not hold in general\footnote{An obvious exception is for smooth harmonic maps into analytic manifolds, since in that case $\textup{rank}(Du)$ is constant in an open dense set.} for the sets $\cB_k(u)$ whenever $k\geq 2$, see e.g.\ \cite{Jin1991}. 
\end{rem}

\subsection{Branch points at free boundary points}\label{subsec:branch}

We now turn our attention to the interaction between $\cB(u)$ and the free boundary $\partial u^{-1}(M)$. We first make the observation that, under rather general circumstances, one can make degree-theoretic arguments to produce branch points on the free boundary. Recall that  $\mathbb B^n\subset \R^n$ denotes the open unit ball.

\begin{prop}\label{prop:monodromy}
Let $n\geq 3$ and let $u\in W^{1,2}_g(\mb B^{n};\R^n\backslash\mb B^n)$ be a minimizing constraint map. Assume that $g\in C^2(\mb S^{n-1};\mb S^{n-1})$ satisfies $\lvert\deg(g)\rvert \geq 2$, and that, for a.e.\ $\e>0$ sufficiently small, the set  $\{|u|\geq 1+\e\}$ is a topological annulus. Then 
$$\bigcup_{k=0}^{n-2} \cB_{k}(u) \cap \partial\{|u|>1\}\neq\emptyset.$$
\end{prop}

\begin{proof}
From Theorem \ref{thm:uni-reg} and Sard's theorem, we deduce that a.e.\ $\e>0$ sufficiently small is a regular value of $|u|$. Thus 
\begin{equation}
    \label{eq:nonvanishingD|u|}
    D|u|= (u\cdot Du)/|u|\neq 0 \quad \text{ on } N_\e := \partial\{|u|>1+\e\},
\end{equation}
 and $N_\e$ is then a smooth manifold.  Since $\lvert \deg (g)\rvert\geq 2$, the map $g\colon \mb S^2\to \mb S^2$ is surjective. Thus, our assumption that the set $\{|u|\geq 1+\e\}$ is a topological annulus for a.e.\ $\e>0$ sufficiently small, essentially amounts to saying that the manifold $N_\e$ is connected. 

Let us suppose, for the sake of contradiction, that 
\begin{equation}
    \label{eq:contraBPs}
    \bigcup_{k=0}^{n-2} \cB_{k}(u) \cap N_0 = \emptyset.
\end{equation}
We write, as usual, $\Pi \circ u = u/|u|$, so that $u=|u|(\Pi\circ u)$ and 
\begin{equation}
    \label{eq:Du-BPs}
    Du = (\Pi\circ u)\otimes D|u| + |u| D(\Pi \circ u).
\end{equation}
Since  $1\leq |u|$ is $C^{1,1}$ by Theorem \ref{thm:uni-reg}, and attains a minimum on $N_0$, we have 
\begin{equation}
    \label{eq:FBBP}
    D|u|=0\quad \text{ on }N_0.
\end{equation}
Notice also that, since $\Pi$ is 0-homogeneous, we have
\begin{equation}
    \label{eq:Pi-0hom}
    u\cdot D (\Pi \circ u )=0
\end{equation}
and so $D(\Pi \circ u)$ is always a singular matrix. Hence, 
combining \eqref{eq:contraBPs}--\eqref{eq:Pi-0hom},  we have 
\begin{equation*}
    \label{eq:rankPi}
    \textup{rank}(D(\Pi \circ u))=n-1 \quad \text{on } N_0,
\end{equation*}
and by continuity the same also holds on $N_\e$, for $\e>0$ small enough.
Let us now choose $\e$ so that \eqref{eq:nonvanishingD|u|} holds.  It follows from \eqref{eq:Du-BPs} and \eqref{eq:Pi-0hom} that $u\cdot Du= |u| D|u|$, therefore
\begin{equation}
    \label{eq:uislocdiff}
    \textup{rank}(Du)=n\quad \text{on } N_\e.
\end{equation}
We have seen in Proposition \ref{prop:deg} that 
$$\deg(u,\mb B^n \cap \{|u|>1\},y) = \deg(g) $$
whenever $1<|y|<\lambda$. Thus, since $N_\e=\{|u|=1+\e\}$ is connected, taking $|y|=1+\e$ we have $u^{-1}(y)\subset N_\e$ and so, according to \eqref{eq:uislocdiff} and the definition \eqref{eq:degu} of the local degree, we see that
$$\deg(u,\mb B^n\cap \{|u|>1\},y) = \sum_{x\in u^{-1}(y)} \textup{sgn}(\det Du(x)) = \deg(g).$$
 We now argue similarly to  the last step in the proof of Proposition \ref{prop:deg}. For $x\in N_\e$, let $\mathcal B_x:=\{\tau^1_x,\dots,\tau^{n-1}_x\}$ be an oriented orthonormal basis for $T_x N_\e$. Thus $D_{\tau^i} |u|=0$ on $N_\e$ and, by \eqref{eq:Du-BPs}, we have
$$D_{\tau} u(x) = Du(x)\tau_x=D(\Pi\circ u)(x)\tau_x.$$
Note that $\nu_x:=\frac{D|u|(x)}{\lvert D|u|(x)\rvert}$ is the unit normal to $N_\e$, which is well-defined by \eqref{eq:nonvanishingD|u|}, and so
$$(D_\nu|u|)(x) = (D|u|)(x) \cdot \nu_x = \lvert (D|u|)(x)\rvert>0.$$
Hence, if we consider the map $p_\e\colon N_\e\to \mb S^{n-1}$ defined by $p_\e := u/(1+\e)$, the last two identities show that, on $N_\e$, we have
$$\textup{sgn}(\det Du) = \textup{sgn}(\det dp_\e),$$
where as usual $dp_\e\colon TN_\e \to T\mb S^{n-1}$ is the differential between tangent spaces.  In particular, we deduce that
\begin{equation}
    \label{eq:degpeps}
    \deg(p_\e) = \sum_{x: p_\e(x)=y/(1+\e)} \textup{sgn}(\det dp_\e(x)) = \deg(g).
\end{equation}
Note that, by \eqref{eq:uislocdiff}, $p_\e\colon N_\e \to \mb S^{n-1}$ is a local diffeomorphism. Since $n\geq 3$, $\mb S^{n-1}$ is simply connected and so, since $N_\e$ is connected, $p_\e$ is actually a \textit{global} diffeomorphism; in particular, we must have $\deg p_\e=1$. However,  this contradicts \eqref{eq:degpeps}, since we assumed that $\lvert\deg(g)\rvert \geq 2.$
\end{proof}

Under the assumption that the non-coincidence set $\{|u|>1\}$, or a small perturbation of it, is a topological annulus, Proposition \ref{prop:monodromy} produces branch points of \textit{low rank} on the free boundary.  This topological condition seems hard to verify in general, but we expect it to hold in any sufficiently symmetric scenario. 

 In any case, we are especially interested in knowing whether there are points with \textit{zero rank} on the free boundary, and it seems that the existence of such points cannot follow from purely topological considerations, such as those used in the proof of Proposition \ref{prop:monodromy}.
Rank zero points are especially relevant to us since, at those points, \eqref{eq:main-sys} becomes a \textit{degenerate} obstacle problem, regardless of what obstacle $M^c$ is; furthermore, for uniformly convex obstacles,  these are the only points at which \eqref{eq:main-sys} is degenerate.  To construct such points, the proof of Proposition \ref{prop:monodromy} naturally leads us to consider axially symmetric maps: if $u\colon \mb B^3\to \R^3\backslash \mb B^3$ is a $C^1$ map which rotates at least twice about the $z$ axis, then a simple calculation reveals that the gradient of $u$ along the $(x,y)$ plane must \textit{vanish} on the $z$ axis, cf. \eqref{eq:Df-complex}.  In other words, the differential of the map $p_\e$ constructed in the proof of Proposition \ref{prop:monodromy} vanishes on the $z$-axis. The construction and complete analysis of points in $\cB_0(u)\cap \partial \{|u|>1\}$ in the axially symmetric case is actually quite involved, and it will occupy the final sections of this paper.

\begin{figure}
\centering
  \includegraphics[width=.65\linewidth]{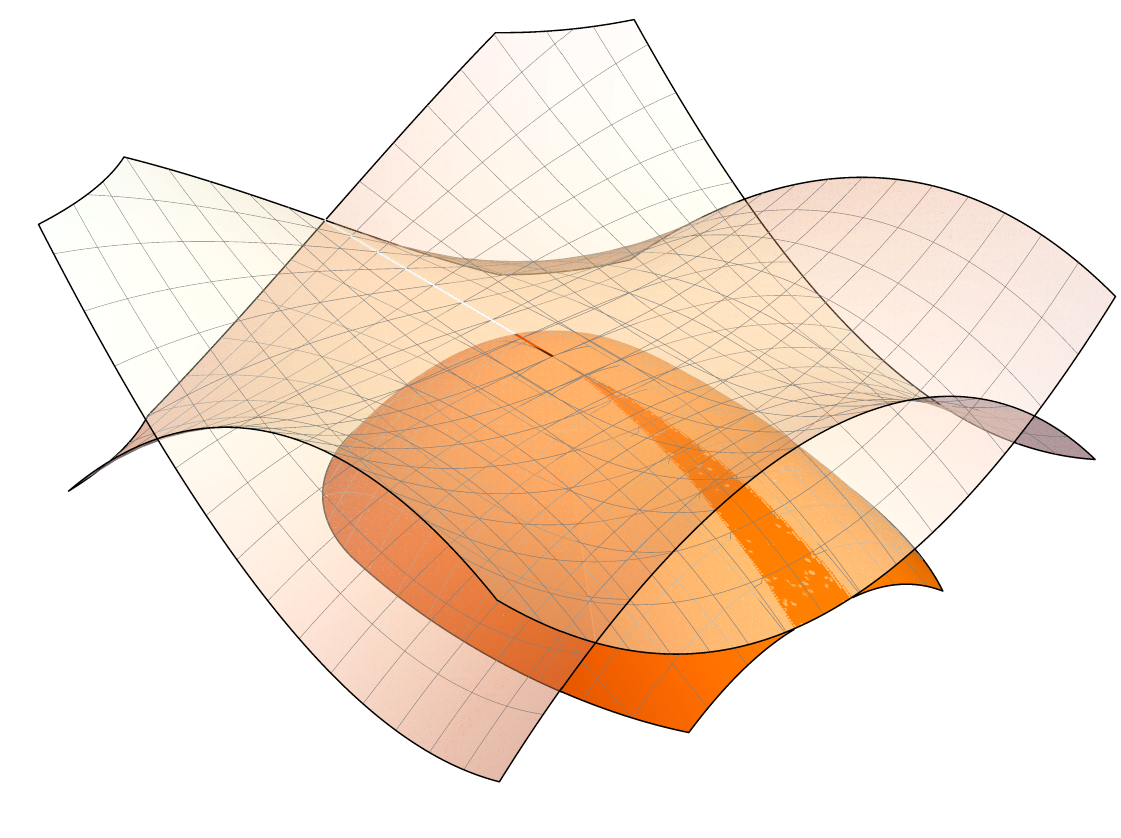}
  \caption{The image of $u_5$,  which touches the obstacle over a cone with tip at the origin.  Note that, over this cone, the obstacle is not convex but it is of class $C^{2,1/2}$, as it coincides with the graph of $z\mapsto -|\re z^{5/2}|$.}
\label{fig:artificial}
\end{figure}

\medskip
We will now show that the existence of rank zero points on the free boundaries is substantially simpler if we allow ourselves to consider non-convex obstacles. 
In fact, we can construct an artificial obstacle that touches the image of the maps $u_k$ defined in \eqref{eq:defuk} over a cone with vertex at 0, cf.\ Figure \ref{fig:artificial}. This example has the caveat that the obstacle is not convex, but on the other hand it has two advantages:
\begin{enumerate}
    \item it shows that the bound $\dim_{\cH}(\cB_0(u) \cap \partial u^{-1}(M))\leq n-2$, which follows from Proposition \ref{prop:dimB0}, is in general optimal;
    \item it shows that branch points can lead to cone singularities on the free boundaries; such singularities do not occur in the classical obstacle problem \cite{Caff98-revisited}.
\end{enumerate}
In the next statement, $[\beta]$ denotes the greatest integer less than or equal to $\beta$. 

\begin{prop}\label{prop:artificialobstacle}
    For each odd integer $k\geq 3$ and each $\alpha \in(0,\frac 1 2)$, there is a $C^{[k/2],\alpha}$-domain $M\subset \R^3$ with compact  complement such that the following holds: the analytic map $u_k$, defined in \eqref{eq:defuk}, satisfies 
    $u_k(\overline {\R^2})\subset \overline M$ and  
 $\partial u^{-1}(M)$ is a cone with a true branch point at its vertex.
\end{prop}

\begin{proof}
We work with polar coordinates $z=re^{i \theta}$ with $(r,\theta)\in (0,\infty)\times[-\pi,\pi]$.
Take a small angle $\theta_k \in (0,\frac{\pi}{8})$ such that $\cos(\frac{k}{2}\theta_k)\geq \frac{1}{2}$, and fix  $\lambda \geq 1$ to be chosen later. Let

\begin{align*}
    C&:=\{r e^{i \theta}: r>0, |\theta|\in [0,\theta_k]\},\\
    E_1&:=\{r e^{i \theta}: r>0, |\theta|\in (\theta_k,\tfrac \pi 2 + \theta_k]\},\\
    E_2&:=\{r e^{i \theta}: r>0, |\theta|\in (\tfrac \pi 2 + \theta_k,\tfrac{3 \pi}{4}]\},\\
    F&:=\{r e^{i \theta}: r>0, |\theta|\in (\tfrac{3\pi}{4},\pi]\}.
\end{align*}
To define the obstacle, we set
$$\vp(z):=\begin{cases}
-\re z^{k/2} & \text{if } z\in C,\\
p_k(\pi(z))-\lambda \dist(z,\partial C)^{[k/2]+\alpha} & \text{if } z\in E_1\cup F,
\end{cases}$$
where $p_k(z)$ is the $[\frac{k}{2}]$-th order Taylor polynomial of $z\mapsto -\re z^{k/2}$ at the origin, and $\pi\colon \R^2\to \overline C$ is the nearest point projection onto the convex set $\overline C$. Note that $\vp\in C_\loc^{[k/2],\alpha}(C\cup E_1\cup F)$, and  $\vp(z) = -\lambda|z|^{[k/2]+\alpha}$ for $z\in F$. Now, one may easily extend $\vp$ to $E_2$ so that $\vp\in C_\loc^{[k/2],\alpha}(\R^2)$. Finally, we let $M$ be the epigraph of $\varphi$, namely $$M=\{(z,x_3)\in \R^2\times \R: \varphi(z)< x_3\}.$$
Note that the image of $u_k$ corresponds to the union of the graphs of the functions $z\mapsto \pm|\re z^{k/2}|$.  Selecting $\lambda > 1$ large, corresponding to the choice of $\alpha\in(0,\frac{1}{2})$, we observe that $\partial M$ intersects the graph of $-|\re z^{k/2}|$ precisely on $C$, and lies below it otherwise.  Thus, $u(\R^2)\subset \overline{M}$ with $u^{-1}(\partial M)= C$.  
\end{proof}

In the next remark, we observe another interesting aspect of the example constructed in Proposition \ref{prop:artificialobstacle}.

\begin{rem}
Consider any integer $k\geq 5$; this restriction ensures that $M$ is at least $C^2$. Note that the map $u_k-u_k(0)=u_k$ vanishes to order $2$ at the origin, since we can compute, in polar coordinates,
$$|u_k(re^{i\theta})|=r^2 \sqrt{ 1+ r^{2k-4} \cos^2(k\theta)}.$$
However, if $k\geq 5$ is odd, the distance $\rho\circ u_k$ of $u_k$ to the obstacle vanishes at least to order $k-1+\alpha >  4$, since
$$u_k(r e^{i \pi/2}) = -r^2(1,0,0)\quad \implies \quad 
\rho\circ u_k(r e^{i\pi/2}) = \lambda\, (r^2)^{[k/2]+\alpha}= \lambda r^{k-1+\alpha}.$$ This shows that, although $u_k-u_k(0)$ vanishes to a fixed order at the origin, $\rho \circ u_k$ can vanish to arbitrarily high order there. 
\end{rem}

The above example is only possible because the obstacle $M^c$ in Figure \ref{fig:artificial} is not uniformly convex. In fact, the next proposition, which builds on the results of Sections \ref{sec:Le} and \ref{sec:reg} and will be useful in Section \ref{sec:axial-bp}, shows that, 
 for uniformly convex obstacles, $\rho \circ u$ vanishes twice as fast as $u-u(x_0)$ at a free boundary point $x_0$; by the above remark, this is false for general obstacles.

\begin{prop}\label{prop:vanishingorders}
    Let $M\subset \R^m$ be a smooth obstacle, and let $u\in W^{1,2}(\Omega;\overline M)$ be a weakly constraint map such that $u \in C(B;\cN(\partial M))$ in a ball $B\subset\Omega$. Then for every $x_0\in B\cap \partial u^{-1}(M)$, there exists $k\in\N$ such that $u - u(x_0)$ vanishes exactly to order $k$ at $x_0$, and $\rho\circ u$ vanishes at least to order $2k$ at $x_0$. In addition, if all the principal curvatures of $\partial M \cap u(B)$ are uniformly positive, such as when $M^c$ is uniformly convex, then $\rho\circ u$ vanishes exactly to order $2k$ at $x_0$. 
\end{prop}

\begin{proof}
Since $u\in C(B;\cN(\partial M))$, it follows from \cite[Theorem A.5]{FKS} that $u\in C_\loc^{1,1}(B)$. By \eqref{eq:main-sys-re}, we can apply the classical results concerning the frequency of solutions to elliptic systems (cf.\ Appendix \ref{app:freq}). 
This yields the existence of $k \in \mathbb N$ and a non-trivial (vector-valued) harmonic polynomial $P_k$ of degree $k$, which is homogeneous around $x_0$, and a radius $r_0 > 0$ small (such that $B_{4r_0}(x_0)\subset B$),  such that 
\begin{equation}\label{eq:u-Ck}
    \|D^j(u - u(x_0) - P_k)\|_{L^\infty(B_r(x_0))} \leq c_\sigma r^{k - j + \sigma},\quad j\in\{0,1\},
\end{equation}
for every $r\in(0,4r_0)$ and for each $\sigma \in (0,1)$, where $c_\sigma > 1$ is independent of $r$.
In particular $u-u(x_0)$ vanishes to order $k\in \N$, proving the first assertion of this proposition. Moreover, replacing $r_0$ with a smaller radius if needed, we also obtain 
        \begin{equation}
            \label{eq:freqbdd}
            N(u, x,r)\leq \frac 3 2 k,
        \end{equation}
        for all $r\in(0,r_0)$ and any $x\in B_{r_0}(x_0)$. Here $N(u,x,r)$ is the frequency formula defined in \eqref{freq}.

We now turn our attention to $\rho \circ u.$ 
Note that, as a consequence of \eqref{eq:u-Ck}, we have
\begin{equation}\label{eq:Du-Ck1}
    \| Du \|_{L^\infty(B_r(x_0))} \leq cr^{k-1}
\end{equation}
for every $r\in(0,r_0)$. 
Also, thanks to \cite[Lemma A.3]{FKS}, \eqref{eq:dist-pde} holds in $B_{4r_0}(x_0)$. Therefore, by \eqref{eq:dist-pde}, \eqref{eq:xitop}, \eqref{eq:hessrho}, Corollary \ref{cor:dist}, and \eqref{eq:Du-Ck1}, we get
    \begin{equation}\label{eq:ru-pde}
    \|\Delta(\rho \circ u)\|_{L^\infty(B_r(x_0))} \leq cr^{2k-2}  
    \end{equation}
    for every $r\in (0,r_0)$. As $\rho \circ u\geq 0$ in $\Omega$, a standard blow-up argument (see \cite{Petrosyan2012}) yields a constant $c_0$, which may depend only $n$ and $k$, such that 
    \begin{equation}\label{eq:ru-C2k}
    \| \rho \circ u \|_{L^\infty(B_r(x_0))} \leq c_0cr^{2k},   
    \end{equation}
    for every $r\in(0,r_0)$. This proves that $\rho \circ u$ vanishes at least to order $2k$ at $x_0$. Moreover, it follows from \eqref{eq:decomp} that $\Pi\circ u - u(x_0)$ vanishes at $x_0$ exactly to order $k$. Therefore, the first part of the proposition is proved.

Assume now that all the principal curvatures of $\partial M\cap u(B)$ are uniformly positive.
Taking a smaller $r_0$ if necessary, we obtain a constant $\kappa > 0$ such that 
\begin{equation}\label{eq:Hessrho-uB}
\Hess\rho(\xi,\xi) \geq \kappa |\xi^\top|^2\quad\text{in }\overline M\cap u(B_{4r_0}(x_0)),   
\end{equation}
cf.\ \eqref{eq:Hessrho}.
The rest of the proof is similar in spirit to that of Lemma \ref{lem:reg}, in that we will show that there is a small $\eta > 0$ such that if $\rho\circ u < \eta r^{2k}$ in $B_{4r}(x_0)$ for some $r > 0$, then $\rho\circ u = 0$ in $B_r(x_0)$. Once this claim is proved, from our choice of $x_0\in B\cap\partial u^{-1}(M)$ we deduce that 
$$
\sup_{B_r(x_0)} (\rho\circ u) \geq \eta r^{2k}
$$ 
for all $r\in(0,r_0/4)$, as desired. 

For the moment, let us fix a scale $r\in(0,r_0/4)$, and consider the rescaling
$$
\tilde u_{x_0,r}(x) := \frac{1}{r^k} (u(x_0 + rx) - u(x_0)).
$$
Note from \eqref{eq:main-sys} that \eqref{eq:main-sys-re} holds with $|a_{ij}^\alpha|\leq |Du|$ a.e.\ in $\Omega$. Thus, by \eqref{eq:main-sys-re}, 
$$
\Delta \tilde u_{x_0,r} = \tilde a_{ij}^\alpha D_\alpha \tilde u_{x_0,r}^j \quad \text{a.e.\ in }B_4,
$$
where, by \eqref{eq:Du-Ck1},  all coefficients satisfy
$$
\|\tilde a_{ij}^\alpha \|_{L^\infty(B_4)} \leq \frac{1}{r^{k-1}}\| a_{ij}^\alpha\|_{L^\infty(B_{4r}(x_0))} \leq \frac{c}{r^{k-1}}\|Du\|_{L^\infty(B_{4r}(x_0))} \leq c.
$$
By \eqref{eq:freqbdd}, we also have 
$$
N(\tilde u_{x_0,r},0,4)  = N(u , x_0,4r) \leq \frac{3}{2}k.
$$
Hence, Theorem \ref{thm:Le-reg} yields 
        \begin{equation}
            \label{eq:DPi-Le}
            \frac{1}{r^{n - \gamma(k-1)}}\int_{B_{2r}(x_0)} |Du|^{-\gamma} \,dx = \int_{B_1} |D\tilde u_{x_0,r}|^{-\gamma} dx \leq c,
        \end{equation}
        for some constants $\gamma>0$ and $c> 1$, independent of $r$. 
        
        We now argue as in the proof of Lemma \ref{lem:reg}. Fix a pair $1\leq t\leq T\leq 2$. We observe from \eqref{eq:Hessrho-uB} that analogously to \eqref{eq:key1} we have 
$$
\frac{c}{(T-t)^2r^2}\int_{B_{Tr}(x_0)} (\rho\circ u)\, dx\geq c_0 \int_{B_{tr}(x_0)\cap u^{-1}(M)} |Du|^2\,dx,
$$
with some $c>1$ large and $0 < c_0 < 1$ small now depending on $\kappa$ above, but independent of $s$, $t$, and $r$. Using \eqref{eq:DPi-Le} in place of \eqref{eq:key2}, we deduce from the generalized H\"older inequality, just as in \eqref{eq:key3} and \eqref{eq:key4}, that 
\begin{equation}\label{eq:key-re}
|B_{tr}(x_0)\cap u^{-1}(M)|^{\frac{2}{p}} \leq \frac{c^{2/\gamma}}{r^{2k-2}}\int_{B_{Tr}(x_0)\cap u^{-1}(M)} |Du|^2\,dx \leq \frac{c}{(T-t)^2 r^{2k}}\int_{B_{Tr}(x_0)} (\rho\circ u) \,dx,
\end{equation}
where $p> 1$ is chosen such that $\frac{1}{p} = \frac{1}{2} + \frac{1}{\gamma}$. Since the last inequality holds for all $1\leq t< T\leq 2$, we can deduce from Lemma \ref{lem:degiorgi} that
\begin{equation}\label{eq:key-re2}
    \sup_{B_{2r}(x_0)} (\rho\circ u) < \eta r^{2k}\quad\implies\quad \rho\circ u = 0\quad\text{in }B_r(x_0),
\end{equation}
for some sufficiently small $\eta > 0$. Here the smallness parameter $\eta$ can be chosen independently of $r$, as all constants appearing in \eqref{eq:key-re} are independent of $r$ (nor of $t$ and $T$). The implication of \eqref{eq:key-re2} leads us to $x_0\not\in B\cap \partial u^{-1}(M)$, which gives us a contradiction. Thus, we must have 
$$
\sup_{B_r(x_0)} (\rho \circ u) \geq \eta r^{2k}. 
$$
This proves that the vanishing order of $\rho\circ u$ is at most $2k$, provided that all principal curvatures of $\partial M\cap u(B)$ are uniformly positive, as desired. 
\end{proof}


\section{Axially symmetric constraint maps: regularity theory}\label{sec:axial-reg}

In this section, we initiate the study of Dirichlet energy minimizers in the class of axially symmetric maps that are constrained to take values outside the unit ball.

The regularity theory of such maps is parallel to that of minimizing constraint maps, which we developed in the last sections, but it is independent of it. Indeed, as we will see in Corollary \ref{cor:not-energy-min} below, minimizers in axially symmetric classes are \textit{not} minimizing constraint maps in general. Nonetheless, in a neighborhood of any regular point, minimizers in axially symmetric classes are also minimizers in the full Sobolev space of admissible maps, cf.\ Appendix \ref{app:unique}.

This section is divided into three parts. In Subsection \ref{sec:defs} we define the basic classes of maps under consideration. In Subsection \ref{sec:propsax} we develop the partial regularity theory for minimizers in axially symmetric classes; the results of this subsection are the analogue of the classical results for minimizing constraint maps presented in Section \ref{sec:known}. Finally, in Subsection \ref{sec:Le-ax}, we prove the main results of this section, Theorems \ref{thm:Le-ax} and \ref{thm:reg-ax}, which are the analogue of Theorems \ref{thm:Le} and \ref{thm:uni-reg}, respectively, in the axially symmetric setting. We emphasize that the proof of the $L^{-\gamma}$ estimate in the axially symmetric setting is rather different from the one presented in Section \ref{sec:Le}, and we believe it is of independent interest.

\subsection{Setup and definitions}\label{sec:defs}

Throughout this section, we will use the following:

\begin{notation}\label{not:axial}
    Set ${e}_3 := (0,0,1)$. We use cylindrical coordinates $(r,\theta,z)$ in the domain. We write $\mathbb D=\{(r,z): r^2 + z^2 < 1, r \geq 0\}$ for the half-disk in the $(r,z)$-plane. We shall write $D' = (\partial_r,\partial_z)$ and define $\Delta' := \partial_{rr} + \partial_{zz} + \frac{1}{r}\partial_r$, so that $\Delta= \Delta' + \frac{1}{r^2} \partial_{\theta\theta}$  is the usual Laplacian in $\R^3$. 
\end{notation}

We now define axially symmetric maps, and the corresponding energy minimizers.

\begin{defn}\label{def:axial}
Fix $k\in \N$.
\begin{enumerate}
    \item\label{it:axial} A measurable map $u \colon \mb B^3\to \R^3$ is said to be {\normalfont $k$-axially symmetric} if it is of the form
    \begin{equation}
        \label{eq:axial}
        u = \varrho\,(\sin\varphi\cos k\theta,\sin\varphi\sin k\theta,\cos\varphi),
    \end{equation}
    for some pair $(\varrho,\varphi):\D\to [0,\infty)\times[0,\pi]$.
\item A map $u\in W^{1,2}(\mb B^3;\R^3\backslash \mb B^3)$ is said to be {\normalfont energy minimizing among $k$-axially symmetric maps} if $u$ is itself $k$-axially symmetric and it satisfies 
$$
\cE[u]\leq\cE[\tilde u]
$$
for all $k$-axially symmetric maps $\tilde u\in W^{1,2}_u(\mb B^3;\R^3\backslash \mb B^3)$. 
\end{enumerate}
\end{defn}

\begin{rem}
    Note that \eqref{eq:axial} implies 
    \begin{equation}
        \label{eq:commrots}
    u\circ R_\theta = R_{k\theta}\circ u,
    \end{equation}
    for every $\theta\in[0,2\pi]$, where $R_\theta$ denotes the rotation about the $z$-axis by angle $\theta$. However, \eqref{eq:commrots} is not sufficient to deduce the representation \eqref{eq:axial}. In fact, \eqref{eq:commrots} allows additional rotation of each horizontal circle in $\R^3$ (i.e.\ a circle of form $\{r = r_0$, $z = z_0\}$) by an amount depending on $(r_0,z_0)$. However, minimizers of the energy among all maps that verify \eqref{eq:commrots} subject to a $k$-axially symmetric boundary condition, in the sense of \eqref{eq:axial}, turn out to be $k$-axially symmetric as well.  We refer the reader to \cite[Lemma 2.8]{GuerraKochLindberg2020c} for a similar argument in a slightly different context.
\end{rem}

In what follows, given any pair of function $(\varrho,\vp)\colon \D \to [1,\infty)\times [0,\pi]$, we define its \textit{energy} through
\begin{equation}\label{eq:E-F}
\cF_k[\varrho,\vp] := \iint_{\mathbb D} \bigg( |D'\varrho|^2 + \varrho^2\bigg( |D'\varphi|^2 + \frac{k^2}{r^2}\sin^2\vp\bigg)\bigg) r\,dr\,dz.
\end{equation}
We call $(\varrho,\varphi)$ admissible if in addition $\cF_k[\varrho,\varphi] < +\infty$. 

\begin{rem}\label{rem:E-F}
Note that $u\in W^{1,2}(\mb B^3;\R^3\setminus\mb B^3)$ is $k$-axially symmetric if and only if the associated pair $(\varrho,\varphi)$ is admissible, since 
$$
\cE[u] = \int_{\mb B^3} |Du|^2\,dx = 2\pi \cF_k[\varrho,\varphi]. 
$$
Thus $u \in W^{1,2}(\mb B^3;\R^3\setminus \mb B^3)$ is $\cE$-minimizing among all $k$-axially symmetric maps if and only if the associated pair $(\varrho,\varphi)$ is $\cF_k$-minimizing among all admissible pairs,  i.e.
$$
\cF_k[\varrho,\varphi] \leq \cF_k[\tilde\varrho,\tilde\varphi]
$$
for every admissible pair $(\tilde\varrho,\tilde\varphi)$ satisfying $\supp(\varrho -\tilde\varrho)\cup\supp(\varphi -\tilde\varphi) \Subset \D$. 
\end{rem}

\subsection{Partial regularity theory}
\label{sec:propsax}

We now turn to the partial regularity theory for energy minimizers among $k$-axially symmetric maps.  Throughout this subsection we will always write $u\in W^{1,2}(\mb B^3;\R^3\backslash \mb B^3)$ for a fixed energy minimizer among $k$-axially symmetric maps.

Let us begin with the Euler--Lagrange system for the energy minimizers, and their regularity away from the $z$-axis, i.e.\ away from $\{r=0\}.$

\begin{lem}\label{lem:EL-ax}
Let $u\in W^{1,2}(\mb B^3;\R^3\setminus \mb B^3)$ be an energy minimizer among $k$-axially symmetric maps.
Then $u\in C_\loc^{1,1}(\mb B^3\setminus\{r=0\})$ and it is a weak solution of 
    \begin{equation}\label{eq:main-ax}
    -\Delta u = |Du|^2 u \chi_{\{|u|=1\}}\quad\text{in }\mb B^3.
    \end{equation}
Moreover, $|u|-1 \in W^{1,2}(\mb B^3)$ is nonnegative and weakly subharmonic. 
\end{lem}

\begin{proof}
Throughout the proof we work on a fixed domain $\Omega\Subset \D\setminus\partial\D = \{(r,z): r^2 + z^2 < 1, r >0\}$. Note from Remark \ref{rem:E-F} that $(\varrho,\varphi)$ is a minimizer of the energy $\cF_k$ among all admissible pairs. With $\varphi$ fixed, we see that $\varrho$ minimizes $\cF_k[\cdot,\varphi]$, i.e., $\cF_k[\varrho,\varphi] \leq \cF_k[\tilde\varrho,\varphi]$ for every admissible pair $(\tilde\varrho,\varphi)$ with $\supp(\varrho-\tilde\varrho) \Subset\D$. Note that the coefficient of $\varrho^2$ in $\cF_k[\varrho,\varphi]$ is nonnegative and belongs to $L^1(\Omega)$, due to the constraint $\varrho \geq 1$. Let $\beta\in C^1((-\infty,\infty))$ be such that  $\beta = 0$ in $(-\infty,1)$, $\beta = 1$ in $[2,\infty)$ and $\beta' \geq 0$ in $[0,\infty)$. Then $(\varrho + \e\beta((\varrho-1)/\e)\psi,\varphi)$ is an admissible pair, for all small $\e> 0$, for each $\psi\in C_c^1(\Omega)$. Hence, the minimality of $(\varrho,\varphi)$ implies that  
\begin{equation}\label{eq:vr-subharm}
\Delta'\varrho \geq \varrho \bigg(|D'\varphi|^2 + \frac{k^2}{r^2} \sin^2\varphi\bigg)\chi_{\{\varrho > 1\}} \geq 0 \quad\text{in }\Omega,    
\end{equation}
in the weak sense. Since the coefficients of $\Delta'$ are bounded away from zero and infinity in a neighborhood of $\Omega$, the local $L^\infty$-estimate \cite[Theorem 8.17]{GT} for weak subsolutions applies and so $\varrho \in L_\loc^\infty(\Omega)$.

Let us now consider the variation of $\cF_k[\varrho,\cdot]$ with $\varrho$ fixed. Although an admissible $\tilde\varphi$ is constrained to taking values in $[0,\pi]$, we can actually consider \textit{unconstrained} variations  in deriving the Euler--Lagrange equation for $\varphi$. Indeed, given any $\tilde\varphi \colon \D\to \R$ satisfying $\cF_k[\varrho,\tilde\varphi] < \infty$, the pair $(\varrho,\max\{\min\{\tilde\varphi,\pi\},0\})$ becomes admissible, and 
$\cF_k[\varrho, \max\{\min\{\tilde\varphi,\pi\},0\}] \leq \cF_k[\varrho, \tilde\varphi].$
Thus, the minimality of $(\varrho,\varphi)$ implies $\cF_k[\varrho, \varphi] \leq \cF_k[\varrho, \varphi + \e\psi]$ for any $\e> 0$ small, for each $\psi\in C_c^\infty(\Omega)$, from which we deduce that 
\begin{equation}\label{eq:vp-pde}
D'\cdot (r\varrho^2 D'\varphi)  = \frac{k^2\varrho^2}{2r}\sin 2\vp\quad\text{in }\Omega
\end{equation}
in the weak sense. 
Note from our earlier discussion that $\varrho\in L_\loc^\infty(\Omega)$. Moreover, since $\cF_k[\varrho,\varphi] < + \infty$ and $\varrho\geq 1$ in $\D$, the right-hand side of \eqref{eq:vp-pde} belongs to $L_\loc^2(\Omega)$. Thus, the higher integrability estimate \cite[Theorem 2]{Meyers1963} applies to \eqref{eq:vp-pde} and we see that $\vp \in W_\loc^{1,p}(\Omega)$ for some exponent $p = p(\Omega) > 2$. Since $\Omega$ is a two-dimensional domain, the Sobolev embedding implies $\varphi\in C_\loc^{0,\bar\sigma}(\Omega)$, with $\bar\sigma := 1-2/p$. 

Let us return to $\varrho$. Noting that $(\varrho +\e\psi,\varphi)$ is an admissible pair for any $\e>0$ small, for each $\psi \in C_c^\infty(\Omega)$ nonnegative, we obtain from the minimality of $(\varrho,\varphi)$ that  
\begin{equation}\label{eq:vr-sup}
\Delta'\varrho\leq \varrho \bigg(|D'\varphi|^2 + \frac{k^2}{r^2} \sin^2\varphi\bigg)\quad\text{in }\Omega,
\end{equation}
in the weak sense. Since $\varphi \in W_\loc^{1,p}(\Omega)$ with some $p > 2$ and $\varrho\in L_\loc^\infty(\Omega)$, the right-hand side of \eqref{eq:vr-sup} belongs to $L_\loc^{p/2}(\Omega)$. By \eqref{eq:vr-sup} and \eqref{eq:vr-subharm}, we obtain from \cite[Theorem 9.11]{GT} that $\varrho \in W_\loc^{2,p/2}(\Omega)$. 
Thus, since once again $\Omega$ is a two dimensional domain, the Sobolev embedding implies $\varrho \in C_\loc^{0,2\bar\sigma}(\Omega)$ (recall $\bar\sigma = 1-2/p$). Hence $(\varrho +\e\psi,\varphi)$ is an admissible pair for any $\e> 0$ small and for each $\psi\in C_c^\infty(\Omega\cap\{\varrho > 1\})$, i.e.\ $\psi$ can change sign in the region where $\varrho>1$. Thus, the minimality of $(\varrho,\varphi)$ yields equality in \eqref{eq:vr-sup} in $\Omega\cap\{\varrho>1\}$ in the weak sense. However, as  $\varrho$ has weak second derivatives in $\Omega$, we also have $\Delta'\varrho = 0$ a.e.\ in $\Omega\cap\{\varrho = 0\}$. The last two observations combined yield  
\begin{equation}\label{eq:vr-pde-E}
\Delta'\varrho = \varrho \bigg(|D'\varphi|^2 + \frac{k^2}{r^2} \sin^2\varphi\bigg)\chi_{\{\varrho > 1\}} \quad\text{a.e.\ in }\Omega. 
\end{equation} 

We now return once again to $\vp$.  Since $\varrho\in C_\loc^{0,2\bar\sigma}(\Omega)$, we can apply the H\"older estimate for the first derivatives \cite[Theorem 8.33 and Corollary 8.35]{GT} to \eqref{eq:vp-pde} and deduce that $\varphi \in C_\loc^{1,2\bar\sigma}(\Omega)$. This in turn yields that the right-hand side of \eqref{eq:vr-pde-E} belongs to $L_\loc^\infty(\Omega)\cap C_\loc^{0,2\bar\sigma}(\Omega\setminus\partial\{\varrho>1\})$, so by the interior $C^{1,1}$-estimate \cite[Theorem 2.5]{Caffarelli1980} for the classical obstacle problem we obtain $\varrho \in C_\loc^{1,1}(\Omega)$.  Plugging this information into \eqref{eq:vp-pde}, we deduce again from the higher-order regularity theory \cite[Theorem 9.19]{GT} that $\vp \in C_\loc^{2,\sigma}(\Omega)$ for every $\sigma\in(0,1)$. In particular, the Euler--Lagrange equation \eqref{eq:vp-pde} for $\varphi$ holds in $\Omega$ in the classical sense. 

Since $\Omega$ was an arbitrary domain compactly contained in $\D\setminus\partial\D$, the above observations show that $(\varrho,\varphi) \in C_\loc^{1,1}(\D\setminus\partial\D)\times C_\loc^{2,\sigma}(\D\setminus\partial\D)$. Thus, in light of \eqref{eq:axial}, we must have $u\in C_\loc^{1,1}(\mb B^3\setminus\{r=0\})$. Moreover, as the Euler--Lagrange equations \eqref{eq:vp-pde} and \eqref{eq:vr-pde-E} hold a.e.\ in $\D\setminus\partial\D$, a direct computation yields that the Euler--Lagrange system \eqref{eq:main-ax} holds in $\mb B^3\setminus\{r=0\}$ in the strong sense. Since $u\in W^{1,2}(\mb B^3;\R^3)$ and the $z$-axis has zero hamronic capacity, we deduce that \eqref{eq:main-ax} holds in the weak sense. Finally, the fact that $|u| - 1$ is weakly subharmonic in $\mb B^3$ follows either directly from \eqref{eq:vr-subharm} or from Lemma \ref{lem:dist-subharm}.
\end{proof}

\begin{rem}
    It is an interesting problem to decide whether, in the setting of Lemma \ref{lem:EL-ax}, we have $\varphi \in C^{2,1}_{\loc}(\mathbb D)$, see \cite[Theorem 1]{ABB1} for the one-dimensional case and  \cite[Corollary 2.7]{FKS} for the two-dimensional case.
\end{rem}

Lemma \ref{lem:EL-ax} shows that $u$ is regular away from the $z$-axis. We now prove the following $\e$-regularity theorem on the $z$-axis, which is the analogue of Theorem \ref{thm:e-reg} and \cite[Lemma 4.1]{HKL} in our setting:

\begin{thm}\label{thm:e-reg-ax}
    There are absolute constants $\e_0 > 0$ and $c>1$, such that if  $B_{4s}(z_0e_3)\subset \mb B^3$ and $E(u,z_0e_3,2s) \leq \e_0^2$, then $u \in C^{1,1}(B_s(z_0e_3))$ and 
    $$
    \| D^j u\|_{L^\infty(B_s(z_0e_3))} \leq \frac{c}{s^j},\quad j=1,2.
    $$
\end{thm}

\begin{proof}
The idea follows essentially from \cite[Lemma 4.1]{HKL}. Translating and rescaling $u$ in the domain, we can assume that $z_0=0$, $s=1$, and $u$ is defined in $4\mb B^3$. We first note that
$$|D'(\cos \vp)| \leq \lvert \sin \vp\rvert (|\partial_r \varphi|+ |\partial_z \vp|) \leq \frac 1 2\bigg(|D'\vp|^2 + \frac{\sin^2 \vp}{r^2}\bigg) r $$
and so, since $\varrho\geq 1$ and $k\geq 1$, we have
\begin{equation}
\label{eq:boundcosphi}
\int_{\mb D} |D'(\cos \vp)| \, d x \leq \cF_k[\varrho, \vp] = \frac{E(u,0,1)}{2\pi} \leq \frac{\e_0^2}{2\pi}.
\end{equation}
By Fubini's Theorem and the fact that $\cF_k[\varrho,\vp]<\infty$, for a.e.\ $R>0$ we have that $\vp|_{S_{R}}$ is continuous, where we write $S_{R} = \{(r,z)\in \mb D:r^2+z^2=R^2\}$.  Similarly, using \eqref{eq:boundcosphi}, we can further choose $R$ so that
\begin{equation}
\label{eq:oscphi}
\osc_{S_R} \cos \vp \leq c\,\e_0^2,
\end{equation}
for a universal constant $c.$ Since $\cos \vp(0,[-1,1])\subseteq \{0,\pi\}$, we can assume without loss of generality that $\vp(0,R)=\vp(0,-R)=0$, the case where $\vp(0,\pm R)=\pi$ being totally analogous. We thus see from \eqref{eq:oscphi} that, if $\e_0$ is chosen small enough, then the image of $\mb B^3$ under $u$ lies in an upper half space, i.e.
$$u(\mb B^3)\subset (\R^3\backslash \mb B^3) \cap \{x_3>\tfrac 1 2\}.$$
In other words, $u(\mb B^3)$ lies in a star-shaped region and so, by \cite[Theorem 1]{F1} (see also Appendix \ref{app:unique} and in particular Theorem \ref{thm:fuchs} therein), we conclude that $u\in C^{1,\sigma}(\mb B^3)$ for any $\sigma\in (0,1)$. The sharp regularity $u\in C^{1,1}(\mb B^3)$ then follows from the standard theory for the scalar obstacle problem, as in \cite{FKS}.
\end{proof}

We now develop further properties of $u$. In analogy with Lemma \ref{lem:monot}, we have the following monotonicity formula:

\begin{lem}\label{lem:monot-ax}
     Let $z_0\in (-1,1)$ and $0<t<s < 1- |z_0|$ be given. Then 
     $$E(u,z_0{e_3},t) \leq E(u,z_0{e_3},s),$$ with equality if and only if $u$ is $0$-homogeneous about $z_0{e_3}$ in $B_s(z_0{e_3})\setminus B_t(z_0{e_3})$.
\end{lem}

\begin{proof}
    Note that the minimizing property over the class of $k$-axially symmetric maps yields the stationarity over domain variations for balls centered on the $z$-axis. The rest of the proof is identical to the standard one, see e.g.\ \cite[Appendix A]{FGKS}. 
\end{proof}

As in Section \ref{sec:known}, Lemma \ref{lem:monot-ax} allows us to define the normalized energy density $E(u,z_0 e_3, 0^+)$ through the limit $E(u,z_0 e_3, r)$ as $r\searrow 0.$

We now turn to the compactness properties of minimizers in the axially symmetric class, which is an analogue of Lemma \ref{lem:strong} in the axially symmetric setting.

\begin{lem}\label{lem:strong-ax}
    Let $\{u_i\}_{i=1}^\infty\subset W^{1,2}(\mb B^3;\R^3\backslash \mb B^3)$ be a bounded sequence of energy minimizers among $k$-axially symmetric maps. Then there is a map $u\in W^{1,2}(\mb B^3;\R^3\backslash \mb B^3)$, which is energy minimizer among $k$-axially symmetric maps, such that $u_i\to u$ strongly in $W_\loc^{1,2}(\mb B^3;\R^3)$ along a subsequence.  
\end{lem}

As in the classical case \cite{L}, the key idea is to consider a suitable axially symmetric extension of the boundary values, see  \cite[Lemma 4.1]{G} and also \cite[Proof of Theorem 4.2]{HKL} in the setting of axially symmetric harmonic maps. Although the extension here gets slightly more involved due to the presence of the non-constant radial component $\varrho$, most of the argument is similar to the aforementioned literature. To maintain the technicalities to a minimum and keep the reading flow, we shall present the proof in Appendix \ref{app:strong-ax}.

As a corollary of Lemma \ref{lem:EL-ax} and Theorem \ref{thm:e-reg-ax}, and using Lemmas  \ref{lem:monot-ax} and \ref{lem:strong-ax}, we can run the usual dimension reduction argument to deduce the following:

\begin{cor}\label{cor:dim-ax}
We have $u\in C^{1,1}_\loc(\Omega\setminus \Sigma(u))$.    The singular set $\Sigma(u)$ consists of discrete points on the $z$-axis, which are characterized by $E(u,z_0 e_3, 0^+)\geq \e_0^2$.
\end{cor}

From Theorem \ref{thm:e-reg-ax} and Lemma \ref{lem:strong-ax}, arguing just as in Section \ref{sec:dist}, we also obtain:

\begin{cor}\label{cor:subharm-ax}
    The function $|u|-1$ is nonegative, subharmonic and  continuous on $\mb B^3$, with a  universal modulus of continuity.   
\end{cor}

\subsection{Regularity near free boundaries}\label{sec:Le-ax}

We now turn to the regularity theory of energy minimizers among $k$-axially symmetric maps near their free boundaries. A crucial ingredient in our analysis is the uniqueness and structure of blow-ups in the axially symmetric class, see \cite{BCL} and \cite[\S 6.3]{HKL}.  In what follows,  we use the following standard terminology:

\begin{defn}
Given $z_0\in (-1,1),$ we call $\phi_{z_0}$ a {\normalfont tangent map} of $u$ at $z_0{e_3}$ if it is a weak $W_\loc^{1,2}$-limit of a sequence $u_{r_j,z_0e_3}:= u(z_0e_3 + r_j\cdot)$, where $r_j\to 0$ as $j\to \infty$. 
\end{defn}
 Combining Lemmas \ref{lem:monot-ax} and \ref{lem:strong-ax} with Corollary \ref{cor:dim-ax},  it follows that $z_0e\in \Sigma(u)$ if and only if there is a non-constant tangent map $\phi_{z_0}$ of $u$ at $z_0e$, which is necessarily 0-homogeneous and minimizing among $k$-axially symmetric maps. In our setting, there are only four different non-constant tangent maps:

\begin{prop}\label{prop:tan-ax}
Let $u\in W^{1,2}(\mb B^3;\R^3\backslash \mb B^3)$ be an energy minimizer among $k$-axially symmetric maps. 
For every point $z_0{e_3}\in \Sigma(u)$, there exists a unique tangent map $\phi_{z_0}$ which, when restricted to $\mathbb S^2$, is one of the four maps 
$$
\pm S^{-1}\circ \psi_k \circ S \quad \text{ or }\quad \pm S^{-1}\circ \overline{\psi_k} \circ S,
$$
where $\psi_k\colon \C\to \C$ is the map $w\mapsto w^k$, and $S\colon \mathbb S^2\to \C$ is the stereographic projection from $e_3$.
Moreover, there is a modulus of continuity $\omega$, depending only on $k$, $\dist(z_0e_3,\partial\Omega)$, and $\|u\|_{W^{1,2}(\Omega)}$, such that 
$$
|u(x) - \phi_{z_0} (x- z_0e_3)| \leq \omega(|x-z_0e_3|),
$$
for every $x\in\Omega$ with $|x-z_0e_3|\leq \frac{1}{2}\dist(z_0e_3,\partial\Omega)$. 
\end{prop}

\begin{proof}
According to Corollaries \ref{cor:dim-ax} and \ref{cor:subharm-ax},  we have $\Sigma(u)\subset \{|u|=1\}\cap \{r=0\}$. Thus, by Lemmas \ref{lem:monot-ax} and \eqref{lem:strong-ax}, at any point in $\Sigma(u)$, all tangent maps of $u$ at that point are 0-homogeneous harmonic maps into $\mb S^2$, which further minimize the energy among $k$-axially symmetric harmonic maps. The first claim now follows from the classification of such maps from \cite[\S 6.3]{HKL}. The second claim, concerning the existence of a modulus of continuity, is a direct consequence of the compactness of minimizing maps provided by Lemma \ref{lem:strong-ax},  and the fact that there is only a \textit{finite} number of tangent maps. Indeed, since the function $r\mapsto u_{r,z_0 e_3}$ is continuous into $W^{1,2}(\mb B^3)$, the closest tangent map to $u_{r,z_0 e_3}$ varies continuously with $r$. But since there is only a finite number of such maps, such a map must eventually always be the same for $r$ sufficiently small.
\end{proof}

\begin{cor}\label{cor:not-energy-min}
Let $k\geq 2$ and let $u\in W^{1,2}(\mb B^3;\R^3\setminus \mb B^3)$ be an energy minimizer among $k$-axially symmetric maps. If $z_0 e_3 \in \Sigma(u)$ then $u$ is not energy minimizing in any neighborhood of $z_0 e_3.$
\end{cor}

\begin{proof}
    As described before the proof of Proposition \ref{prop:tan-ax}, to  any point $z_0e_3\in \Sigma(u)$ one can associate tangent maps $\phi_{z_0}$ of $u$ at $z_0e_3$, which in particular are non-constant 0-homogeneous $k$-axially symmetric maps. If $u$ was minimizing near $z_0 e_3$ (in the full Sobolev space of admissible maps), then according to Lemma \ref{lem:strong} so would $\phi_{z_0}.$ However, by $k$-axial symmetry it is easy to see that 
    $\phi_{z_0}\colon \mb S^2\to \mb S^2$ has $\deg(\phi_{z_0})=k$. Thus, according to \cite[Theorem 7.4]{BCL}, $\phi_{z_0}$ is not a minimizer.
\end{proof}

We can now prove the first main result of this section, which is an extension of Theorem \ref{thm:Le} to the setting of energy minimizers among $k$-axially symmetric maps.

\begin{thm}\label{thm:Le-ax}
Let $u\in W^{1,2}(\mb B^3;\R^3\backslash \mb B^3)$ be energy minimizing among $k$-axially symmetric maps. There are positive constants $s,\gamma$, and $c$, which depend only on $k$ and $\|u\|_{W^{1,2}(\mb B^3)}$, such that if $z_0e_3\in\Sigma(u)$ for some $z_0\in (-1,1)$, then 
    $$
    \int_{B_s(z_0e_3)} |Du|^{-\gamma}\,dx \leq c. 
    $$
\end{thm}

\begin{proof}
Let us write $\Lambda:=\|u \|_{W^{1,2}(\mb B^3)}$ for simplicity. By pre-composing $u$ with a translation along the $z$-axis, which of course preserves the $k$-axial symmetry, we can assume that $z_0=0$.
Furthermore, by Lemma \ref{lem:strong-ax} and Corollary \ref{cor:dim-ax}, we deduce that $B_{\bar s}\cap\Sigma(u) = \{0\}$, for some $\bar s =\bar s(k,\Lambda)$. 

Consider the $0$-homogeneous rescalings $u_s(y) := u(sy)$, and let $\phi_0$ be as in Proposition \ref{prop:tan-ax}. This proposition also yields a modulus of continuity $\omega$, depending only on $k$ and $\Lambda$, such that $\|u_s - \phi_0\|_{L^{\infty}(\mb S^2)} \leq \omega(s)$ on $\mb S^2$ for all $s < \frac{1}{2}$. Since $\phi_0|_{\mb S^2}$ is smooth and $u_s|_{\mb S^2}$ is uniformly close to it, $\varepsilon$-regularity implies that it must also be continuous and therefore $C^{1,1}$. Hence, given $\alpha \in (0,1)$,  by interpolation we can find a new modulus of continuity $\omega$, which still depends only on $k$ and $\Lambda$, such that
$$
\|u_s -\phi_0\|_{C^{1,\alpha}(\mb S^2)} \leq \omega(s).
$$ 
 In light of Proposition \ref{prop:tan-ax}, we can find $\e,\delta > 0$, depending only on $k$, such that
$$
|D\phi_0| \geq 2\delta \quad\text{on }\mb S^2\setminus B_\e(\pm e_3).
$$
Hence, combining the previous two estimates, we can take $\bar s$ smaller if necessary, without changing its dependence on $k$ and $\Lambda$, in such a way that 
$$
|Du_s| \geq \delta\quad\text{on }\mb S^2\setminus B_\e(\pm e_3)
$$
for every $s\in(0,\bar s)$. Thus, taking the cone $\cC_\e := \{0\}\cup (\bigcup_{s > 0} B_{s \e}(\pm se_3))$,
$$
\int_{B_{\bar s}\setminus\cC_\e} |Du|^{-1}\,dx = \int_0^{\bar s} s^2 \,ds \int_{\mb S^2\setminus B_\e(\pm e_3)} |Du_s|^{-1} \,d\cH^2 \leq c\delta^{-1}.
$$
Therefore, by Jensen's inequality, it suffices to find some $\gamma = \gamma(k,\Lambda) \in (0,1)$ such that 
$$
\int_{B_{\bar s}\cap \cC_\e} |Du|^{-\gamma}\,dx \leq c. 
$$
We now follow the arguments in Section \ref{sec:Le}. We first claim, analogously to Lemma \ref{lem:freq}, that
\begin{equation}\label{eq:freq-ax}
\sup_{1\leq t\leq 4} N(u_s, \pm te_3, 1) \leq \lambda_0,
\end{equation}
uniformly for all $s \in (0,2\bar s)$, for some $\lambda_0=\lambda_0(k,\Lambda)$, with possibly smaller $\bar s = \bar s(k,\Lambda)$. This claim can be verified by following the lines of the proof of Lemma \ref{lem:freq} without major modifications. We only need to employ Lemmas \ref{lem:monot-ax}, \ref{lem:strong-ax}, and Corollary \ref{cor:dim-ax}  in place of Lemmas \ref{lem:monot}, \ref{lem:strong}, and Corollary \ref{cor:e-reg}, respectively. Note also that $|u_s - \phi_0| \leq \omega(5s)$ in $B_5$ and $|\phi_0| = 1$ replace Corollary \ref{cor:sup} in the proof; alternatively, this could also be inferred from Corollary \ref{cor:subharm-ax}. Apart from these changes the proof goes through unchanged, and so we shall not repeat the details here.

Let us fix an absolute constant $\tau \in [1,4)$ to be determined later. In light of Proposition \ref{prop:tan-ax}, we have 
$$
|D\phi_0| \leq c \quad\text{in }B_{4\tau\e}(\pm te_3)
$$
uniformly for $t\in[1,4]$.  Choosing $\e>0$ smaller if necessary, we also have 
$$
Du_s \to D\phi_0\quad\text{uniformly in }B_{4\tau\e}(\pm te_3)
$$
as $s\to 0$; here again the rate of convergence depends only on $k$ and $\Lambda$. Thus, 
\begin{equation}\label{eq:Dus-Linf}
|Du_s| \leq 2c\quad\text{in }B_{2\tau\e}(\pm te_3),
\end{equation}
uniformly for all $s\in(0,2\bar s)$, with a possibly smaller $\bar s = \bar s(k,\Lambda)$. Thus, in view of \eqref{eq:main-ax} and \eqref{eq:freq-ax}, we can also employ Lemma \ref{lem:freq-monot} to find $\delta,\sigma \in (0,1)$, which depend only on $k$ and $\Lambda$, such that 
\begin{equation}\label{eq:Dus-double}
|\{ |Du_s| \geq \delta \|Du_s\|_{L^\infty(3Q)}\}\cap Q|\geq (1-\sigma)|Q|,
\end{equation}
for every cube $Q$ with $3Q\subset B_{2\tau\e}(\pm t e_3)$. With \eqref{eq:Dus-double} at hand, we can proceed with the Calder\'on-Zygmund cube decomposition argument, just as in the proof of Theorem \ref{thm:Le}, to get that 
\begin{equation}\label{eq:LeBta}
\int_{B_{\tau\e}(\pm te_3)} |Du_s|^{-\gamma}\,dx \leq c,
\end{equation}
uniformly for all $s\in(0,2\bar s)$ and $1\leq t\leq 4$. 

To conclude, we choose $\tau \geq 1$ such that $(B_4\setminus B_1)\cap\cC_\e \subset \cup_{1\leq t\leq 4} B_{\tau \e}(\pm te_3)$. By \eqref{eq:LeBta}, 
\begin{equation}\label{eq:LeCe}
    \int_{(B_4\setminus B_1)\cap\cC_\e} |Du_s|^{-\gamma}\,dx \leq c,
\end{equation}
uniformly for $s\in (0,2\bar s)$. Choose $\bar\ell$ such that $2^{-\bar\ell}\geq \bar s > 2^{-\bar\ell -1}$. As $\bar s$ depends only on $k$ and $\Lambda$, so does $\bar\ell$. Then by \eqref{eq:LeCe}, 
$$
\begin{aligned}
\int_{B_{\bar s}\cap\cC_\e} |Du|^{-\gamma}\,dx & \leq \sum_{\ell = \bar\ell}^\infty  \int_{(B_{2^{-\ell+2}}\setminus B_{2^{-\ell}})\cap\cC_\e} |Du|^{-\gamma}\,dx \leq \sum_{\ell =\bar\ell}^\infty 2^{-\ell(n + \gamma)} c \leq c. 
\end{aligned}
$$
This finishes the proof. 
\end{proof}

Based on Theorem \ref{thm:Le-ax}, and as shown  in Section \ref{sec:reg}, we can establish that energy minimizers among $k$-axially symmetric maps are regular around their free boundaries.

\begin{thm}\label{thm:reg-ax}
    Let $u \in W^{1,2}(\mb B^3;\R^3\backslash \mb B^3)$ be an energy minimizer among $k$-axially symmetric maps. Then for every $x_0 \in \mb B^3$ lying in the closure of $\{|u|>1\}$, $u \in C^{1,1}(B_\delta(x_0))$ with the estimate 
    $$
    \|u \|_{C^{1,1}(B_\delta(x_0))} \leq c,
    $$
    where $\delta > 0$ and $c>1$ depend only on $1-|x_0|$, $k$ and $\|u\|_{W^{1,2}(\mb B^3)}$. 
\end{thm}

 This result follows from Theorem \ref{thm:Le-ax} by using the properties of energy minimizers in the axially symmetric class established in Subsection \ref{sec:propsax},  exactly as in the proof of Theorem \ref{thm:uni-reg}: since the argument is identical, we omit it.

\section{Axially symmetric constraint maps: branch points on free boundaries}\label{sec:axial-bp}

The goal of this section is to prove Theorem \ref{thm:axial}, concerning the structure of the free boundary points in $\cB_0(u)$, for an energy minimizer $u\in W^{1,2}(\mb B^3;\R^3\backslash\mb B^3)$ in the class of $k$-axially symmetric maps. 
Uunless otherwise stated, in what follows we assume that $u$ is such a map and that
$k\geq 2.$ As in \eqref{eq:axial}, we then write
$$u=\varrho\,(\sin \varphi \cos k\theta,\sin \vp \sin k \theta, \cos \vp),$$
for some pair $(\varrho,\varphi)\colon \mb D\to[0,\infty)\times [0,\pi]$, with $\varrho = \varrho (r, z) $ and $\varphi = \varphi (r,z)$. We will also consistently use Notation  \ref{not:axial} and,  to keep the notation as simple as possible, we will write $v \colon \mb B^3\to\Ss^2$ for the projection map of $u$ onto the sphere $\Ss^2$. Specifically,  $v=\Pi\circ u$, in the notation of Section \ref{sec:prelim}:   
\begin{equation}\label{eq:v}
v  := \Pi \circ u = \frac{u}{\varrho} = (\sin\varphi\cos (k\theta),\sin\varphi\sin(k\theta),\cos\varphi)\quad\text{in }\mb B^3. 
\end{equation}
In particular, $|Dv|$ is constant in $\theta$, and 
\begin{equation}\label{eq:b}
|Dv|^2 = (\partial_r\varphi)^2 + (\partial_z\varphi)^2 + \frac{k^2}{r^2}\sin^2\varphi\quad \text{in }\mb B^3\setminus\{r = 0\}. 
\end{equation}
As mentioned in Section \ref{subsec:branch}, our motivation for considering maps $u$ as described above is that any free boundary point on the $z$-axis is necessarily a branch point in $\cB_0(u).$ Indeed, by Theorem \ref{thm:reg-ax}, $u$ is $C^{1,1}$ around the free boundary,  and thus $D\varrho=0$ there, since $\varrho$ attains its minimum on the free boundary. On the other hand,  we have $Dv=0$ at any regular point on the $z$-axis. Indeed,  using the representation \eqref{eq:v}, we see that, for a fixed $z$, the map $f\colon \C \to \C$ defined by $f(w)=f(r e^{i\theta}) :=(v^1+i v^2)(r,\theta,z)$ satisfies
\begin{equation}
\label{eq:Df-complex}
 \partial_{w} f = \bigg(\cos \varphi\, \partial_r \varphi+\frac{k}{r} \sin\varphi \bigg)e^{i (k-1)\theta},\quad \partial_{\bar w} f = \bigg(\cos \varphi\,\partial_r \varphi - \frac{k}{r} \sin\varphi\bigg)e^{i (k+1)\theta}.
\end{equation}
Since $k\geq 2$, both exponentials $e^{i(k\mp 1)\theta}$ cannot be extended continuously to  $w=0$ and so, at any point where $v$  (hence $f$)  is $C^1$, their coefficients in \eqref{eq:Df-complex} must converge to $0$ as $
|w| = r\to 0^+$. Since $|\cos \varphi(0,z)|=1$ and $\sin\varphi(0,z)=0$, we conclude that $Dv$ vanishes at all regular points on the $z$-axis. Thus
$$\partial \{|u|>1\} \cap \{r=0\} \subseteq \cB_0(u).$$
The results of this section can be seen as a much more precise version of the above calculation.  Indeed, we will perform a fine blow-up analysis of $u$ near a free boundary point on the $z$-axis. The structure of $u$ at such a point is quite rich and the next theorem, which is the main result of this section, summarizes the main conclusions of our analysis.

\begin{thm}[Blow-up analysis on the free boundary]\label{thm:bp-precise} Given $k \geq 2$, let  $u\in W^{1,2}_g(\mb B^3;\R^3\backslash \mb B^3)$ be an energy minimizer in the class of $k$-axially symmetric maps. Set $\varrho=|u|$ and $v = u/|u|$ as in \eqref{eq:v}, and assume that $g\in C^2(\mb S^2;\lambda\mb S^2)$ for some $\lambda>1$, and that $g$ is $k$-axially symmetric with $\lvert\deg(g)\rvert=k.$

Let $z_0 e_3\in \partial \{\varrho>1\}$ be an arbitrary free boundary point on the $z$-axis. Then:
\begin{enumerate}
    \item\label{it:vanishorderv} 
    $v\mp e_3$ vanishes at $z_0e_3$ exactly to order $k$, and there are constants $C,a,\sigma,s > 0$ such that 
    \begin{equation}\label{eq:v-C2a}
    \big| D^j\big( v - (ar^k\cos k \theta,ar^k\sin k\theta,\pm e_3)\big)\big| \leq C r^{k - j +\sigma}, \quad j = 0,1,2,
    \end{equation}
    whenever $r^2 + (z-z_0)^2 < s^2$;
    \item\label{it:vanishorderrho} $\varrho-1$ vanishes exactly to order $2k$ at $z_0e_3$; also, if we define the sequence of rescalings
    $$ \varrho_s(r,z):=\frac{\varrho(s r,z_0+s z)-1}{a s^{2k}},$$
    then, for every sequence $s_i\to 0^+$, up to a subsequence there is a blow-up $\tilde\varrho$ such that $\varrho_{s_i}\to \tilde \varrho$ in  $W_\loc^{2,q}(\R^3)$ for any $q\in [1,\infty)$. Moreover, $\tilde \rho$ is a non-constant, $k$-axially symmetric, $(2k)$-homogeneous, non-negative solution of
    \begin{equation}\label{eq:tvr-pde}
        \Delta \tilde\varrho = k^2 r^{2k-2}\chi_{\{\tilde\varrho>0\}}
    \quad\text{in }\R^3;
    \end{equation}
    \item\label{it:singfb} 
    there exists a modulus of continuity $\omega$ and a small constant $s_0 > 0$ such that 
$$
\dist_\cH \big(B_s(z_0e_3)\cap \partial\{\varrho > 1\}, \{r = 0\}\cup (Z_{2k-1}\setminus\{z= 0\})\big) < s\omega(s),
$$
for all $s \in (0,s_0)$, where $Z_\ell$ is the nodal set of the zonal harmonics of degree $\ell$. In particular, the free boundary $\partial\{\varrho>1\}$ cannot be represented as a $C^1$-graph near any point of $\cB_0(u)$.
\end{enumerate}
\end{thm}

Before proceeding with the blow-up analysis required to prove Theorem \ref{thm:bp-precise}, let us show how it easily implies Theorem \ref{thm:axial}.

\begin{proof}[Proof of Theorem \ref{thm:axial}]
Let $k\geq 2$, $\lambda > 1$, and consider a $k$-axially symmetric map $g \in C^2(\Ss^2; \Ss^2)$ with $\lvert\deg(g)\rvert= k$. Let $u\in W^{1,2}(\mb B^3;\R^3\setminus\mb B^3)$ be energy minimizing among all $k$-axially symmetric maps in $W^{1,2}_{\lambda g}(\mb B^3;\R^3\setminus\mb B^3)$. Note that, by axial symmetry, the image of the $z$-axis under $u$ is contained in the $z$-axis; on the other hand, the image of the $z$-axis under $u$ is also disconnected, due to the constraint. Hence,  $u$ is necessarily discontinuous, i.e. $\Sigma(u)\neq\emptyset$. 

By Corollary \ref{cor:dim-ax}, $\Sigma(u)$ is a discrete set on the $z$-axis. Fix any $z_*e_3 \in \Sigma(u)$. By Theorem \ref{thm:reg-ax}, we can find a ball $B_*$ centered at $z_*e_3$ such that $B_*\subset\{\varrho = 1\}$, where $\varrho$ is as in \eqref{eq:vr}. Hence, we can consider the largest open and connected set $U\subset\{\varrho = 1\}$ containing $z_*e_3$. Choose $z_0e_3\in \partial U$. Then $z_0e_3\in \partial\{\varrho > 1\}$ by the definition of $U$. By Theorem \ref{thm:reg-ax} again, we can find a small ball $B$ centered at $z_0e_3$ such that $u \in C^{1,1}(B)$. 
Taking the radius of $B$ smaller if necessary, we deduce from Corollary \ref{cor:uniqueness} that $u$ is the unique minimizing constraint map in $W^{1,2}_u(B;\R^3\setminus\mb B^3)$. 

We now verify the claims in Theorem \ref{thm:axial}; note that the boundary condition $\lvert \deg(g)\rvert=k$ guarantees that $u$ is necessarily non-constant. Claim \eqref{it:nontrivFB} follows from the choice of $z_0e_3\in\partial U$, which readily implies $z_0e_3\in \partial({\rm int}\{\varrho = 1\}) \neq\emptyset$. Claim
\eqref{it:vanishorderDu} follows from Theorem \ref{thm:bp-precise}\eqref{it:vanishorderv}-\eqref{it:vanishorderrho}. Finally, claim 
\eqref{it:singFB} follows from Theorem \ref{thm:bp-precise}\eqref{it:singfb}.
\end{proof}

The proof of Theorem \ref{thm:bp-precise} is rather involved and is divided over the next three subsections. 

\subsection{Proof of Theorem \ref{thm:bp-precise}(\ref{it:vanishorderv})}

We begin our analysis with a simple consequence of Lemma \ref{lem:EL-ax}.

\begin{lem}\label{lem:vr-C11-vp-C2a}
We have $\varrho\in C_\loc^{1,1}(\mb B^3\setminus\Sigma(u))$, $v\in C_\loc^{2,\sigma}(\mb B^3\setminus\Sigma(u);\mathbb S^2)$, and $\varphi\in C_\loc^{2,\sigma}(\mb B^3\setminus\Sigma(u))$ for every $\sigma\in(0,1)$. Also, the pair $(\varrho,v)$ solves 
\begin{equation}\label{eq:main-sys-ax-re}
\begin{dcases}
    \Delta\varrho = \varrho |Dv|^2\chi_{\{\varrho>1\}}\\
    \Delta v + \frac{2}{\varrho} D\varrho\cdot Dv = -|Dv|^2 v
\end{dcases}
\quad\text{in }\mb B^3,
\end{equation}
in the strong sense. 
\end{lem}

\begin{proof}
Note that \eqref{eq:main-sys-ax-re} follows directly from \eqref{eq:main-ax}; in fact, the computation works without $k$-axial symmetry. Concerning the regularity statements, let $z_0e_3\in \mb B^3\setminus\Sigma(u)$. By Lemma \ref{lem:EL-ax}, there is a ball $B\subset \mb B^3$ centered at $z_0e_3$ such that $u\in C^{1,1}(B)$. Since $\varrho = |u|\geq 1$ in $\mb B^3$, we also have $\varrho \in C^{1,1}(B)$. We then have $v\in C_\loc^{2,\sigma}(B)$ for every $\sigma\in(0,1)$: this can be deduced from the second equation in \eqref{eq:main-sys-ax-re}, cf.\ \cite[Theorem 2.4]{FKS}. 

Finally, we deal with the regularity of $\varphi$. Writing $v=(v^1,v^2,v^3) \in C_\loc^{2,\sigma}(B;\R^3)$, we see that $\varphi = \arccos v^3$ in $\mb B^3$. Then, since $\arccos$ is continuous in $[-1,1]$, and $v^3 \in C^{2,\sigma}(B)\subset C(B)$, we have $\varphi\in C(B)$. Since $\varphi(\{r = 0\}\cap B) \subseteq \{0,\pi\}$, shrinking the radius of $B$ if necessary, we may assume without loss of generality that $\varphi \in [0,\frac{\pi}{6}]$; the other case is symmetric by replacing $\varphi$ with $\pi-\varphi$. Then we may also write $\varphi = \arcsin (v^1/(\cos k \theta))$ in $B\cap\{|\theta|<\frac \pi 2\} $. As $\arcsin$ is analytic in $(-1,1)$, $|v^1|\leq \sin\varphi \leq 1/2$ in $B$, and $v^1 \in C^{2,\sigma}(B)$, then $\varphi \in C^{2,\sigma}(B\cap \{|\theta|<\frac \pi 2)$ for every $\sigma\in(0,1)$. But $\vp$ does not depend on $\theta$, so this last restriction on $\theta$ can be removed, completing the proof. 
\end{proof}

\begin{proof}[Proof of Theorem \ref{thm:bp-precise}\eqref{it:vanishorderv}]
Assume that $z_0e_3\in\mb B^3\setminus\Sigma(u)$. By Lemma \ref{lem:vr-C11-vp-C2a} we have, for some $s$ small, that $\varrho\in C^{1,1}(B_s(z_0e_3))$ and $v\in C^{2,\sigma}(B_s(z_0e_3))$ for every $\sigma\in(0,1)$. Since $\varrho\geq 1$ in $\mb B^3$ by the constraint,it follows that $(1/\varrho)D\varrho \in C^{0,1}(B_s(z_0e_3))$. Moreover, by Lemma \ref{lem:vr-C11-vp-C2a}, 
$|Dv|^2 \in C^{1,\sigma}(B_s(z_0e_3))$ for every $\sigma\in(0,1)$. 

Let $\tilde v \colon B_s(z_0e_3) \to\R^2$ be the map consisting of the first two components of $v$, i.e., 
\begin{equation}\label{eq:tv}
\tilde v := ( \sin\varphi\cos k\theta,\sin\varphi \sin k\theta). 
\end{equation}
Since the boundary datum $g\colon \mb S^2\to \lambda \mb S^2$ satisfies $\lvert \deg(g) \rvert=k$, we see from Proposition \ref{prop:deg} that $u$ is surjective onto the annulus $\lambda \mb B^3\setminus \mb B^3$. In particular, we deduce that $\sin\vp$ is not identically zero in $B_s(z_0 e_3)$, since $\vp$ satisfies the unique continuation property by \eqref{eq:vp-pde}, and hence $\tilde v$ is also not identically zero in $B_s(z_0 e_3)$.
We observe that each component of $\tilde v$ is a $C^{2,\sigma}$-solution of the second line of \eqref{eq:main-sys-ax-re}, 
where we think of $(1/\varrho)D\varrho$ and $|Dv|^2$ as coefficients of the equation. Therefore, it follows from \cite[Remark 3]{A} that each component of $\tilde v$ vanishes to finite order. Applying \cite[Theorem 1.5]{Han1994}  to each component of $\tilde v$, we see that there is a non-trivial homogeneous harmonic polynomial $P = (P^1,P^2)$ of degree $d$ such that  
\begin{equation}\label{eq:tv-P-C2a}
|D^j (\tilde v - P)| \leq C (r^2 + (z-z_0)^2)^{\frac{d-j + \sigma}{2}},\quad j\in \{0,1,2\},
\end{equation}
whenever $r^2 + (z-z_0)^2 < s^2$; here we take $s$ smaller if necessary, and the constant $C$ may depend on $\sigma$ but not on $(r,z-z_0)$. 

Since $P$ is homogeneous of degree $d$, \eqref{eq:tv-P-C2a} implies that 
\begin{equation}\label{eq:P}
P  = (p\cos k\theta,p\sin k\theta),
\end{equation}
where $p$ is a homogeneous polynomial of degree $d$ in $(r,z)$. Moreover, since $P$ is non-trivial, so is $p$. In addition, being $\varphi$ a polar angle of $v$, it must satisfy $\sin\varphi \geq 0$ in $\mb B^3$. So, \eqref{eq:tv}, \eqref{eq:tv-P-C2a}, and \eqref{eq:P} altogether imply 
$$p \geq 0 \quad \text{in } \{r\geq 0\}.$$ 
Then, the structure of homogeneous harmonic polynomials \cite[\S 2.1]{HL-book} shows that we must have $d = k$, and 
$p (r,z) = ar^k,$
for some $a > 0$; in other words, $P$ is a \textit{zonal harmonic}. Noting that $\tilde v = P = 0$ on $\{r = 0\}$, it follows from \eqref{eq:tv-P-C2a} that 
\begin{equation}\label{eq:tv-P-C2a-re}
|D^j (\tilde v - P)| \leq C r^{k-j + \sigma},\quad j\in \{0,1,2\},
\end{equation}
whenever $r^2 + (z-z_0)^2 < s^2$. 

Choosing $s$ smaller if necessary, we can also assume that $\sin^2\varphi \in [0,\frac{1}{2}]$ in $B_s(z_0e_3)$. Thus, \eqref{eq:tv} and \eqref{eq:tv-P-C2a-re} together yield
\begin{equation}\label{eq:Dvp-sinvp}
(\partial_r\varphi)^2 + (\partial_z\varphi)^2 + \frac{k^2}{r^2} \sin^2\varphi \leq 2|D\tilde v|^2 \leq Cr^{2k-2},
\end{equation}
whenever $r^2 + (z-z_0)^2 < s^2$ and $r>0$. Also, since $v(z_0e_3) = \pm e_3$, and $\sin^2\varphi \in [0,\frac{1}{2}]$ in $B_s(z_0e_3)$, we have
$$
|v^3 \mp 1| \leq C\sin^2\varphi, \quad |Dv^3| \leq \sin\varphi|D\varphi|,\quad \text{and}\quad |D^2v^3| \leq |D\vp|^2 + \sin^2\varphi |D^2\varphi|^2\qquad \text{in $B_s(z_0e_3)$.}
$$
Therefore, \eqref{eq:Dvp-sinvp} along with $\varphi \in C^{2,\sigma}(B_s(z_0e_3))$ yield
\begin{equation}\label{eq:v3-1-C2a}
|D^j (v^3 \mp 1)| \leq Cr^{2k - j}, \quad j =0,1,2,
\end{equation}
whenever $r^2 + (z-z_0)^2 < s^2$. Combining \eqref{eq:tv-P-C2a-re} with \eqref{eq:v3-1-C2a}, we arrive at \eqref{eq:v-C2a}.
\end{proof}

\subsection{Proof of Theorem \ref{thm:bp-precise}(\ref{it:vanishorderrho})}

We begin by noting that, by Theorem \ref{thm:bp-precise}\eqref{it:vanishorderv}, the $(2k-2)$-homogeneous rescaling of the function $|Dv|^2$, which depends only on $(r,z)$ as seen in \eqref{eq:b}, converges to $k^2a^2 r^{2k-2}$, which is again $(2k-2)$-homogeneous. The following corollary of Theorem \ref{thm:bp-precise}\eqref{it:vanishorderv} shows that this convergence occurs at an almost-Lipschitz rate. 

\begin{cor}\label{cor:v-C2a}
    Given $z_0e_3\in\mb B^3\setminus\Sigma(u)$ and $\sigma\in(0,1)$, there exist $C,s_0 > 0$  such that 
    \begin{equation}
        \big|\big( r\partial_r + (z-z_0)\partial_z  - (2k-2)\big)|Dv|^2 \big| \leq Cr^{2k-2} s^\sigma,
    \end{equation}
    whenever $r^2 + (z-z_0)^2 \leq s^2$ and $s\in(0,s_0)$. 
\end{cor}

\begin{proof} 
    By \eqref{eq:v-C2a}, we have
    \begin{equation}\label{eq:Dv-C1a}
    \big| D^j\big( |Dv|^2 - k^2 a^2 r^{2k-2}\big) \big| \leq C r^{2k - j -2 +\sigma}, \quad j = 0,1,
    \end{equation}
    whenever $r^2 + (z-z_0)^2 < s^2$. Since $\big(r\partial_r + (z-z_0)\partial_z\big) r^{2k-2} = (2k-2)r^{2k-2}$, the conclusion follows from \eqref{eq:Dv-C1a} and the triangle inequality. 
\end{proof}

We now observe that, due to \eqref{eq:main-sys-ax-re}, the function $\varrho - 1$, which corresponds to the distance map from $\mb B^3$ to $\Ss^2$, must vanish exactly to order $2k$ at every free boundary point. 

\begin{lem}[Vanishing order of $\varrho$]\label{lem:vr-C2k}
    Given $z_0e_3\in \partial\{\varrho>1\}$, $\varrho - 1$ vanishes at $z_0e_3$ exactly with order $2k$, and there exist $s_0 > 0$ small and $C>1$ large such that
    \begin{equation}\label{eq:vr-C2k}
    \frac{1}{C} s^{2k}\leq \sup_{B_s(z_0e_3)}\big(\varrho - 1 + s|D\varrho|\big) \leq Cs^{2k},
    \end{equation}
    for all $s \in (0,s_0)$.  
\end{lem}

\begin{proof}
    By Lemma \ref{lem:EL-ax}, $u\in W^{1,2}(\mb B^3;\R^3\setminus\mb B^3)$ is a weakly constraint map. Taking a small $s_0$ we have $u\in C^{1,1}(B_{s_0}(z_0e_3))$. The conclusion then follows directly from Proposition \ref{prop:vanishingorders}.
\end{proof}

Given any $z_0\in (-1,1)$, let us define the following Weiss-type \cite{Weiss1999} energy  
\begin{equation}\label{eq:W}
    \begin{aligned}
        W(z_0,s) := \frac{1}{s^{4k+1}}\int_{B_s(z_0e_3)} \big(|D\varrho|^2+ 2(\varrho-1)\varrho|Dv|^2\big)\,dx - \frac{2k}{s^{4k+2}} \int_{\partial B_s(z_0e_3)} (\varrho - 1)^2 \,d\cH^2;
    \end{aligned}
\end{equation}
recall from \eqref{eq:b} that $|Dv|$ is a function of $(r,z)$.  In the scalar obstacle problem, the analogous energy is monotone. Here we have the following similar result: 

\begin{lem}[Weiss-type almost monotonicity formula]\label{lem:weiss}
    Given any $z_0e_3\in\partial\{\varrho>1\}$ and every $\sigma\in(0,1)$, there is $s_0 > 0$ such that $s\mapsto W(z_0,s) + Cs^\sigma$ is non-decreasing in $s \in (0,s_0)$, where $C>1$ may depend on $s_0$ and $\sigma$. In particular, $W(z_0,0^+):=\lim_{s\to 0} W(z_0,s)$ exists. 
\end{lem}

\begin{proof}
By Theorem \ref{thm:reg-ax}, we have $z_0e_3\not\in\Sigma(u)$. Let $s_0$ be a small radius to be determined and note that, from Lemma \ref{lem:vr-C11-vp-C2a}, $\frac{d}{ds}{W}(z_0,s)$ exists for every $s\in (0, s_0)$.

Let $a > 0$ be as in \eqref{eq:v-C2a}. In what follows, we fix $s\in(0,s_0)$. For every $s > 0$ small, consider 
\begin{equation}\label{eq:vrs-Fs}
\varrho_{s} (r,z) := \frac{\varrho(sr,z_0 + sz) - 1}{a^2s^{2k}},\quad F_{s} (r,z) := \frac{|Dv|^2(sr,z_0 + sz) }{a^2s^{2k-2}}. 
\end{equation}
In view of \eqref{eq:main-sys-ax-re}, $\varrho_{s}$ solves  
\begin{equation}\label{eq:vrs-pde}
    \Delta\varrho_{s} = (1 + s^{2k}\varrho_{s})F_{s}\chi_{\{\varrho_{s}>0\}}\quad\text{in }B_{1/s_0},
\end{equation}
in the strong sense. Let us also write $\dot\varrho_{s} = \frac{d}{ds}\varrho_{s}$ and $\dot F_{s} = \frac{d}{ds}F_{s}$, so that
\begin{equation}\label{eq:dvrs-dFs}
\dot\varrho_{s}  = \frac{(r\partial_r + z\partial_z- 2k)\varrho_{s}}{s},\quad \dot F_{s} = \frac{(r\partial_r + z\partial_z- (2k-2))F_{s}}{s}.    
\end{equation}
Note that, thanks to Corollary \ref{cor:v-C2a} and Lemma \ref{lem:vr-C2k},
\begin{equation}\label{eq:dvrs-dFs2}
   \varrho_s + s|\dot\varrho_s| +  F_s + s^{1-\sigma} |\dot F_{s}| \leq \frac{C}{a^2}\quad\text{in }\mb B^3,
\end{equation}
for every $s\in(0,s_0)$. Also, we can rewrite 
\begin{equation}\label{eq:W-re}
\begin{aligned}
\frac{1}{a^4}W(z_0,s) & = \int_{\mb B^3} \big( |D\varrho_s|^2 + 2\varrho_{s}(1 + s^{2k}\varrho_{s})F_s\big) \,dx  - 2k \int_{\Ss^2} \varrho_{s}^2\,d\cH^2. 
\end{aligned}
\end{equation}
Differentiating both sides of \eqref{eq:W-re} in $s$, integrating by part, and using \eqref{eq:vrs-pde} for $\varrho_s$, we compute  
\begin{equation}\label{eq:W'}
\begin{aligned}
\frac{1}{a^4}\frac{d}{ds}W(z_0,s) &= 2\int_{\mb B^3}  \big(D\varrho_s\cdot D\dot\varrho_s + \dot\varrho_{s} (1 + s^{2k}\varrho_{s}) F_{s} \big) \,dx - 4 k \int_{\Ss^2} \varrho_{s} \dot\varrho_{s}\,d\cH^2 \\
&\quad + 2\int_{\mb B^3}  \varrho_{s} \frac{d}{ds}\big( (1 + s^{2k}\varrho_{s})F_{s}\big)\,dx \\
& = 2\int_{\Ss^2} (r\partial_r\varrho_s + z\partial_z\varrho_s - 2k\varrho_s)\dot\varrho_s \,d\cH^2 \\
&\quad + 2 \int_{\mb B^3} \varrho_s \big((1+ s^{2k}\varrho_s)\dot F_s + s^{2k-1}(2k\varrho_s + s\dot\varrho_s)F_s\big)\,dx\\
& \geq  2 \int_{\Ss^2} s\dot\varrho_s^2 \,d\cH^2 - \frac{C}{a^4} s^{\sigma - 1} \geq - \frac{C}{a^4} s^{\sigma - 1},
\end{aligned}
\end{equation}
where, in the second last inequality, we used \eqref{eq:dvrs-dFs} and \eqref{eq:dvrs-dFs2}. By \eqref{eq:W'}, we obtain 
$$
\frac{d}{ds} (W(z_0,s) + Cs^\sigma) \geq 0,
$$
for every $s\in(0,s_0)$, which finishes the proof. 
\end{proof}

\begin{proof}[Proof of Theorem \ref{thm:bp-precise}\eqref{it:vanishorderrho}]
Let $s_i\to 0^+$ be given 
and define $\tilde\varrho_i := \varrho_{s_i}$, and $\tilde F_i := F_{s_i}$, where $\varrho_s$ and $F_s$ are as in \eqref{eq:vrs-Fs}. By Theorem \ref{thm:bp-precise}\eqref{it:vanishorderv}, we have 
\begin{equation}\label{eq:Fs-lim}
    \tilde F_i \to k^2r^{2k-2}\quad\text{in }C_\loc^{1,\sigma}(\R^3).
\end{equation}
In addition, due to Theorem \ref{thm:bp-precise}\eqref{it:vanishorderv} and Lemma \ref{lem:vr-C2k}, a standard compactness argument yields, up to a non-relabeled subsequence,
\begin{equation}\label{eq:tvrj-tvr}
\tilde\varrho_i \to \tilde\varrho \quad\text{in }C_\loc^{1,\sigma}\cap W_\loc^{2,q}(\R^3),
\end{equation}
for every $\sigma\in(0,1)$ and $q\in[1,\infty)$, for some non-trivial $\tilde\varrho\in C_\loc^{1,1}(\R^3)$. Since $\tilde\varrho_i$ is axially symmetric and nonnegative, so is its limit $\tilde\varrho$. The fact that $\tilde \varrho$ satisfies 
 \eqref{eq:tvr-pde} follows from \eqref{eq:Fs-lim} and \eqref{eq:tvrj-tvr}, along with the regularity $\tilde\varrho \in C_\loc^{1,1}(\R^3)$. 

Thus, it remains to prove the $(2k)$-homogeneity of $\tilde\varrho$, and this follows from the almost monotonicity of the Weiss-type energy, as we now detail. By Lemma \ref{lem:weiss}, the limit $W(z_0,0^+)$ exists. However, analogously to \eqref{eq:W-re}, we also have, for each $t > 0$, that
\begin{equation}\label{eq:W-re2}
    \frac{1}{a^4} W(z_0,s_it) = \frac{1}{t^{4k+1}} \int_{B_t}( |D\tilde\varrho_i|^2 + 2\tilde\varrho_i(1+(s_it)^{2k}\tilde\varrho_i)\tilde F_i)\,dx - \frac{2k}{t^{4k+2}}\int_{\partial B_t} \tilde\varrho_i^2\,d\cH^2. 
\end{equation}
Thus, passing to the limit in \eqref{eq:W-re2} along with \eqref{eq:Fs-lim} and \eqref{eq:tvrj-tvr}, we deduce from Lemma \ref{lem:weiss} that 
\begin{equation}\label{eq:tW-c}
    \frac{1}{a^4}W(z_0,0^+) = \lim_{j\to \infty} W(z_0,s_it) = \tilde W(t),
\end{equation}
for all $t > 0$, where the energy $\tilde W$ is defined as 
\begin{equation}\label{eq:tW}
\tilde W(t) := \frac{1}{t^{4k+1}}\int_{B_t} (|D\tilde\varrho|^2 + 2\tilde\varrho k^2r^{2k-2})\,dx - \frac{2k}{t^{4k+2}}\int_{\partial B_t} \tilde\varrho^2\,d\cH^2.
\end{equation}
Using  \eqref{eq:tvr-pde}, one can following essentially verbatim the proof of Lemma \ref{lem:weiss} to show that 
$$
\frac{d}{dt}\tilde W(t) = 2 \int_{\Ss^2} t\dot{\tilde\varrho}_t^2 \,d\cH^2,\qquad \tilde\varrho_t(r,z)=\frac{\tilde\varrho(tr,tz)}{t^{2k}},
$$
cp. \eqref{eq:W'}. Since \eqref{eq:tW-c} implies that $\tilde W(t)$ is constant, we deduce that $\tilde\varrho$ is $(2k)$-homogeneous, as desired.
\end{proof}

\subsection{Proof of Theorem \ref{thm:bp-precise}(\ref{it:singFB})}\label{sec:blow-upFB}
To complete our analysis, it remains to analyze the free boundary for $\tilde \varrho$. According to Theorem \ref{thm:bp-precise}\eqref{it:vanishorderrho},  $\tilde\varrho$ satisfies \eqref{eq:tvr-pde}, and a crucial feature of this equation is that the forcing term on the right-hand side is \textit{independent} of $z$. In particular, $\partial_z \tilde \varrho$ is a $(2k-1)$-homogeneous $k$-axially symmetric harmonic polynomial on each component of $
\{\tilde \varrho>0\}$. In other words, $\partial_z\tilde \varrho$ is a \textit{zonal harmonic} of degree $(2k-1)$ on each connected component of $\{\tilde \varrho >0\}$. 

Recall that the zonal harmonics of degree $\ell$ are constant multiples of
$$x\mapsto |x|^\ell P_\ell(\cos \vp),$$
where $P_\ell$ is the Legendre polynomial of degree $\ell$.  The polynomial $ t\mapsto P_\ell(t)$ has $\ell$ distinct zeros  in $(-1,1)$, arranged symmetrically about $t = 0$, and the nodal set of $P_\ell$ is given by 
$$Z_\ell=\{P_\ell(\cos \vp)=0\}.$$ We are of course specially interested in the case where $\ell=2k-1$ is odd: in this case, the intersection of  $Z_\ell$ with $\mb S^2$ is  the union of $\ell$-latitude circles. We refer the interested reader to \cite[\S 2.1]{HL-book} for further details.

The previous discussion leads us to consider the connected components of $\R^3\backslash Z_\ell$.  Thus, for $-1 \leq t_1 < t_2 \leq 1$,  consider the axially symmetric region 
\begin{equation}
\label{eq:conicalregion}
\cC_{t_1,t_2} := \{t_1 \sqrt{r^2 +z^2}< z < t_2 \sqrt{r^2 + z^2}\},
\end{equation} 
and an axially symmetric $2k$-homogeneous solutions $q$ to 
\begin{equation}\label{eq:q-pde}
\Delta q = k^2 r^{2k-2} \quad \text{in }\cC_{t_1,t_2}.
\end{equation}
By the axial symmetry of $q$,  we see that 
\begin{equation}\label{eq:q}
q(r,z) = (r^2 + z^2)^k y\bigg(\frac{z}{\sqrt{r^2 +z^2}}\bigg), 
\end{equation}
where $y\colon(t_1,t_2)\to\R$ is a solution to 
\begin{equation}\label{eq:y-ode}
(1 - t^2) \ddot y - 2t \dot y + 2k(2k+1) y = k^2 (1-t^2)^{k-1},\quad t_1 <t < t_2. 
\end{equation}
The following result gives a description of the solutions to this ODE.

\begin{lem}\label{lem:q2k}
  Let $y(t)$ be a solution of \eqref{eq:y-ode} and assume that $k\geq 2$.  Then
  \begin{equation}\label{eq:y}
y  = c_1 P_{2k} + c_2 Q_{2k} + p_k,     
\end{equation}
where $P_{2k}$ is the Legendre polynomial of degree $2k$, $Q_{2k}$ is the Legendre function of the second kind of degree $2k$, and $p_k$ is the particular solution to \eqref{eq:y-ode} given by
\begin{equation}\label{eq:q2k}
p_k(t) = \sum_{{i} = 0}^k a_{2{i}}t^{2{i}},
\end{equation}
under the recurrence relation
\begin{equation}\label{eq:ak}
a_0=0, \qquad a_{2{i} + 2} = \frac{-(2k-2{i})(2k + 2{i} + 1)a_{2{i}} + (-1)^{i} k^2 {{k-1}\choose{{i}}}}{(2{i}+2)(2{i}+1)},\quad {i} = 0,1,2,\cdots.
\end{equation}
Moreover $\dot p_k(t) \neq 0$ for every $t\in (-1,1)\setminus\{0\}$ where $p_k(t) = 0$. 
\end{lem}

\begin{proof}
The claims concerning \eqref{eq:y}--\eqref{eq:ak} are elementary and are left to the interested reader to verify, so let us turn to the last claim.
    Suppose, towards a contradiction, that $p_k(t_0) = \dot p_k(t_0) = 0$ for some $0< |t_0| < 1$. In view of \eqref{eq:y-ode} and $t_0\neq 0$, we obtain  
\begin{equation}\label{eq:ddq2k}
    \ddot p_k(t_0) = k^2 (1-t_0^2)^{k-2}.
\end{equation}
    Using the recurrence relation \eqref{eq:ak}, we can compute from \eqref{eq:q2k} (and $k\geq 2$) that 
    \begin{equation}\label{eq:ddq2k-re}
    \begin{aligned}
    k^2 (1-t_0^2)^{k-2} & = \sum_{i=0}^{k-1}(2i+2)(2i+1)a_{2i+2}t_0^{2i} \\
    &= \sum_{i=0}^{k-1}(2k - 2i)(2k+2i+1)a_{2i}t_0^{2i} + k^2 \sum_{i=0}^{k-1} (-1)^i {{k-1}\choose{i}}t_0^{2i},    
    \end{aligned}        
    \end{equation}
    whence the binomial theorem yields 
    \begin{equation}\label{eq:ddq2k-re2}        
    k^2 t_0^2 (1-t_0^2)^{k-2} = \sum_{i=0}^{k-1}(2k-2i)(2k+2i+1)a_{2i}t_0^{2i} = \sum_{i=1}^k (2k-2i)(2k+2i+1)a_{2i}t_0^{2i};
    \end{equation}
    the last identity follows from the choice $a_0 = 0$ and $(2k-2i)(2k+2i+1)a_{2i} = 0$ when $i =k$. 
    However, $p_k(t_0) = \dot p_k(t_0) = 0$ along with \eqref{eq:q2k} and $t_0\neq 0$ yields 
    \begin{equation}\label{eq:q2k-dq2k}
    \sum_{i=1}^k a_{2i}t_0^{2i-2} = \sum_{i=1}^k 2i a_{2i}t_0^{2i-2} = 0.
    \end{equation}
    By \eqref{eq:q2k-dq2k}, we can first divide by $t_0^2 \neq 0$ the leftmost and rightmost sides of \eqref{eq:ddq2k-re2} and then proceed further as 
    \begin{equation}\label{eq:ddq2k-re3}
    \begin{aligned}
        k^2 (1-t_0^2)^{k-2} &= \sum_{i=1}^k(2k-2i)(2k+2i + 1)a_{2i} t_0^{2i - 2}
        = - \sum_{i=1}^k 2i(2i + 1) a_{2i}t_0^{2i-2}\\
        & = - \sum_{i=1}^k 2i(2i - 1) a_{2i}t_0^{2i-2} = - k^2 t_0^2(1-t_0^2)^{k-2},
    \end{aligned}
    \end{equation}
where the derivation of the last identity follows simply from the first line of \eqref{eq:ddq2k-re}. However, \eqref{eq:ddq2k-re3} implies $k^2(1+t_0^2)(1-t_0^2)^{k-2} = 0$, which is incompatible with $|t_0| < 1$. This finishes the proof. 
\end{proof}

Finally, the next lemma gives a description of the possible free boundaries in the blow-up limit.

\begin{lem}\label{lem:branch}
Let $\tilde\varrho$ be any blow-up limit as in Theorem \ref{thm:bp-precise}\eqref{it:vanishorderrho}. Then its free boundary $\partial\{\tilde\varrho>0\}$ is contained either in the $z$-axis, i.e.\ $\{r = 0\}$, or inside $\{0\}\cup (Z_{2k-1}\setminus\{z =0\})$. 
\end{lem}

\begin{proof}
The simplest case is when the harmonic function $\partial_z\tilde\varrho$ vanishes identically in $\R^3$. Then it is easy to compute that the only such non-trivial, $2k$-homogeneous, axially symmetric solution to \eqref{eq:tvr-pde} is 
$$\tilde\varrho (r,z) = \frac{1}{4}k^2 r^{2k}.$$ This is precisely the case where $\{\tilde\varrho = 0\} = \{r = 0\}$. Hence, in what follows, we assume that $\partial_z\tilde\varrho \not\equiv 0$ in $\R^3$.

As we discussed above, by differentiating both sides of \eqref{eq:tvr-pde}, and noting that $\tilde\varrho = |D\tilde\varrho| = 0$ on $\{\tilde\varrho = 0\}$ and that $\tilde\varrho$ is a non-trivial $2k$-homogeneous axially symmetric function on $\R^3$, we observe that $\partial_z\tilde\varrho$ is a zonal harmonic of degree $2k-1$ on each connected component of $\{\tilde\varrho > 0\}$ that vanishes continuously on $\partial\{\tilde\varrho > 0\}$. Since $\partial_z\tilde\varrho \not\equiv 0$, the above observation along with $|D\varrho| = 0$ on $\{\varrho = 0\}$ shows us that $\{\tilde\varrho = 0\} \subseteq \{\partial_z\tilde\varrho = 0\}\subseteq Z_{2k-1}$. Therefore, it only remains for us to show that $\partial\{\tilde\varrho > 0\}\cap \{z = 0\}\setminus\{0\} = \emptyset$.

Suppose by contradiction that $\partial\{\tilde\varrho > 0\}\cap \{z = 0\}\setminus\{0\}\neq\emptyset$. By the homogeneity of $\tilde\varrho$, we must have $\{z = 0\}\subseteq \partial\{\tilde\varrho  >0\}$. Let $U$ be the maximal open connected set in $\R^3$ contained in $\{\tilde\varrho > 0\}$ and such that $\partial U$ contains $\{z = 0\}$. Since $\partial_z\varrho$ is a zonal harmonic of degree $2k-1$ in $U$ that vanishes on $\partial U$, $U$ cannot be the entire half-space above or below $\{z = 0\}$. Thus, by the homogeneity of $\tilde\varrho$, $U$ is enclosed between a cone $\{z = t_0 \sqrt{r^2 + z^2} \}$, for some $0< |t_0| < 1$, and the plane $\{z = 0\}$; i.e.\ we have (following the notation in \eqref{eq:conicalregion})
$$U = \cC_{0,t_0} \quad \text{ or } \quad U=\cC_{t_0,0}.$$ 
Then, it follows from our discussions at the beginning of this subsection that
$$
0<  \tilde\varrho = q \quad\text{in }U, \qquad \tilde\varrho = |D\tilde\varrho| = 0 \quad\text{on }\partial U, 
$$
where $q$ satisfies \eqref{eq:q-pde}--\eqref{eq:ak} with either $(t_1,t_2) = (0,t_0)$ or $(t_1,t_2) = (t_0,0)$. We note that, as $0<|t_0|<1$, the Legendre function $Q_{2k}$ of the second kind is regular in $[-|t_0|, |t_0|]$. By our counter assumption and the $2k$-homogeneity of $\tilde\varrho$, we must have $\tilde\varrho  = |D\tilde\varrho| = 0$ at both $(r,z) = (1,0)$ and $(r,z) = (\sqrt{1-t_0^2},t_0)$. Then, by \eqref{eq:q} and \eqref{eq:y},
\begin{equation}\label{eq:P2k-Q2k-pk}
 \begin{cases}
c_1 P_{2k}(t) + c_2 Q_{2k}(t) + p_k(t) = 0 \\
c_1 \dot P_{2k}(t) + c_2 \dot Q_{2k}(t) + \dot p_k(t) = 0
\end{cases}
\quad\text{at }t = 0, t_0  .
\end{equation}
By \eqref{eq:q2k} and \eqref{eq:ak}, we have $p_k(0) = \dot p_k(0) = 0$. On the other hand, $2k$ being an even number implies that $P_{2k}(0) \neq 0$, $\dot P_{2k}(0) = 0$, $Q_{2k}(0) = 0$, and $\dot Q_{2k}(0) \neq 0$. Using these elementary facts in \eqref{eq:P2k-Q2k-pk} for $t = 0$, we deduce $c_1 = c_2 = 0$. Using this information back inside \eqref{eq:P2k-Q2k-pk}, it follows that $p_k(t_0) = \dot p_k(t_0) = 0$. Recalling that $0 < |t_0| < 1$, we arrive at a contradiction to Lemma \ref{lem:q2k}. This shows that $\partial\{\tilde\varrho > 0\}\cap \{z = 0\} = \{0 \}$, as desired. 
\end{proof}

\begin{proof}[Proof of Theorem \ref{thm:bp-precise}\eqref{it:singfb}]
The conclusion follows immediately from Lemma \ref{lem:branch}.
\end{proof}


\appendix


\section{Topological degree}
\label{ap:degree}

In this appendix, we recall some standard properties of the local topological degree, referring the reader to \cite{FG} for further details.

Given $\Omega\subset \R^n$, $u\in C^0(\overline \Omega;\R^n)$, and a point $y\not \in u(\partial \Omega)$, one can associate to it an integer, which we denote by $\deg(u,\Omega,y)$ and refer to as the \textit{local topological degree} of $u$ at $y$ with respect to $\Omega$.  This integer can be defined as follows. If $u\in C^1(\overline \Omega)$ and $y$ is a regular value\footnote{Recall that $y$ is a \textit{regular value} of $u$ if $Du(x)$ is invertible for all points $x\in u^{-1}(y)$.} of $u$,  then
\begin{equation}
\label{eq:degu}
\deg(u,\Omega,y):= \sum_{x\in u^{-1}(y)} \textup{sgn}(\det D u(x));
\end{equation}
note that, by Sard's theorem, this defines $\deg(u,\Omega,\cdot)$ a.e.\ in $\R^n$. For maps $u$ that are just continuous, or at points $y$ which are not regular values, one defines $\deg(u,\Omega,\cdot)$ by approximation, in such a way that $y\mapsto \deg(u,\Omega,y)$ is continuous and $\deg(u_j,\Omega,y)\to \deg(u,\Omega,y)$ if $u_j\to u \text{ uniformly}.$

\begin{thm}\label{thm:degprops}
The topological degree has the following properties:
\begin{enumerate}
    \item\label{it:norm} normalization: $\deg(\id,\Omega,y)$ is 1 if $y\in \Omega$ and 0 if $y\not \in \overline\Omega$;
    \item \label{it:decomp}
    decomposition: if $\Omega_1\cap \Omega_2=\emptyset$ then $\deg(u,\Omega_1\cup \Omega_2, y) = \deg(u,\Omega_1,y) + \deg(u,\Omega_2, y)$;
    \item\label{it:hominv} homotopy invariance: if $h\in C^0([0,1]\times \overline \Omega;\R^n)$ and $y\not \in h([0,1]\times \partial \Omega)$ then $$\lambda\mapsto \deg(h(\lambda,\cdot),\Omega,y)$$ is constant in $[0,1]$;
    \item\label{it:im} existence of solutions: if $\deg(u,\Omega,y)\neq 0$ then $y\in u(\Omega)$;
    \item\label{it:const} constancy: $\deg(u,\Omega,\cdot)$ is locally constant in $\R^n\setminus u(\partial \Omega)$.
\end{enumerate}
\end{thm}

In fact, properties \eqref{it:norm}--\eqref{it:hominv} characterize the topological degree uniquely.

It is also useful to have a notion of topological degree for maps defined on \textit{manifolds without boundary}. If $N,\tilde N\subset \R^m$ are connected oriented compact manifold without boundary, then we can define the degree of a map $g\in C^1(N;\tilde N)$ exactly as in \eqref{eq:degu}, namely through
\begin{equation}
\label{eq:degg}
\deg(g) := \sum_{x\in u^{-1}(y)} \textup{sgn}(\det d g_x),
\end{equation}
where $y$ is a regular value of $g$ and we see $dg_x\colon T_xN\to T_y \tilde N$ as an invertible linear map.  This definition is independent of $y$ and it coincides with the usual algebraic-topological definition of the degree: in algebraic topology,  one says that the homomorphism induced by $g$ between the top homology groups of $N$ and $\tilde N$ corresponds to multiplication by $\deg(g)$. In this paper, we are content with the differential topology definition \eqref{eq:degg}, and we will only apply it in the case where $N=\tilde N=\partial \Omega$ is the boundary of a smooth domain in $\R^n$.


\section{The frequency function for elliptic systems}
\label{app:freq}

In this appendix, we recall some standard properties of the Almgren frequency for solutions of elliptic systems and some of its consequences. In this paper we are often concerned with weak solutions $f\in W^{1,2}(\Omega;\R^m)$ of regular elliptic systems of the type
\begin{equation}
    \label{eq:ellsystem}
    \Delta f^i = a_{ij}^\alpha D_\alpha f^j \quad\text{in }\Omega, 
\end{equation}
for each $1\leq i\leq m$.
Indeed, according to Theorem \ref{thm:e-reg}, any solution of \eqref{eq:main-sys} is a solution of such a system near regular points.

In Section \ref{sec:Le} we recalled the definition of the Almgren frequency function, but other variants of this function are also possible. Following \cite{Garofalo1987,NV2}, we define the generalized frequency through
$$\overline N(f,x_0,r):=\frac{r\int_{B_r(x_0)} \left[|Df|^2+(f-f(x_0))\cdot \Delta f \right]dx}{\int_{\partial B_r(x_0)} |f - (f)_{x_0,r}|^2 \haus}.$$
The generalized frequency is essentially comparable to the standard frequency (see the notation in Section \ref{sec:notation1}), since
\begin{equation}
    \label{eq:compfreq}
    \frac{|N(f,x_0,r)-\overline N(f,x_0,r)|}{N(f,x_0,r)}\leq c r,
\end{equation}
where $c$ depends only on $n$ and $\|a_{ij}^\alpha\|_{L^\infty(\Omega)}$, cf.\ \cite[Proposition 4.5]{NV2}. However, the generalized frequency has the following advantage:

\begin{prop}[Almost monotonicity formula]\label{prop:amf}
Let $f\in W^{1,2}(B_2;\R^m)$ be a non-constant weak solution to \eqref{eq:ellsystem}.
   There are constants $c > 1$ and $r_0 \in (0,1)$, depending only on $n$ and $\|a_{ij}^\alpha\|_{L^\infty(B_2)}$, such that $r\mapsto e^{c r} N(f,x_0,r)$ is non-decreasing for all $r\in (0,r_0)$ and for each $x_0\in B_1$. 
\end{prop}

Proposition \ref{prop:amf} is proved in the \textit{scalar} case in \cite{Garofalo1987,NV2}, i.e.\ when $m=1$; however, as pointed out in \cite{HL}, the same calculations go through for general systems of the type \eqref{eq:ellsystem}.

As a consequence of Proposition \ref{prop:amf}, 
the limit $N(f,x_0,0^+):=\lim_{r\to 0} N(f,x_0,r)$ exists. Moreover, by the existence and uniqueness results for tangent maps proved in \cite{Han1994}, this quantity coincides with the \textit{vanishing order} $d$ of $f-f(x_0)$ at $x_0$, i.e.\ $d$ is the integer characterized by
   $$\limsup_{x\to x_0} \frac{|f(x)-f(x_0)|}{|x-x_0|^d}<\infty, \qquad 
   \limsup_{x\to x_0} \frac{|f(x)-f(x_0)|}{|x-x_0|^{d+1}}=\infty.$$
   In fact, near $x_0$, we may write
   $$f(x)-f(x_0)= P_d(x-x_0)+R(x),$$
   for some non-zero $d$-homogeneous vector-valued harmonic polynomial $P_d$ and a remainder term $R$ which vanishes to higher order, and satisfies the estimates
$$|D^j R(x)|\leq c|x-x_0|^{d - j +\e},\quad j\in\{0,1\},
$$
for every $\e \in (0,1)$, where $c$ here may depend on $\e$. 

In a standard way, Proposition \ref{prop:amf} and \eqref{eq:compfreq} also yield the following result, cf.\ \cite[Lemma 2.2]{HL}:

\begin{thm}\label{thm:bddfreq}
    Let $f \in W^{1,2}(B_2;\R^m)$ be a weak solution of \eqref{eq:ellsystem}. There is a constant $c>1$, depending only on $n$ and $\Lambda$, such that 
    $$N(f,x_0,r)\leq c N(f,0,2) \quad \text{for all } r\in (0,1),$$
    for every $x_0\in B_1$. 
\end{thm}

The frequency controls not only the vanishing order, but also doubling properties of $f$. Here we state and prove a variant of the doubling property of $f$, which is used in Section \ref{sec:Le}. 

\begin{lem}\label{lem:freq-monot}
Let $f \in W^{1,2}(B_2;\R^m)$ be a weak solution of \eqref{eq:ellsystem}.  For every $\lambda > 1$ there exist $\delta,\sigma\in(0,1)$, depending only on $n$, $m$, $\lambda$, and $\|a_{ij}^\alpha\|_{L^\infty(B_2)}$, such that if $N(f,0,2) \leq \lambda$ then 
$$
|\{ |Df|\geq \delta \| Df\|_{L^\infty(3Q)} \} \cap Q| \geq (1-\sigma) |Q|,
$$
for every cube $Q$ with $3Q\subset B_1$.
\end{lem}

\begin{proof}
By Theorem \ref{thm:bddfreq} and a compactness argument we can find a small constant $\delta$, which may depend further on $\lambda$, such that for every cube $Q$ with $3Q\subset  B_1$, 
 \begin{equation}\label{eq:Df-double}
 \| Df \|_{L^\infty(3Q)} \leq \frac{1}{2\delta} \| Df\|_{L^\infty(\frac{1}{2}Q)}. 
 \end{equation}
 The rest of the proof is now classical. By the interior estimates \cite{GT}, $Df^i \in C^{1,\gamma}(2Q)$, for each $\gamma\in(0,1)$, with 
 \begin{equation}\label{eq:Df-Ca}
 [Df^i]_{C^{0,\gamma}(2Q)} \leq c_\gamma \osc_{3Q} f^i \leq 3c_\gamma (\diam Q) \|D f\|_{L^\infty(3Q)},
 \end{equation}
 where $c_\gamma$ depends only on $n$, $m$, $\|a_{ij}^\alpha\|_{L^\infty(B_2)}$, and $\gamma$. On the other hand, by \eqref{eq:Df-double}, there exists a point $y_0\in \frac{1}{2} Q$ such that 
 $$
 |Df(y_0)| \geq 2\delta \|Df \|_{L^\infty(3Q)}. 
 $$
 Hence, by \eqref{eq:Df-Ca} with $\gamma = 1/2$, we can find a small scale $r > 0$, depending only on $c_{1/2}$ and $\delta_\lambda$, such that 
 $$
 |Df| \geq \delta \|Df \|_{L^\infty(3Q)}\quad\text{in } Q_s(y_0),
 $$
 where $s = \frac{1}{\sqrt n}r \diam Q$. Since $y_0 \in \frac{1}{2} Q$, we can choose $r$ small so that  $Q_s(y_0)\subset Q$. Thus, the assertion of the lemma holds with $\sigma = 1- 2^{-n}r^n \in (0,1)$. 
 \end{proof}

 As a result of Lemma \ref{lem:freq-monot}, we obtain an $L^{-\gamma}$-estimate for the energy density of solutions to regular elliptic systems. This can be considered as a simpler version of Theorem \ref{thm:Le}. 

\begin{thm}\label{thm:Le-reg}
    Let $f \in W^{1,2}(B_2;\R^m)$ be a weak solution of \eqref{eq:ellsystem}.  For every $\lambda > 1$ there exist $\gamma > 0$ and $c > 1$, depending only on $n$, $m$, $\lambda$, and $\|a_{ij}^\alpha\|_{L^\infty(B_2)}$, such that if $N(f,0,2) \leq \lambda$, then
$$
\int_{B_1} |Df|^{-\gamma}\,dx \leq c. 
$$
\end{thm}

\begin{proof}
    The proof is essentially the same as that of Theorem \ref{thm:Le}, but  simpler. In fact, as $\Sigma(f) = \emptyset$, the critical scale $r_f$ as in Lemma \ref{lem:crit} becomes irrelevant. So, we only need to use Lemma \ref{lem:freq-monot} instead of Lemma \ref{lem:gen}, and do not need to invoke Lemma \ref{lem:above} in the argument. We shall omit the details. 
\end{proof}


\section{Proof of Lemma \ref{lem:strong-ax}}\label{app:strong-ax}

Here we prove Lemma \ref{lem:strong-ax}, which asserts the strong $W^{1,2}$-convergence of the energy minimizers among $k$-axially symmetric maps. Our proof follows essentially from \cite[Lemma 4.1, Theorem 4.2]{G} (see also \cite{HKL}), which deals with harmonic maps in the axially symmetric setting. The involvement of the radial component $\varrho$ needs extra attention, and we shall present the full detail below. 

\begin{proof}[Proof of Lemma \ref{lem:strong-ax}]
Throughout the proof, $c$ will denote a positive constant independent of $i$, and it may differ at each occurrence. Moreover, we shall use the spherical coordinate systems in the target space, to simplify the exposition; e.g., in such a system, \eqref{eq:axial} can be written as
$$
u(r,\theta,z) = (\varrho(r,z), k\theta,\varphi(r,z)). 
$$
Furthermore, let us write $D' = (\partial_r,\partial_z)$.

Upon a $0$-homogeneous rescaling, we may assume that $\{u_i\}_{i=1}^\infty\subset W^{1,2}(\Omega;\R^3\backslash \mb B^3)$  is a bounded sequence for a concentric ball $\Omega$ strictly containing $\mb B^3$. Then 
there exists a map $u\in W^{1,2}(\Omega;\R^3)$ such that 
\begin{equation}\label{eq:wk-ax}
u_i\rightharpoonup u\quad \text{weakly in $W^{1,2}(\Omega;\R^3)$},
\end{equation}
along a subsequence. By the Rellich compactness theorem, up to non-relabelled subsequences, $u_i\to u$ strongly in $L^2(\Omega;\R^3)$ and $u_i\to u$ a.e.\ in $\Omega$. The almost everywhere convergence ensures that $u$ is $k$-axially symmetric and that $|u|\geq 1$ a.e.\ in $\Omega$. 

By \eqref{eq:wk-ax} and  $\mb B^3\Subset\Omega$, we deduce from the Fubini theorem (upon slightly changing the radius of $\Omega$ through a further 0-homogeneous rescaling) that
\begin{equation}\label{eq:DukDu-L2}
\int_{\mb S^2} \bigg(|{D} u_i|^2 + |{D} u|^2 + \frac{|u_i - u|^2}{\lambda_i^2}\bigg) \,d\cH^2 \leq c,
\end{equation}
for a sequence $\lambda_i\to 0$.  
Let $(\varrho_i,\vp_i)$ be the representation of $u_i$, and let $(\varrho,\vp)$ be that of $u$. By \eqref{eq:DukDu-L2}, there is $\lambda_i'\to 0 $ satisfying 
\begin{equation}\label{eq:lk'}
\frac{1}{2}\lambda_i^{3/4} \leq \lambda_i' \leq 2\lambda_i^{3/4}, 
\end{equation}
such that at $p_i^\pm := (\lambda_i', \pm\sqrt{1 - (\lambda_i')^2})\in \mathbb D$, we have 
\begin{equation}\label{eq:pkpm}
|{D}'(\varrho_i,\vp_i)(p_i^\pm)|^2 + |{D}'(\varrho,\vp)(p_i^\pm)|^2 + \frac{|(\varrho_i-\varrho,\vp_i-\vp)(p_i^\pm)|^2}{\lambda_i^2}  \leq c. 
\end{equation}
On the other hand, due to Lemma \ref{lem:EL-ax}, $\varrho_i$ is weakly subharmonic in $\Omega$ and, since $\mb B^3\Subset\Omega$ and $\{u_i\}$ is bounded in $L^2(\Omega)$,
\begin{equation}\label{eq:vrk}
\varrho_i\leq c\quad\text{a.e.\ in }\mb B^3;
\end{equation}
now since $\varrho_i\to \varrho$ a.e.\ in $\Omega$, \eqref{eq:vrk} yields
\begin{equation}\label{eq:vr}
\varrho \leq c \quad\text{a.e.\ in }\mb B^3. 
\end{equation}

Let $v\in W^{1,2}(\mb B^3;\R^3\backslash \mb B^3)$ be an arbitrary $k$-axially symmetric map with $v=u$ on $\mb S^2$. As in the proof of \cite[Lemma 4.1]{G}, let us define $w_i \colon \mb B^3\to\R^3$ as follows. Consider the annulus $\cA_i := \mb B^3 \setminus (1-\lambda_i)\mb B^3$, and the polar region $\cC_i := \{(r,\theta,z) \in\cA_i: r\leq \lambda_i'\sqrt{r^2+z^2}\}$. Set $\tilde\varrho_i \colon\overline{\cA_i} \to \R$ as the linear interpolation between $\varrho_i|_{\mb S^2}$ and $\varrho|_{\mb S^2}$, i.e., 
\begin{equation}
    \label{eq:defvarrho}
    \tilde\varrho_i(r,z) := \lambda_i^{-1}(1 - |(r,z)|)\varrho\bigg(\frac{(r,z)}{|(r,z)|}\bigg) + \lambda_i^{-1}(|(r,z)| + \lambda_i - 1) \varrho_i \bigg(\frac{(r,z)}{|(r,z)|}\bigg),
\end{equation}
and define $\tilde\vp_i \colon \overline\cA_i\to \R$ in the same way. Clearly, $\tilde\rho_i\geq 1$ and $0\leq\tilde\vp_i\leq \pi$ in $\overline{\cA_i}$. Next, we set $\tilde w_i \colon \cC_i \to \R^3$ as the $0$-homogeneous extension of $(\tilde\rho_i,k\theta,\tilde\vp_i)|_{\partial\cC_i}$ with respect to the center $(0,0,\pm(1-\frac{\lambda_i}{2}))$. Finally, define $w_i\colon \overline {\mb B^3}\to \R^3$ by
\begin{equation}\label{eq:wk}
w_i(r,\theta,z)  := \begin{cases}
    v((1-\lambda_i)^{-1}r,k\theta,(1-\lambda_i)^{-1}z) & \text{in } \mb B^3\setminus\cA_i, \\
    (\tilde\varrho_i(r,z),k\theta,\tilde\vp_i(r,z)) & \text{in }\cA_i\setminus\cC_i, \\
    \tilde w_i (r,\theta,z) & \text{in }\cC_i,\\
    u_i(r,\theta,z) & \text{in } \mb S^2.
\end{cases}
\end{equation}
We refer the reader to Figure \ref{fig:extension} for a depiction of this construction.
\begin{figure}
    \centering
    \includegraphics[scale=0.25]{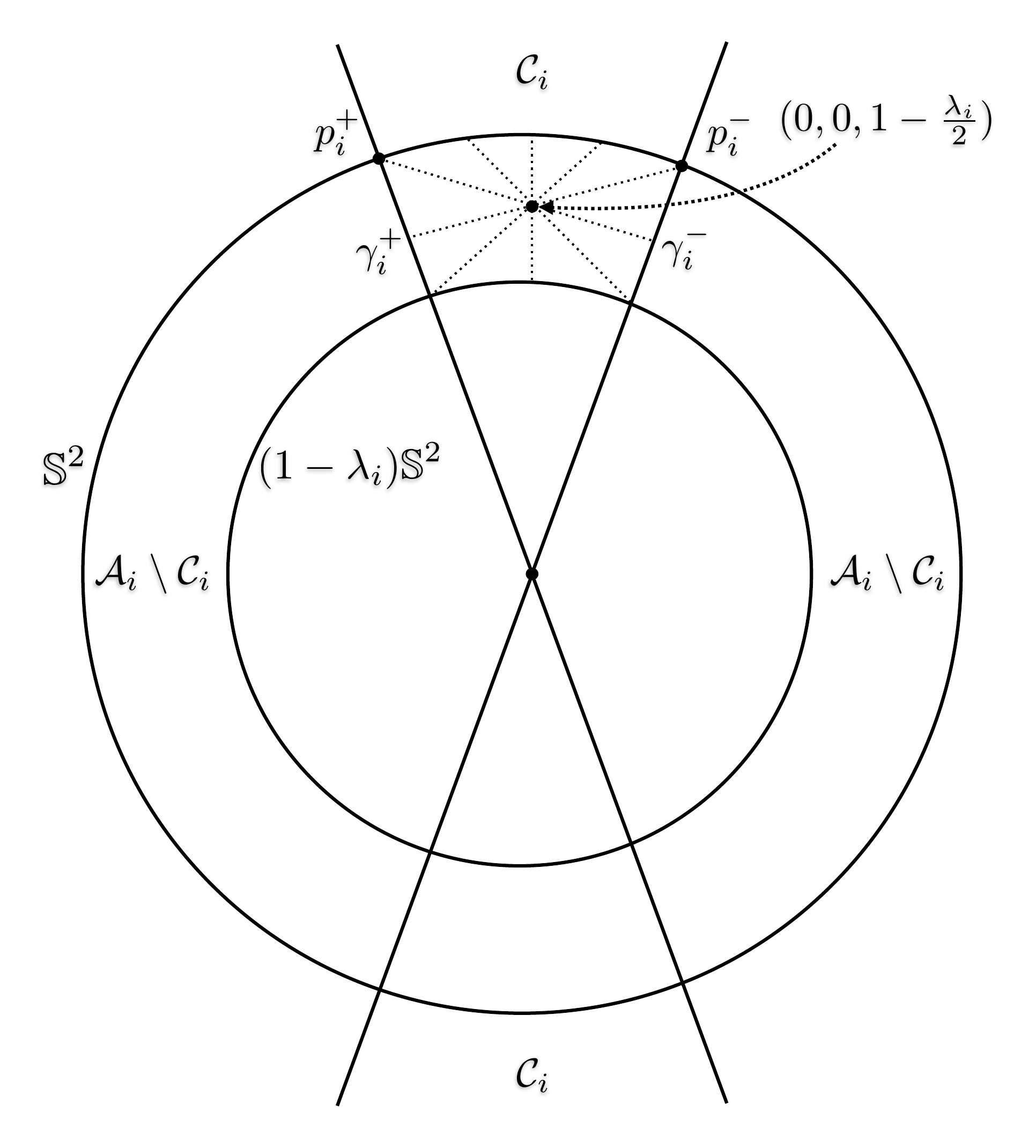}
    \caption{A cross section in the $(r,z)$-plane of the construction in Lemma \ref{lem:strong-ax}.}
    \label{fig:extension}
\end{figure}

Note that $w_i \in W^{1,2}(\mb B^3;\R^3\backslash \mb B^3)$. We claim that 
\begin{equation}\label{eq:Dwk-L2}
\int_{\cA_i} |{D} w_i|^2 \,dx \leq c \lambda_i^{\frac{1}{4}}. 
\end{equation}
Once the claim is proved, then we may use the minimality of $u_i$ over the class of $k$-axially symmetric maps in $W^{1,2}(\mb B^3;\R^3\backslash\mb B^3)$ to deduce that 
$$
\int_{\mb B^3} |{D} u_i|^2\,dx \leq \int_{\mb B^3} |{D} w_i|^2\,dx = (1 + o(1))\int_{\mb B^3} |{D} v|^2\,dx + \int_{\cA_i} |{D} w_i|^2\,dx = \int_{\mb B^3} |{D} v|^2\,dx + o(1). 
$$
However, by \eqref{eq:wk-ax} and the lower semicontinuity of the energy, we also have
$$
\int_{\mb B^3}|{D} u|^2\,dx \leq \liminf_{i\to\infty} \int_{\mb B^3}|{D} u_i|^2\,dx.
$$
The last two inequalities show that $u$ is energy minimizing among $k$-axially symmetric maps, since $v$ was an arbitrary admissible $k$-axially symmetric map with $v=u$ on $\mb S^2$. In fact, taking $v=u$,  the last two inequalities show that
$$
 \int_{\mb B^3} |{D} u|^2\,dx=\lim_{i\to\infty} \int_{\mb B^3} |{D} u_i|^2\,dx,
$$
and so by the weak convergence $u_i \rightharpoonup u$ in $W^{1,2}(\mb B^3;\R^3)$, we deduce that $u_i \to u$ strongly in $W^{1,2}(\mb B^3;\R^3)$, as desired. 

Therefore, we are only left with justifying the claim in \eqref{eq:Dwk-L2}. Denote by $\cD_i$ the domain in the $(r,z)$-plane for which $\cA_i$ is obtained by the revolution of $\cD_i$ about the $z$-axis. Following the first three general (in)equalities at \cite[Page 286]{G}, we obtain 
\begin{equation}\label{eq:Dtvrk-pt}
|{D}'\tilde\varrho_i|^2(r,z) \leq c \bigg[|{D}'\varrho_i|^2 + |{D}'\varrho|^2 + \frac{|\varrho_i - \varrho|^2}{\lambda_i^2} \bigg]\bigg(\frac{(r,z)}{|(r,z)|}\bigg),
\end{equation}
for $(r,z)\in\cD_i$, whence 
\begin{equation}\label{eq:Dtvrk-L2}
   \begin{aligned}
   \iint_{\cD_i} |{D}'\tilde\varrho_i|^2 r\,dr\,dz &\leq c\lambda_i \int_{\partial \mathbb D} \bigg(|{D}'\tilde\varrho_i|^2 + |{D}'\varrho|^2 + \frac{|\varrho_i -\varrho|^2}{\lambda_i^2}\bigg)r\,d\cH^1 \\
   &\leq c\lambda_i\int_{\partial B} \bigg(|{D} u_i|^2 + |{D} u|^2 + \frac{|u_i - u|^2}{\lambda_i^2}\bigg)\,d\cH^2 \leq c\lambda_i. 
   \end{aligned}
\end{equation}
Analogously, we can also compute that 
\begin{equation}\label{eq:Dtvpk-L2}
\iint_{\cD_i} |{D}'\tilde\varphi_i|^2 r\,dr\,dz \leq c\lambda_i. 
\end{equation}
On the other hand, since $(r,z) \in \cD_i$ implies $r > \lambda_i'|(r,z)|\geq \lambda_i'(1-\lambda_i)$, it follows from \eqref{eq:lk'} that 
\begin{equation}\label{eq:sintvpk-L2}
\iint_{\cD_i} \frac{\sin^2\tilde\varphi_i}{r}\,dr\,dz \leq c\lambda_i^{-\frac{3}{4}} |\cD_i| \leq c \lambda_i^{\frac{1}{4}}. 
\end{equation}
However, since \eqref{eq:vrk} and \eqref{eq:vr} yield
\begin{equation}\label{eq:tvrk}
\esssup_{\cA_i}\tilde\varrho_i \leq c,
\end{equation}
combining \eqref{eq:Dtvrk-L2}, \eqref{eq:Dtvpk-L2}, \eqref{eq:sintvpk-L2}, and \eqref{eq:tvrk} altogether gives   
\begin{equation}\label{eq:DwkAkCk-L2}
 \int_{\cA_i\setminus\cC_i} |{D} w_i|^2\,dx = \iint_{\cD_i} \bigg( |{D}'\tilde\varrho_i|^2 + \tilde\varrho_i^2 \bigg( |{D}'\tilde\varphi_i|^2 + \frac{k^2}{r^2}\sin^2\tilde\varphi_i \bigg) \bigg) r\,dr\,dz \leq c\lambda_i^{\frac{1}{4}}.
\end{equation}

The $L^2$-estimate of $|{D} w_i|$ over $\cC_i$ is similar to the second paragraph of \cite[Page 287]{G}. Note that the lateral part of $\partial\cC_i$ is the revolution of the union of the line segments $\gamma_i^\pm := \{t p_i^\pm: 1-\lambda_i \leq t\leq 1\}$. Hence, by the definition \eqref{eq:defvarrho} of $\tilde \varrho_i$,  \eqref{eq:Dtvrk-pt} and \eqref{eq:pkpm}, we get
\begin{equation}\label{eq:Dtvrk-Ck}
|{D}'\tilde\varrho_i|^2(r,z) \leq c\bigg[|{D}'\varrho_i|^2 + |{D}'\varrho|^2 + \frac{1}{\lambda_i^2}|\varrho_i-\varrho|^2\bigg](p_i^\pm) \leq c, 
\end{equation}
$\cH^1$-a.e.\ for $(r,z)\in\gamma_i^\pm$. A similar argument for $\tilde\vp_i$, using also \eqref{eq:lk'}, gives
\begin{equation}\label{eq:Dtvpk-Ck}
|{D}'\tilde\varphi_i|(r,z) + \frac{k^2}{r^2}\sin^2\tilde\varphi_i (r,z)\leq c + c\lambda_i^{-\frac{3}{2}},
\end{equation}
$\cH^1$-a.e.\ for $(r,z)\in\gamma_i^\pm$. Recalling from \eqref{eq:wk} the construction of $w_i|_{\partial\cC_i}$, it follows from \eqref{eq:DukDu-L2}, \eqref{eq:lk'}, \eqref{eq:tvrk}, \eqref{eq:Dtvrk-Ck}, and \eqref{eq:Dtvpk-Ck} that 
\begin{equation}\label{eq:DwkCk-L2}
\int_{\cC_i} |{D} w_i|^2 \,dx \leq |\cC_i| ( c + c \lambda_i^{-\frac{3}{2}} ) \leq c\lambda_i'\lambda_i ( 1 +  \lambda_i^{-\frac{3}{2}} ) \leq c\lambda_i^{\frac{1}{4}}.
\end{equation}
Consequently, the claim in \eqref{eq:Dwk-L2} follows from \eqref{eq:DwkAkCk-L2} and \eqref{eq:DwkCk-L2}, which completes the proof. 
\end{proof}


\section{Uniqueness of continuous weakly constraint maps}\label{app:unique}

The purpose of this appendix is to establish the uniqueness of weakly constraint maps into the exterior of a ball, provided that the images lie in a small neighborhood of a point on the boundary of the ball. Our analysis follows the classical argument of J\"ager--Kaul \cite{JK} for harmonic maps, who consider very general target manifolds without boundary. Here we give a short, self-contained proof for weakly constraint maps, which exploits the symmetry of the ball, although we believe that the same uniqueness result should hold for general target constraints. 

To be precise, let $\Omega$ be a bounded domain in $\R^n$, let $\mb B^m\subset \R^m$ be the open unit open ball centered at the origin, and let $p$ be an arbitrary point on $\mb S^{m-1}$. We say that a map $u\in W^{1,2}(\Omega; \R^m\backslash\mb B^m)$ is \textit{a weakly constraint map} if it verifies 
\begin{equation}\label{eq:main-sys-ball}
\Delta u = - |Du|^2u\chi_{\{|u| = 1\}}\quad\text{in }\Omega
\end{equation}
in the weak sense. For sufficiently regular maps, this system is equivalent to saying that $u$ is a solution of the variational inequality for the Dirichlet energy \cite{CM,D}.

\begin{prop}\label{prop:unique}
    Let $u_i \in W^{1,2}\cap C(\Omega;\R^m\backslash\mb B^m)$ be weakly constraint maps, for each $i\in\{1,2\}$, such that $u_i(\overline\Omega)\subset B_r(p)$ for some $r\in(0,1)$, and $p \in \mb S^{m-1} $. If $u_1 - u_2 \in W_0^{1,2}(\Omega;\R^m)$, then $u_1 = u_2$ in $\Omega$. 
\end{prop}

\begin{proof}
Our proof essentially follows from the idea of \cite{JK} on harmonic maps. However, we chose to present the argument since our system \eqref{eq:main-sys-ball} involves a characteristic term, which may behave badly when taking the difference between two solutions. Here, the exact structure of the obstacle being a ball turns out to resolve this issue quite smoothly.

    Let us introduce 
    $$
    \psi:= \frac{1}{2} |u_1 - u_2|^2 \quad\text{and}\quad \psi_i := \frac{1}{2} |u_i - p|^2,
    $$
    and define 
    $$
    w := \frac{\psi}{(1-\psi_1)(1-\psi_2)}.
    $$
    The basic idea of the proof is that, up to a conformal change of the  metric on the domain (which depends on $\psi_1,\psi_2$), the function $\psi$ is subharmonic, and so it satisfies a maximum principle.
    
    As $u_i\in C(\Omega;\R^m\backslash\mb B^m)$, it follows from  Theorem \ref{thm:e-reg} that $u_i\in C_\loc^{1,1}(\Omega;\R^m\backslash\mb B^m)$. Hence, we have $\psi,\psi_1,\psi_2, w\in W^{1,2}\cap C_\loc^{1,1}(\Omega)$. By \eqref{eq:main-sys-ball}, a direct computation yields 
    $$
    \Delta \psi = |D(u_1-u_2)|^2 - (1 - u_1\cdot u_2)\sum_{i=1}^2 |Du_i|^2\chi_{\{|u_i|=1\}}
    $$
    a.e.\ in $\Omega$. Recall that, if $f$ is a Sobolev function, then for every $k\in \N$ we have $D^kf=0$ a.e.\ in $f^{-1}(0)$, thus
    \begin{equation}\label{eq:psi-pde-re}
    \Delta \psi = 0 \quad\text{a.e.\ on }\psi^{-1}(0). 
    \end{equation}
    On the other hand, note that $|D\psi|^2 \leq 2\psi |D(u_1-u_2)|^2$ in $\Omega$. Also, by the constraint, we have $\psi \geq 1 - u_1\cdot u_2$ in $\Omega$. Using these two inequalities in the above equation, we derive that
\begin{equation}\label{eq:psi-pde}
\Delta\psi \geq \frac{|D\psi|^2}{2\psi} - \psi \sum_{i=1}^2 |Du_i|^2\chi_{\{|u_i|=1\}}\quad\text{a.e.\ in }\Omega\setminus\psi^{-1}(0).
\end{equation}
Again by the constraint we have $\psi_i \geq 1 - p\cdot u_i$ in $\Omega$, and so by \eqref{eq:main-sys-ball} we similarly derive that  
\begin{equation}\label{eq:psii-pde}
\begin{aligned}
\Delta\psi_i &= |Du_i|^2 - (1- p\cdot u_i)|Du_i|^2\chi_{\{|u_i|=1\}} \geq (1 - \psi_i) |Du_i|^2 \chi_{\{|u_i|=1\}}\quad\text{a.e.\ in }\Omega. 
\end{aligned}
\end{equation}
One may compare \eqref{eq:psi-pde-re}, \eqref{eq:psi-pde}, and \eqref{eq:psii-pde} with \cite[(27), (29)]{JK}. 

As in the proof of \cite[Theorem A]{JK}, let us take 
$$
\phi := -\sum_{i=1}^2\log(1-\psi_i),
$$
so that $w = e^\phi \psi$ in $\Omega$; note from the assumption $u_i (\overline\Omega)\subset B_r(p)$ for some $r\in(0,1)$ that $\phi\in C(\overline\Omega)\cap C_\loc^{1,1}(\Omega)$ and $0\leq \phi \leq c$ in $\Omega$. Define the self-adjoint elliptic operator $\cL \equiv \ddiv (e^{-\phi}\nabla\cdot )$; the ellipticity is uniform over $\Omega$ as $1\geq e^{-\phi}\geq e^{-c}$ there. 

As in the proof of \cite[Lemma 1]{JK}, we can compute that 
$$
\cL w = \cL (e^\phi \psi) = \Delta \psi + \psi \Delta\phi + D\psi\cdot D\phi\quad\text{a.e.\ in }\Omega.
$$
As above, we have
\begin{equation}\label{eq:Lw1}
\cL w = 0\quad\text{a.e.\ in }\psi^{-1}(0).
\end{equation}
On the other hand, away from $\psi^{-1}(0)$, we have 
\begin{equation}\label{eq:Lw2}
\begin{aligned}
\psi\Delta \phi + D\psi\cdot D\phi &=\sum_{i=1}^2 \bigg(\frac{\psi|D\psi_i|^2}{(1-\psi_i)^2} + \frac{\psi\Delta\psi_i}{1-\psi_i} + \frac{ D\psi\cdot D\psi_i}{1-\psi_i}\bigg)
\geq  - \frac{|D\psi|^2}{2\psi} + \psi \sum_{i=1}^2 \frac{\Delta\psi_i}{1-\psi_i},
\end{aligned}
\end{equation}
where in the last inequality we used Young's inequality
$$\frac{D\psi\cdot D\psi_i}{1-\psi_i}\leq \frac{\psi|D\psi_i|^2}{(1-\psi_i)^2} + \frac{|D\psi|^2}{4\psi}.$$
Thus, we deduce from \eqref{eq:psi-pde}, \eqref{eq:psii-pde}, and \eqref{eq:Lw2} that 
\begin{equation}\label{eq:Lw3}
\cL w\geq  \Delta\psi - \frac{|D\psi|^2}{2\psi} + \psi \sum_{i=1}^2 \frac{\Delta\psi_i}{1-\psi_i} \geq 0 \quad\text{a.e.\ in }\Omega\setminus \psi^{-1}(0).
\end{equation}
Combining \eqref{eq:Lw1} with \eqref{eq:Lw2}, we conclude that 
$$
\cL w \geq 0\quad\text{a.e.\ in }\Omega. 
$$
By the (uniform) ellipticity of $\cL$, the weak maximum principle applies to $w \in W^{1,2}\cap C_\loc^{1,1}(\Omega)$. Finally, if $u_1 - u_2 \in W_0^{1,2}(\Omega;\R^m)$, then $w \in W_0^{1,2}(\Omega)$, so the weak maximum principle implies $w = 0$ in $\Omega$, i.e., $u_1 = u_2$ in $\Omega$ as desired. 
\end{proof}

We also note here the following result, which is the analogue for constraint maps of the results of Hildebrant--Kaul--Widman for harmonic maps \cite{Hildebrandt1977}:

\begin{thm}\label{thm:fuchs}
    Let $u\in W^{1,2}(\Omega;\R^m\backslash\mb B^m)$ be a weakly constraint map.
    If $u(\overline\Omega)\subset B_r(p)$ for some $r\in(0,1)$ and some $p\in \mb S^{m-1}$, then $u$ is continuous.
\end{thm}

\begin{proof}
Since $0<r<1$, $B_r(p)\cap \mb S^{m-1}$ is a graph in the direction of $p$. Thus we can employ \cite[Theorem 1]{F1} to deduce that $u$ is continuous. Note that, although \cite[Theorem 1]{F1} is stated for minimizers, the short proof only relies on the Euler--Lagrange system \eqref{eq:main-sys-ball} and the results of \cite{Hildebrandt1977a}, which hold for solutions of elliptic systems.
\end{proof}

As an immediate corollary of the above two results we obtain that a weakly constraint map which  takes values in a small neighborhood of $p$ is the unique weakly constraint map with those boundary values: 

\begin{cor}\label{cor:uniqueness}
    Let $u \in W^{1,2}(\overline\Omega;\R^m\backslash\mb B^m)$ be a weakly constraint map such that $u(\overline\Omega)\subset B_r(p)$ for some $r\in(0,1)$ and some $p\in \mb S^{m-1}$.  Then $u$ is the unique weakly constraint map with its own boundary data, and so in particular it is minimizing.
\end{cor}

\begin{proof}By Theorem \ref{thm:fuchs}, $u$ is continuous. 
    Let $v\in W^{1,2}_u(\Omega; \R^m\backslash\mb B^m)$ be an arbitrary weakly constraint map. We claim that $v\in B_r(p)$ a.e.\ in $\Omega$.  Indeed,  \eqref{eq:psii-pde} shows that $f:= |v - p|^2$ is a weakly subharmonic function in $\Omega$. However, as $v \in u + W_0^{1,2}(\Omega;\R^m)$ and $u(\overline\Omega)\subset B_r(p)$, we have $f\leq r^2$ on $\partial \Omega$, so the weak maximum principle implies that $f \leq  r^2$ a.e.\ in $\Omega$, i.e.\ $v\in B_r(p)$ a.e.\ in $\Omega$. The claim now follows from Proposition \ref{prop:unique}.
\end{proof}


 \section*{Acknowledgements}
A.G. was supported by Dr.~Max Rössler, the Walter Haefner Foundation, and the ETH Zürich Foundation. Additionally, he is thankful to Uppsala University and KTH for their hospitality during a visit where part of this research was conducted.
S.K. was supported by a grant from the Verg Foundation.
H.S. was supported by the Swedish Research Council (grant no.~2021-03700).

\section*{Declarations}

\noindent {\bf  Data availability statement:} All data needed are contained in the manuscript.

\medskip
\noindent {\bf  Funding and/or Conflicts of interests/Competing interests:} The authors declare that there are no financial, competing or conflict of interests.


\bibliographystyle{abbrv}
\bibliography{library.bib}

\end{document}